\documentclass[12pt]{amsart}

\usepackage{latexsym,amsfonts,amsthm,amsmath,amssymb,amscd}

\newtheorem{theorem}{Theorem}[section]

\newtheorem{corollary}[theorem]{Corollary}

\newtheorem{lemma}[theorem]{Lemma}

\newtheorem{proposition}[theorem]{Proposition}

{\theoremstyle{definition}

       \newtheorem{definition}[theorem]{Definition}
       \newtheorem{remark}[theorem]{Remark}
       \newtheorem{example}[theorem]{Example}
       \newtheorem{parrafo}[theorem]{{\!}}  }

\numberwithin{equation}{theorem}

 \newcommand{\cali}{{\mathcal {I}}}
\newcommand{\calo}{{\mathcal {O}}}

\newcommand{\Fnu}{\boldsymbol{\nu}}

\DeclareMathOperator{\Max}{\underline{Max}}
 \DeclareMathOperator{\ord}{ord}
 \DeclareMathOperator{\Reg}{Reg}
\DeclareMathOperator{\Sing}{Sing}
\DeclareMathOperator{\Spec}{Spec}
\DeclareMathOperator{\word}{w-ord}
\DeclareMathOperator{\vord}{v-ord}
\def\O{{\mathcal{O}}}

\title{ An Introduction to constructive Desingularization.}

\author{O. Villamayor U.}

\date{June 2005}

\begin{document}

\maketitle
\section{Introduction}

\small Resolution of singularities has been an area of intense
research since the late eighties. Particularly in simplification
of the theory, but also in the task of implementations.

In these notes, intended for non-specialist, we present this new
approach to the subject. So here we prove two important Theorems
of algebraic geometry over fields of characteristic zero:

1) Desingularization (or Resolution of singularities).

2) Embedded Principalization or Log-Resolution of ideals.

Both results, stated in Theorems \ref{classical} and
\ref{principalization}, are due to Hironaka. We focus here on the
proof in \cite{EncVil99}, which is more elementary than that of
Hironaka. In fact, it avoids the use of Hilbert Samuel functions,
and of normal flatness.

Theorem \ref{principalization}, of Embedded Principalization,
plays a fundamental role in the study of morphisms, and
particularly on the elimination of base points of linear systems.

 Hironaka's proof of both theorems is {\it existential}; he
proves that every singular variety, over a field of characteristic
zero, can be desingularized. Our proof of the theorems is {\it
constructive}, in the sense that we provide an algorithm to
achieve such desingularization. We refer to \cite{GabSch98} and to
\cite{FP} for two computer implementations. Bodn\'ar-Schicho's
implementation available at
\begin{center}
 http://www.risc.uni-linz.ac.at/projects/basic/adjoints/blowup
\end{center}

 There are
several other proofs of these two theorems, which also provide an
algorithm: \cite{BM97}, \cite{BV2}, \cite{EH}, \cite{Villa89}, and
\cite{WL}.

It is natural to ask why is it interesting to study algorithms of
resolution of singularities. Usually we simply need to know the
formulation of a theorem in order to apply it. But sometimes a
proof of a theorem can be strong enough to be useful as a tool.
This is the purpose of developing algorithms to achieve resolution
of singularities; a theorem with many applications in algebraic
geometry. A very natural application arises, for example, when we
want to classify singularities by the way that they can be
desingularized. To this end it is not enough to know that
singularities can be resolved, it is necessary to have an explicit
manner to resolve them. This is an advantage of a constructive (or
algorithmic) proof over an existential proof.

These notes are written as an introduction to the subject, and
includes the contents of various one weeks courses on the subject
(see also \cite{Villa05}).
Resolution of singularities is based on
a peculiar form of induction.
 In the case of resolution
 of hypersurfaces this form of induction was stated clear and explicitly by Abhyankar, in what is called a Tschirnhausen transformation.

  We will focus on this point in Part I, where we discuss examples of this form of induction,
 with some indication on how it provides inductive invariants. These
 invariants are gathered in our {\em resolution functions}, and we
 prove the two Main Theorems \ref{classical} and
\ref{principalization} by extracting natural properties from
 these functions.
 In Part II we prove results which were motivated through examples in the first Part. In Part III we
 introduce the resolution functions.
 A mild technical aspect appears in Part II, where the behavior of
 derivations and monoidal transformations are discussed. But essentially the first
 three parts are intended to provide a conceptual (non-technical) and self-contained introduction to desingularization.

Technical aspects appear in the last Part IV, where we present the
algorithm in
 full detail. This will allow the reader to understand also other algorithms, and will hopefully
 encourage the search for new ones.

 These notes follow the notation
 in \cite{Villa92} (basic objects, and general basic objects). In that paper we prove that the algorithm of
 desingularization in \cite{Villa89} is equivariant, and that it
 also desingularizes schemes in \'etale topology.
But the algorithm in \cite{Villa89} (and in
 \cite{Villa92}) provide embedded desingularizations which makes
 use of Hironaka's invariants (of Hilbert-Samuel functions); whereas
 in these notes we discuss an algorithm in which such invariants are avoided.
 Hence the outcome
 is, in general, a different embedded desingularization.
 It turns out, however, that both algorithms coincide when it
 comes to the case of embedded desingularization of hypersurfaces.
 For this reasons we refer to the examples in \cite{Villa92}, such
 as the desingularization of the Whitney Umbrella, or for examples that illustrate
 equivariance of the desingularization of embedded hypersurfaces.

The algorithm in these notes is also equivariant, and also extends
to \'etale topology. However we do not study these properties in
these introductional notes, and we refer to \cite{BEV} and
\cite{EncVil97:Tirol} for the study of these and of further
properties of this proof. Among these further properties discussed
in those cited papers, there is a new and remarkable formulation
of embedded desingularization, with a strong algebraic flavor,
obtain in \cite{BV2} (see \ref{anab} in these notes).

We finally refer to the notes of D. Cutkosky \cite{Cutkosky}, H.
Hauser \cite{Hauser}, and K. Matsuki \cite{Matsuki}, for other
introductions to desingularization theorems.

I thank Ana Bravo and Augusto Nobile for suggestions and
improvements on the notes.

\section{First definitions and formulation of Main Theorem.}
 The set of regular points, of a reduced scheme of finite type over a field, is a dense open set.

\begin{definition}\label{buenosair} We say that a birational morphism of reduced irreducible schemes
\begin{equation}\label{lopri}
X  \stackrel{\pi}{\longleftarrow} X'
\end{equation}
is a desingularization of $X$ if:

i) $\pi$ defines an isomorphism over the open set $U=Reg(X)$ of
regular points.

ii) $\pi$ is proper, and $X'$ is regular.

\end{definition}

We will prove the existence of desingularizations, over fields of
characteristic zero, by proving a theorem of {\em embedded
desingularization} in Theorem \ref{classical}. There we view an
irreducible scheme as a closed subscheme in a smooth scheme $W$.

Let $W_1  \stackrel{{\pi}}{\longleftarrow} W_2$ be a proper
birational morphism of smooth schemes of dimension $n$. If a
closed point $x_2\in W_2$ maps to $x_1\in W_1$, there is a linear
transformation of n-dimensional tangent spaces, say $T_{W_2,x_2}
\to T_{W_1,x_1}$. The set of points $x_2\in W_2$ for which
$T_{W_2,x_2} \to T_{W_1,x_1}$ is not an isomorphism defines a
hypersurface $H$ in $W_2$, called the jacobian or exceptional
hypersurface. It turns out that there is an open set $U\subset
W_1$ such that $U \stackrel{\pi}{\leftarrow}\pi^{-1}(U)$ is an
isomorphism, and $\pi^{-1}(U)=W_2- H$.
 Examples of proper
birational morphisms of this kind are the monoidal
transformations, defined by blowing up a closed and smooth
subscheme $Y$ in a smooth scheme $W_1$. In such case
$H=\pi^{-1}(Y)$ is a smooth hypersurface. Let \tiny
\begin{equation}
\label{ressol}
\begin{array}{cccccccc}
W_{0} & \longleftarrow & (W_{1}, E_1=\{H_1\})& \longleftarrow &
(W_{2}, E_2=\{H_1,H_2\})&\cdots  &\longleftarrow &
(W_{r},E_{r}=\{H_1,H_2,..,H_r \})\\ Y & & Y_1 & & Y_2& &  &
\end{array}
\end{equation}
\small
 be a composition of monoidal transformations, where each  $Y_j\subset
W_j$ is closed and smooth, $H_{j}\subset W_{j}$ is the exceptional
hypersurface of $W_{j-1}\leftarrow W_{j}$ (the blow up at
$Y_{j-1}$), and where $\{H_1,H_2,..,H_r \}$ denote the strict
transforms of the $H_i's$ in $W_r$. The composite $W_0 \leftarrow
W_r$ is a proper birational morphism of smooth schemes, and
$H=\cup_{1\leq i\leq r} H_i$ is the exceptional hypersurface.

\begin{theorem}[{\bf Embedded Resolution of
Singularities}] \label{classical}

 Given $W_0$ smooth over a field $k$ of characteristic zero, and $X_0\subset W_0$ closed and reduced,
 there is a sequence (\ref{ressol}) such that
\begin{enumerate}
\item[(i)] $\cup_{i=1}^rH_i$ have normal crossings in $W_r$.

 \item[(ii)]
 $W_0-\Sing (X_0)\simeq W_r\setminus\cup_{i=1}^rH_i$, and hence it induces a square diagram
 \tiny
\begin{equation*}
\begin{array}{ccc}
W_0    & \stackrel{\Pi_r}{\longleftarrow} & W_{r}
\\
 \cup & & \cup \\
  X_0 &\stackrel{\overline{\Pi}_r}{\longleftarrow} & X_r \\
\end{array}
\end{equation*}
\small of proper birational morphisms, where $X_r$ denotes the
strict transform of $X_0$. \item[(iii)]  \( X_{r} \) is regular
and has normal crossings with \( E_{r}=\cup_{i=1}^rH_i \).
\end{enumerate}
\end{theorem}

  In particular
\(\Reg(X_0)\cong \overline{\Pi}^{-1}_r(\Reg(X_0))\subset X_{r}\)
and $ X_0 \stackrel{\overline{\Pi}_r}{\longleftarrow}  X_r$ is a
desingularization (\ref{buenosair}).

\begin{theorem}[{\bf Embedded Principalization of
ideals}] \label{principalization}

Given \( I\subset\O_{W_{0}} \), a non-zero sheaf of ideals, there
is a sequence (\ref{ressol}) such that:
\begin{enumerate}
\item[(i)] The morphism \(W_0\leftarrow W_r\) defines an
isomorphism over \(W_0\setminus V(I)\). \item[(ii)] The sheaf
$I\calo_{W_r}$ is invertible and supported on a divisor with
normal crossings, i.e.,
\begin{equation}
{\mathcal {L}}=I\calo_{W_r}=\cali(H_1)^{c_1}
\cdot\ldots\cdot\cali(H_s)^{c_s},
\end{equation}
where \(E^{\prime}=\{ H_1,H_2,\dots, H_s\}\) are regular
hypersurfaces with normal crossings,  \(c_i\geq 1\) for
\(i=1,\ldots,s\), and $E^{\prime}=E_r$ if $V(I)$ has no components
of codimension 1.
\end{enumerate}
\end{theorem}

\bigskip

\centerline {\bf Part I}

Throughout these notes $W$ will denote a smooth scheme of finite
type over a field $k$ of characteristic zero.
 We first recall here some definitions used in the formulation of the
previous theorems.

\begin{definition}
\label{Defnc} Fix \( y \in W \), and let $\{x_1,\dots,x_d\}$ be a
regular system of parameters (r.s. of p.) in the local regular
ring $\O_{W,y}$.

1) $ Y (\subset W )$, defined by $I(Y)\subset \O_W$, is {\bf
regular at} $y \in Y$, if there is a r. s. of p. such that
$I(Y)_y=<x_1,...,x_s>  \mbox{ in } \O_{W,y}.$

2)A set \( \{H_{1},\ldots,H_{r}\} \) of hypersurfaces in \( W \)
has {\it normal crossings at} $y$ if there is a  r.s. of p. such
that $\cup H_i=V(\langle x_{j_1}\cdot x_{j_2}\cdots
x_{j_s}\rangle)$ locally at $y$, for some $j_i\in \{1,\dots r \}$.

3) A closed subscheme $Y$ has normal crossings with $E$ at $y$, if
there is a r.s. of p. such that, locally at $y$: $$
I(Y)_y=<x_1,...,x_s> \mbox{ and } \cup H_i =V(\langle x_{j_1}\cdot
x_{j_2}\cdots x_{j_s} \rangle ).$$
\end{definition}

 $Y$ is said to be {\it regular} if it is regular at any point;
and \( E=\{H_{1},\ldots,H_{r}\} \) is said to have normal
crossings if the condition holds at any point.

\begin{remark}\label{bsas}  If
\begin{equation*}
\label{presol}
\begin{array}{ccc}
W_{0} & \stackrel{\pi}{\longleftarrow} & W_{1} \supset
H=\pi^{-1}(Y),\\
 Y & &  \\
\end{array}
\end{equation*}
denotes a monoidal transformation with a closed and regular center
$Y (\subset W_{0})$ , then:

 1) $\pi$ is proper and $W_1$
smooth.
\

2) $H=\pi^{-1}(Y)$ is a smooth hypersurface in $W_1$.
\

3) $W_{0}-Y \cong W_1-H$ (i.e. $\pi$ is birational).
\end{remark}

\begin{definition}\label{defoforder}
The {\it order} of an non-zero ideal $ J$ in a local regular ring
$(R,M)$ is the biggest integer $b \geq 0$ such that $J\subset
M^b$.
\end{definition}
\begin{remark}
 Assume that $Y$ in \ref{bsas} is irreducible with generic point $y \in W$, and
let $h \in W_1$ be the generic point of $H$. Note that
 $ \O_{W,y }$ is a local regular ring, and that $ \O_{W_1,h }$ is a discrete valuation ring. Let  $M_y$ denote the maximal ideal of $ \O_{W,y}$.

 Set $W_0 \longleftarrow W_1$ and $H\subset W_1$ as above. Then, for an ideal $J\subset  \O_{W }$, the following are
equivalent:

  \

 a) $J_y  \subset  M_y^b$ (i.e. the order of J at $\calo_{W,y}$ is $ \geq b$)

  \

 b)$ J \O_{W_1 }= I(H)^b\cdot J_1$ for some $J_1$ in $\O_{W_1 }$.

  \

c) $J\O_{W_1 } $ has order $  \geq b$ at $\O_{W_1,h }$.

 \

\end{remark}

\begin{definition}\label{deftransid}
Given a sheaf of ideals $J\subset \O_X$ and a morphism of schemes,
 $ X \leftarrow Y$, the sheaf of ideals $ J \O_{Y }$ is called
the {\em total transform} of $J$ in $Y$. In the previous remark we
considered the total transform by a monoidal transformation, and
we do not assume $b$ to be the order of $J$ at the generic point
of $Y$. When such condition holds, then $b$ is the highest integer
for which an expression  $ J \O_{W_1 }= I(H)^b\cdot J_1$ can be
defined; and in such case $J_1$ is called the {\em proper
transform} of $J$.
\end{definition}

The following result will be used to ensure that $E_r$
 has normal crossings in a sequence of monoidal
transformations (\ref{ressol}).

\begin{proposition}\label{hacer} Let $W$ be smooth over $k$, and let $E=\{H_1,\dots,H_s\}$ be a set of smooth
hypersurfaces with normal crossings. Assume that $Y (\subset W)$
is closed, regular, and has normal crossings with
$E=\{H_1,\dots,H_s\}$, and set the monoidal transformation

\begin{equation*}
\begin{array}{ccc}
(W , E=\{H_1,\dots,H_s\})   & \stackrel{\pi}{\longleftarrow} &
(W_{1}, E_1=\{H'_1,\dots,H'_s,H_{s+1}=\pi^{-1}(Y) \}) \\
 Y & &  \\
\end{array}
\end{equation*}
where $H'_i$ denotes the strict transform of $H_i$. Then $E_1$ has
normal crossings in $W_1$.
\end{proposition}



\section{Examples: Tschirnhausen and a form of induction on resolution problems.}
\label{tchirn}

A variety, or an ideal, is usually presented by equations in a
certain number of variables. A key point in resolution problems is
to argue by induction on the number of variables involved. In
order to illustrate the precise meaning of this form of induction
we first consider the polynomial $f= Z^2+2\cdot X \cdot Z+ X^2+X
\cdot Y^2 \in k[Z,X,Y]$, defining a hypersurface $\mathbb{X}
\subset \mathbb {A}^3_k$, where $k$ denotes here an algebraically
closed field of characteristic zero. We will see that all points
in this hypersurface are of multiplicity at most two.

{\bf Question}: How to describe the closed set of points of
multiplicity 2?, say $\mathcal{F}_2\subset \mathbb{X}$.

Recall first two definitions:

\begin{definition}
Set $p \in \mathbb{X}=V(\langle f \rangle ) \subset
\Spec(k[Z,X,Y])$. We say that the hypersurface $\mathbb{X}$ has
{\it multiplicity} $b$ at $p$, or that $p$ is a {\it b-fold point}
of the hypersurface, if $\langle f \rangle $ has order $b$ at the
local regular ring $k[Z,X,Y]_p $ (\ref{defoforder}). We will
denote by $\mathcal{F}_b$ the set of points in $ \mathbb{X}$ with
multiplicity $b$.

\end{definition}
 There are now two ways in which we can address our question.

Approach 1): Consider the extension of the ideal $J=\langle f
\rangle $, say:
$$J(1)=\langle f, \frac{\partial^{} f}{\partial X},
\frac{\partial^{} f}{\partial Y},\frac{\partial^{} f}{\partial Z}
\rangle  .$$ Clearly $V(J(1))=\mathcal{F}_2$. In fact, by taking
Taylor expansions at any closed point $q$ we conclude that $q \in
V(J(1))$ if and only if the multiplicity of $\mathbb{X}$ at $q$ is
at least 2. Note also that $\mathbb{X}$ has no closed point of
multiplicity higher than $2$ since $\frac{\partial^{2}
f}{\partial^2 Z}$ is a unit. So the hypersurface $ \mathbb{X}$ has
only closed points of multiplicity one and two.

As for the non-closed points of $\mathbb{X}$, recall first that in
a polynomial ring any prime ideal is the intersection of all
maximal ideals containing it. On the other hand the multiplicity
defines an upper-semi-continuous function on the hypersurface. So
the multiplicity at a non-closed point, say $y \in \mathbb{X}$,
coincides with the multiplicity at closed points in an non-empty
open set of the closure $\overline{y}$. This settles our question.

\

\begin{parrafo}\label{app2}Approach 2) (linked to the previous): Set
$Z_1=Z+X$. At $k[Z_1,X,Y] = k[Z,X,Y]$:
\begin{equation}\label{aceituna}
 f= Z_1^2+X \cdot  Y^2.
\end{equation}

2i) Note first that $Z_1 \in J(1)$, and hence $\mathcal{F}_2
\subset \overline{W}$, where $\overline{W}=V(Z_1)$ is a smooth
hypersurface. \

2ii) Set $J^*=\langle X\cdot Y^2 \rangle \subset
\O_{\overline{W}}$. We claim that $\mathcal{F}_2\subset
\overline{W}$ is also defined as the set of points $q \in
\overline{W}$ where the order of $J^*$, at the local regular ring
$\O_{ \overline{W},q}$, is at least 2. \

In fact, if $q \in \Spec( k[Z,X,Y])$ is a point (a prime ideal) of
order 2, then $J(1)\subset q$, so $$ Z_1 \in q \subset
k[Z_1,X,Y].$$ It is clear that among the prime ideals containing
$Z_1$, those where $ Z_1^2+X \cdot Y^2$ has order 2, are those
where $X\cdot Y^2$ has order at least 2. So the claim follows by
setting $ \overline{W}=V(Z_1)$ and $J^*=\langle X \cdot Y^2
\rangle \subset \O_{ \overline{W}}$ as before.
\end{parrafo}

\begin{parrafo}We will see that the answer to our earlier Question, provided in
Approach 2, is better adapted to resolution problems, at least
over fields of characteristic zero.

We started by asking for those points where the ideal $\langle f
\rangle \subset k[Z,X,Y]$ has order at least 2. So we fixed an
ideal $J$ ($J= \langle f \rangle $ in this case), and a positive
integer $b$ ($b=2$ in this case), and we considered the closed set
$\mathcal{F}_2$ of points where this ideal has order 2. We ended
up with a new ideal, $J^*=\langle X\cdot Y^2 \rangle $ in the ring
of functions in $ \overline{W}$, where $$
\overline{W}=(Spec(k[X,Y])=)Spec (k[Z_1,X,Y]/ \langle Z_1 \rangle
)\subset Spec(k[X,Y,Z]), $$ together with an integer $b_1=2$,
describing the same closed set $\mathcal{F}_2$, but involving one
variable less.
\end{parrafo}
\begin{definition}\label{couple} Fix a scheme $W$, smooth over a field of
characteristic zero. A {\em couple} will be an ideal $J \subset
\O_W$ and an integer $b$, and will be denoted by $(J,b)$.

 The {\em set described by the couple} will be the set of points
$\{ x \in W / \nu_x(J) \geq b\}$, where $\nu_x(J) $ denotes the
order of $J$ at the local regular ring $\O_{W,x}$.
\end{definition}

\begin{parrafo}\label{lip} The set described by the couple $(J=\langle Z_1^2+X\cdot Y^2
\rangle ,2)$ in $\mathbb{A}^3_k$ is included in a smooth
hypersurface $\overline{W}=V(Z_1)$. The dimension of
$\overline{W}$ is of course one less then that of $W$. This
inclusion is called the {\em local inductive principal}. Note that
this closed set is also defined by the couple $(J^*,2)$ ($J^*=
\langle X\cdot Y^2 \rangle \subset \O_{\overline{W}}$).
\end{parrafo}
\begin{example} \label{example2}
The fact that $J^* \subset \O_{\overline{W}}$ is principal just a
coincidence of the previous example. Let now $\mathbb{Y} \subset
\mathbb{A}^3_k$ be the hypersurface defined by $g= Z^3+X\cdot Y^2
\cdot Z+X^5 \in k[Z,X,Y]$. Define $$J(2)=\langle g, \
\frac{\partial^{} g}{\partial x_{i}} , \ \frac{\partial^{2}
g}{\partial x_{i}\partial x_j}
       \ / \mbox{ where } x_1=X, x_2=Y, x_3=Z \rangle$$
so $V(J(2))=F_3$ is the set of points of multiplicity at least 3.
The pattern of this equation is $$Z^3+a_2\cdot Z+a_3 \mbox{ with }
a_2,a_3 \mbox{ in } k[X,Y].$$ One can check that $Z \in J(2)$, and
that $\mathbb{Y}$ has at most points of multiplicity 3 since
$\frac{\partial^{3} g}{\partial^3 Z}$ is a unit.

We can argue as in Approach 2 to show that if $q \in \Spec(
k[Z,X,Y])$ is a point (a prime ideal) of multiplicity 3, then
$J(2) \subset q$. So $$ Z \in q \subset k[Z,X,Y],$$ and among all
prime ideals $q$ containing $Z$, the polynomial $Z^3+X \cdot Y^2
\cdot Z+X^5$ has order 3 at $k[Z,X,Y]_ q$ if and only if $X \cdot
Y^2$ has order at least 2, and $X^5$ has order at least 3. In fact
$Z$ has order one at $k[Z,X,Y]_ q$, and $Z, X,$ and $Y$ are
independent variables.

 Set now $\overline{W}=V(Z)$, $\overline{a_2}=\overline{X \cdot Y^2}$,
$\overline{a_3}=\overline{X^5}$ (the class of $a_2$ and $a_3$ in
$\O_{\overline{W}}$), and note that $$ F_3=\{x\in \overline{W}/
\nu_x(\overline{a_2}) \geq 2; \nu_x(\overline{a_3}) \geq 3\};$$
where $\nu_x(\overline{a_i}) $ denotes the order of
$\overline{a_i}$ at the local regular ring $\O_{\overline{W},x}$.

Set
\begin{equation} \label{f1}
(J^*, 6), \mbox{ where } J^*=\langle
(\overline{a_2})^3,(\overline{a_3})^2 \rangle \subset
\O_{\overline{W}}.
\end{equation}

Finally check that $F_3 \subset \overline{W}$ (local inductive
principal (\ref{lip})), and note that we use this fact to show
that the closed set $F_3$ is also defined by the couple $(J^*,
6)$.
\end{example}

\begin{remark} {\bf Transformations of couples and stability of inductive principal.}\label{trcpp}

Let $\mathbb{Y} \subset \mathbb{A}^3_k$ be the hypersurface
defined by $g= Z^3+X\cdot Y^2 \cdot Z+X^5 \in k[Z,X,Y]$, as in
Example \ref{example2}. The origin $ \overline{0}\in
\mathbb{A}^3_k$ is clearly a point of the closed set defined by
$(J,3)$. We now define:
\begin{equation}\label{eq1}
  \mathbb{A}^3_k \longleftarrow {W_1}
\end{equation}
as the blowup at $\overline{0}$. Let $\overline{W}_1$ be the
strict transform of $\overline{W}$, $\mathbb{Y}_1$ the strict
transform of $\mathbb{Y}$, and $H$ the exceptional hypersurface.
By restriction of the morphism to the subschemes we obtain
\begin{equation}\label{eq11}
  \overline{W} \longleftarrow \overline{W}_1,
\end{equation}
 which is also the monoidal transformation at the point
$\overline{0} \in \overline{W}$, with exceptional hypersurface
$\overline{H}=H \cap \overline{W}_1$.

Note that there is a well defined factorization of the form
\begin{equation}\label{eq2}
  J\O_{W_1}=I(H)^3 \cdot J_1
\end{equation}
 for an ideal $J_1 \subset \O_{W_1}$, defined in terms of (\ref{eq1}); and a factorization
 \begin{equation}\label{eq3}
  J^* \O_{\overline{W}_1}=I(\overline{H})^6 \cdot J^*_1
\end{equation}
 for  $J^*_1 \subset \O_{\overline{W}_1}$, defined in terms of (\ref{eq11}). These factorizations
hold because $\overline{0}$ is a point of the closed set defined
by $(J,3)$, thus of the closed set in $\overline{W}$ defined
by $(J^*, 6)$.

Since $\overline{0}$ is a point of order 3 of $J$ (a point of
multiplicity 3 of the hypersurface $\mathbb{Y}$), $J_1 \subset
\O_{W_1}$ is the
 ideal defining the strict transform $\mathbb{Y}_1$.

\

{\bf Claim:} The set of 3-fold points of the hypersurface
$\mathbb{Y}_1$, or say the closed set of points defined by $(J_1,
3)$, is included in $\overline{W}_1$ and coincides with the closed
set defined by $(J^*_1, 6)$.

In other words, we claim that the role played by $\overline{W}$
and $(J^*, 6)$ for the hypersurface $\mathbb{Y}$ (the local
inductive principal (\ref{lip})), is now played by
$\overline{W}_1$ and $(J^*_1, 6)$ for the hypersurface
$\mathbb{Y}_1$. We call this the {\em stability} of the local
inductive principal.

\

To check this claim note first that $W$ can be covered by three
charts: $$U_{X}=Spec( k[Z/X,X,Y/X])= \mathbb{A}^3_k$$
$$U_{Y}=Spec( k[Z/Y,X/Y,Y])= \mathbb{A}^3_k$$ $$U_{Z}=Spec(
k[Z,X/Z,Y/Z])= \mathbb{A}^3_k$$

The morphism: $ \mathbb{A}^3 \longleftarrow U_Y=Spec(
k[Z/Y,X/Y,Y])= \mathbb{A}^3_k$, induced by (\ref{eq1}), is defined
by the inclusion $ k[Z,X,Y] \to  k[Z/Y,X/Y,Y]$. \

At this chart $I(H)=\langle Y \rangle $, the factorization in
(\ref{eq2}) is $$g= Z^3+X \cdot Y^2\cdot Z+X^5= Y^3\cdot
((Z/Y)^3+(X/Y)\cdot Y \cdot (Z/Y)+(X/Y)^5 \cdot Y^2),$$ and
$I(\overline{W_1} \cap U_Y )= \langle Z/Y \rangle $.

Note that $g_1=(Z/Y)^3+(X/Y)\cdot Y \cdot (Z/Y)+(X/Y)^5 \cdot Y^2
\in k[Z/Y,X/Y,Y]$ has the same general pattern as $g$, namely:
$(Z/Y)^3+ b_2 \cdot (Z/Y)+b_3$, with $b_2, b_3$ in $k[X/Y,Y]$. So
the same argument applied to $g$ asserts that:

\

{\bf 1)} The set of 3-fold points of $ \mathbb{Y}_1 \cap U_Y$ is
included in $V( \langle Z/Y) \rangle )$, or say in
$$\overline{W}_1 \cap U_Y= \Spec( k[Z/Y,X/Y,Y]/\langle Z/Y \rangle
)=\Spec( k[X/Y,Y]).$$

\

{\bf 2)}  The set of 3-fold points $\mathbb{Y}_1$ in $U_Y$ is the
closed set in $\overline{W}_1 \cap U_Y$ defined by $( \mathcal{A}
,6)$, where $$\mathcal{A} =\langle
(\overline{b_2})^3,(\overline{b_3})^2 \rangle \subset k[X/Y,Y].$$
 \

We are finally ready to address the main property of our form of
induction in the number of variables, namely the compatibility of
induction with transformations. To this end note that $$
\overline{W} \longleftarrow  \overline{W}_1 \cap U_Y$$ is defined
by $ k[X,Y] \to  k[X/Y,Y]$, and the transform of the couple
$(J^*,6)$ in (\ref{f1}), defined in (\ref{eq3}), is such that $$
J^*_1 \O_{(\overline{W}_1 \cap U_Y)}= \mathcal{A}.$$

A similar argument applies for $ \mathbb{A}^3 \longleftarrow U_X$.
To study our claim for $ \mathbb{A}^3 \longleftarrow W_1$ it
suffices to check at the charts $U_X, U_Y$. In fact, $U_X \cup
U_Y$ cover all of $W_1$ except for one point (the origin at
$U_Z=\mathbb{A}^3$), which is not a point of $\mathbb{Y}_1$. So
$U_Z$ can be ignored for our purpose.
\end{remark}

\begin{parrafo}\label{slip}
 {\bf Summarizing: Stability of inductive principal.} Our previous discussion showed that the set of 3-fold points of $\mathbb{Y} \subset
\mathbb{A}^3 $ (defined by $g= Z^3+X\cdot Y^2 \cdot Z+X^5 \in
k[Z,X,Y]$) is included in a smooth hypersurface $\overline{W}$
(defined by $Z \in k[Z,X,Y]$)(\ref{lip}). From this fact we
conclude that the set is also defined by $(J^*,6)$, where $J^*$ is
an ideal in the surface $\overline{W}$. The property that links
$\overline{W}$ with 3-fold points of $\mathbb{Y}$ goes beyond this
fact. A transformation at a 3-fold point of $\mathbb{Y}$ defines a
strict transform $\mathbb{Y}_1$. It also induces a transformation
$\overline{W} \longleftarrow \overline{W}_1$, together with a
transformation of $(J^*,6)$, say $(J^*_1,6)$. $\overline{W}_1$ is
the strict transform of $\overline{W}$, and the property is that
the set of three fold points of $\mathbb{Y}_1$ is included in
$\overline{W}_1$. This is what we call the stability of the
inductive principal. Furthermore, $(J^*_1,6)$ defines the closed
set of 3-fold points of $\mathbb{Y}_1$. In particular, if $J^*_1$
would not have points of order $6$ (which is not the case in our
example), then $\mathbb{Y}_1$ would not have 3-fold points. Here
we have analyzed this stability for one quadratic transformation,
but it turns out that the same argument holds for any sequence of
monoidal transformations: Defining a sequence of transformations,
say
\begin{equation}\label{}
\begin{array}{cccccccc}
 & \mathbb{A}^3  & \stackrel{\pi_1}{\longleftarrow} & W_{1} & \stackrel{\pi_2}{\longleftarrow}
 &\ldots& \stackrel{\pi_r}{\longleftarrow}&W_k,\\ & \mathbb{Y} & & \mathbb{Y}_1 &&&&\mathbb{Y}_k\\
\end{array}
\end{equation}
where each $\pi_{i+1}$ is a  blow-up at a closed and smooth
centers included in the 3-fold points of $ \mathbb{Y}_i$, the
strict transform of $ \mathbb{Y}_{i-1}$, is {\em equivalent} to
the definition of a sequence of transformations
\begin{equation}\label{hayquever}
\begin{array}{cccccccc}
 & \overline{W}  & \stackrel{\pi_1}{\longleftarrow} & \overline{W}_1 & \stackrel{\pi_2}{\longleftarrow}
 &\ldots& \stackrel{\pi_r}{\longleftarrow}&\overline{W}_k.\\ & (J^*,6) & & (J_1^*,6) &&&&(J_k^*,6)\\
\end{array}
\end{equation}
where each $J^*_i\subset \O_{ \overline{W}_1}$, and $(J^*_i,6)$ is
defined in terms of $(J^*_{i-1},6)$ as in (\ref{eq3}). Moreover,
each $\overline{W}_i $ is a smooth hypersurface in $W_i$, and the
closed set defined by $(J^*_i,6)$ in the hypersurface
$\overline{W}_i $ is the set of 3-fold points of $\mathbb{Y}_i$.
In particular, if the second sequence is defined with the property
that $J^*_k $ has no points of order $6$ in $\overline{W}_k$, then
the hypersurface $\mathbb{Y}_k$ has at most points of multiplicity
2.

This is induction on the dimension of the ambient space, where the
lowering of the highest order of an ideal in a smooth scheme of
dimension 3 is {\em equivalent} to a related problem in a smooth
scheme of dimension 2. This property of the smooth hypersurface
$\overline{W}$ will be discussed in Section \ref{scirop}.
\end{parrafo}

\begin{parrafo}{\bf Tschirnhausen.}\label{tschirn1} Set $f=Z^b+a_1Z^{b-1}+\cdots +a_b  \in k[Z,X_1,..,X_n]$, with $a_i \in
k[X_1,..,X_n]$ for $i=1,\dots , b$. If the characteristic of $k$
is zero set $Z_1=Z+\frac{1}{b}a_1$. Check that
$k[Z,X_1,\dots,X_n]=k[Z_1,X_1,\dots,X_n]$, and that
$f=Z_1^b+c_2Z_1^{b-2}+\cdots+c_b$, with $c_i \in k[X_1,\dots,X_n]$
and $c_1=0$. One can argue as in Example \ref{example2}, to show
that the $b$-fold points of $\mathbb{Y}$ are included in the
hypersurface $\overline{W}=V(Z_1)(\subset \mathbb{A}^{n+1})$(local
inductive principle (\ref{lip})). Furthermore, $\overline{W}$ will
have the stability property discussed above, where the role of
$(J^*,6)$ in Remark \ref{trcpp} (in (\ref{hayquever})) is now
played by $(J^*,b!)$, where
$$J^*=\langle c_i^{\frac{b!}{i}}, i=2,3,\dots,b \rangle \subset
\O_{\overline{W}}.$$
\end{parrafo}


\section{Resolution functions and the main resolution theorems.}

Our proofs of the two main theorems \ref{classical} and
\ref{principalization} will be constructive, as opposed to the
original existential proofs of Hironaka. We introduce here the
notion of resolution algorithm, or resolution functions.
Constructive resolutions will be defined in terms of these
functions, and the main purpose in this Section is to show how
both proofs follow easily from natural properties of these
functions.

\begin{parrafo}\label{comenttrans} In \ref{example2} we study the transform of a hypersurface in
$\mathbb{A}^3$ by a monoidal transformation at a 3-fold point.
Note that (\ref{eq2}) is an example of a {\em proper transform} of
an ideal, as defined in \ref{deftransid}. However the ideal $J^*$
has order 9 at the center of the monoidal transformation, so
$J_1^*$ in  (\ref{eq3}) is not a proper transform. This shows that
our form of induction will lead us to transformations, defined by
expressions of the form $ J \O_{W_1 }= I(H)^b\cdot J_1$, even when
$b$ is not the highest possible integer in such expression.
\end{parrafo}
We have defined {\em couples} as pairs $(J,b)$, where $J \subset
\O_W$ is a non-zero sheaf of ideals, and $b\in N$ is a positive
integer. We introduce now two notions related to couples:

$\bullet$ The {\bf closed set} attached to $(J,b)$:
$$\Sing(J,b)=\{ x\in W / \nu_x(J_x)\geq b \},$$ namely the set of
points in $W$ where $J$ has order at least $b$. This is closed in
$W$ (see \ref{PropOrdDelta}, ii)).

\bigskip

$\bullet${\bf Transformation} of $(J,b)$:

 Let $Y \subset \Sing(J,b)$ be a closed and smooth subscheme, and let
\begin{equation*}
\begin{array}{ccccc}
 & W & \stackrel{\pi}{\longleftarrow} & W_{1} \supset
H=\pi^{-1}(Y)& \\ & Y & &  &\\
\end{array}
\end{equation*}
be the monoidal transformation at $Y$. Since $Y \subset
\Sing(J,b)$, the total transform $J\O_{W_1}$ can be expressed as a
product: $$J\O_{W_1}= I(H)^b J_1 (\subset \O_{W_1})$$ for a
uniquely defined $J_1 $ in $ \O_{W_1}$. The new couple  $(J_1, b)$
is called the {\em transform} of $(J, b)$, and the transformation
is denoted by
\begin{equation}
\begin{array}{ccccc}
 & W & \stackrel{\pi}{\longleftarrow} & W_{1} & \\ & (J,b) & & (J_1,b) &\\
\end{array}
\end{equation}

A sequence of transformations will be denoted as

\begin{equation}\label{sectransp}
\begin{array}{cccccccc}
 & W & \stackrel{\pi_1}{\longleftarrow} & W_{1} & \stackrel{\pi_2}{\longleftarrow}
 &\ldots& \stackrel{\pi_k}{\longleftarrow}&W_k.\\ & (J,b) & & (J_1,b) &&&&(J_k,b)\\
\end{array}
\end{equation}
Note that in such case
\begin{equation}\label{express}
J \O_{W_k }= I(H_1)^{c_1}\cdot I(H_2)^{c_2} \cdots
I(H_k)^{c_k}\cdot J_k
\end{equation}
for suitable exponents $c_2,\dots,c_k$, and $c_1=b$. Furthermore,
all $c_i = b$ if for any index $ i<k$ the center $Y_i $ is not
included in $ \cup_{j \leq i}H_j \subset W_i$ (the exceptional
locus of $W \longleftarrow W_i$).

\bigskip

\begin{example}\label{onexample1}
The ideal $J=< x^2-y^5> \subset k[x,y] $ has a unique 2-fold point
at the origin $(0,0) \in \mathbb{A}^2$. Let $ W=\mathbb{A}^2
\longleftarrow W_1$ be the blow up at the origin. The strict
transform of the curve has a unique 2-fold point, say $q \in W_1$.
Set $W_1 \longleftarrow W_2$ by blowing-up $q$. This defines a
sequence,
\begin{equation}
\begin{array}{cccccc}
 & W & \stackrel{}{\longleftarrow} & W_{1} & \stackrel{}{\longleftarrow}& W_2\\ & (J,b) & & (J_1,b) &&(J_2,b)\\
\end{array}
\end{equation}
 Here $J \O_{W_2 }= I(H_1)^{2}\cdot I(H_2)^{4}\cdot J_2$ provides
 an expression of the total transform of $J$ involving $J_2$.
\end{example}
\begin{remark}\label{onexample11}
The ideal $J_1$ in the previous example is the proper transform of
$J$, and $J_2$ is the proper transform of $J_1$ (Def
\ref{deftransid}). In particular $J_2$ does not vanish along $H_1$
or $H_2$. Recall however that this is not a general fact as
indicated in \ref{comenttrans}. Set now $K=J$, and note the same
sequence as before defines $(K,1) ; \ (K_1,1); \ (K_2,1)$ and $K
\O_{W_2 }= I(H_1)^{1}\cdot I(H_2)^{2}\cdot K_2$.

In this case the ideal $K_2$ does vanish along the exceptional
hypersurface $H_i$, in fact  there is a unique and well defined
expression, say
\begin{equation}\label{ecuabarra}
K_2=I(H_1)^a\cdot I(H_2)^b \cdot \overline{K}_2
\end{equation}
in $\O_{W_2}$, so that $\overline{K}_2$ does not vanish along the
exceptional hypersurfaces. It follows from \ref{onexample1} that
$a=1$, $b=2$ and $\overline{K}_2=J_2$.

\end{remark}

\begin{definition}\label{resolpairs}
 Fix $J \subset \O_W$, $W$ smooth over a field of characteristic zero, and a couple $(J,b)$. A sequence of
transformations as in (\ref{sectransp}) is said to be a {\em
resolution} of $(J,b)$ if:

i) $\Sing(J_k,b)= \emptyset$.

ii) The exceptional locus of $W \longleftarrow W_k$, namely $
\cup_{1 \leq i \leq k} H_i$, is a union of hypersurfaces with
normal crossings.

\end{definition}

\begin{parrafo}\label{normalc1}
We define a {\em pair}, denoted by $(W, E=\{H_1,.., H_r\})$, to be
a set of smooth hypersurfaces $H_1,.., H_r$ with normal crossings
in a smooth scheme $W$. \

Let $W \longleftarrow W_1$ be a monoidal transformation at a
closed an d smooth center $Y$. If $Y$ has normal crossings with
$\cup H_i$, we say that $Y$ is permissible for the pair $(W,E)$,
and that
$$(W,E=\{H_1,.., H_r\})\longleftarrow (W_1,E_1=\{H_1,.., H_r,
H_{r+1}\})$$ is a {\em transformation of pairs} (see  Prop
\ref{hacer}). \

We define a {\em basic object} to be a pair $(W, E=\{H_1,..,
H_r\})$ together with a couple $(J,b)$, with the condition that
$J_x\neq 0(\subset \O_{W,x})$) at any point $x\in W$. We indicate
this basic object by $$ (W,(J,b),E).$$

 If a smooth center $Y$
defines a transformation of the pair $(W,E)$, and in addition $Y
\subset \Sing(J,b)$, then a transform of the couple $(J,b)$ is
defined. In this case we say that $$ (W,(J,b),E) \longleftarrow
(W_1,(J_1,b),E_1)$$ is a {\em transformation} of the basic object.
A sequence of transformations
\begin{equation}\label{transfuno}
(W,(J,b),E) \longleftarrow (W_1,(J_1,b),E_1)\longleftarrow \cdots
\longleftarrow (W_s,(J_s,b),E_s)
\end{equation}
 is a {\em resolution} of the basic object if
$\Sing(J_s,b)=\emptyset.$

In such case
\begin{equation}\label{traboooo}
J\cdot \O_{W_s}=I(H_{r+1})^{c_1}\cdot I(H_{r+2})^{c_2}\cdots
I(H_{r+s})^{c_s
}\cdot J_s
\end{equation}
for some integer $c_i$, where $J_s$ is a sheaf of ideals of order
at most $b-1$, and the $H_j$ have normal crossings.

\begin{definition}
\label{semicon} Let \(X\) be a topological space, and \((T, \geq
)\) a totally ordered set.  A function \(g: X \rightarrow T\) is
said to be {\it upper semi-continuous} if: {\bf i)} \(g\) takes
only finitely many values, and, {\bf ii)} for any $\alpha \in T$
the set $$\{x\in X \ / g(x)\geq \alpha \}$$ is closed in $X$.

Then largest value achieved by \(g\) will be denoted by
\begin{equation*}
\max g.
\end{equation*}
Clearly the set
\begin{equation*}
    \Max g=\{ x\in X : g(x)= \max g
\}
\end{equation*}
    is a closed subset of \(X\).
\end{definition}

\end{parrafo}

\begin{parrafo}\label{resfunct}{\bf Resolution functions.} We now show why {\em constructive} resolutions of basic objects
will lead us to simple proofs of both Main Theorems
\ref{classical} and \ref{principalization}.

In \ref{app2} we defined an upper semi-continuous function, say
$h_3: Spec(k[Z_1,X,Y]\to \mathbb{Z}$, defined by taking order of
the ideal $J=<Z_1^2+X \cdot  Y^2>$. It was shown that $\max
h_3=2$, and that $\Max h_3 (=\mathcal{F}_2) \subset \overline{W}$,
where $\overline{W}=V(Z_1)$ is a smooth hypersurface isomorphic to
$Spec(k[X,Y])$. Furthermore, an ideal $J^*=\langle X\cdot Y^2
\rangle \subset \O_{\overline{W}}$ was attached to $\Max h_3$. We
may take now $h_2: Spec(k[X,Y])\to \mathbb{Z}$, defined by taking
order of the ideal $J^*$, so that $\Max h_2$ is included in a
smooth hypersurface; and ultimately define a function $h_1$ with
values at $\mathbb{Z}$.

In this frame of mind it is conceivable to assign a copy of
$\mathbb{Z}$ for each dimension, namely $\mathbb{Z}\times
\mathbb{Z} \times \mathbb{Z}$, with lexicographic order, and a
function, say $h=(h_3,h_2,h_1)$ with values at this ordered set,
so that $h$ is upper semi-continuous. This is not exactly the way
we will proceed, but we will define a totally ordered set for each
dimension, and then take the product of copies of this set, one
for each dimension.

 We will fix an integer $d$, and define a totally ordered set $(I^d,\geq)$. Moreover, for
 any basic object
\begin{equation*}
B: (W,(J,b),E),
\end{equation*}
 dimension of $W=d$, an upper semi-continuous function $f_B:
\Sing(J,b) \to I^d$ will be defined with the property that $\Max
f_B $ is a smooth subscheme of $\Sing(J,b)$, and a permissible
center for the pair $(W,E)$. Thus, a transformation of the basic
object can be defined with center $\Max f_B $.

In this way a unique sequence (\ref{transfuno}) is defined
inductively, by setting centers $\Max f_{B_i}$. In addition, this
sequence defined by the functions will be a resolution of the
basic object. In fact, for some index $s$ (depending on $B$)
$\Sing(J_s,b)=\emptyset $.

In other words, the set $(I^d,\geq)$ will be fixed, and the
functions on this set defined so as to provide a resolution for
any basic object of dimension d. We now state the properties that
will hold for such sequence:
\end{parrafo}

{\bf Properties:}

{\bf P1)} For each $l$, $\Max f_l$ is closed regular and has
normal crossings with $\cup_{H_i\in E_l} H_i$.

\bigskip

{\bf P2)} For some index $k_0$, depending on the basic object $B$,
$\Sing(J_{k_0},b)=\emptyset$.
\bigskip

If $p \in \Sing(J_k,b)$, and $p \notin \Max f_k$, then $p$ can be
identified with a point in $W_{k+1}$. Furthermore, $p \in
Sing(J_{k+1},b)$, and:

{\bf P3)}  $f_k (p)=f_{k+1} (p).$

\bigskip

Of particular interest will be the case of basic objects with
$b=1$. In such case $\Sing(J_0,1)$ is the underlying topological
space of $V(J_0)$ (the subscheme defined by the sheaf of ideals
$J_0$).

{\bf P4)} There is fixed value $R\in I^d$, and whenever $p\in
\Sing(J_0,1)$ is a point where the subscheme defined by $J_0$ is
smooth, then $f_0(p)=R$ (where $f_0:\Sing(J_0,1) \to I^d$).

\vspace{5.mm}

The definition of  $(I^d,\geq)$, and of the functions $f$, will be
discussed in Part III, and studied exhaustively in Part IV. We now
prove our two Main Theorems \ref{classical} and
\ref{principalization} using the the properties of resolution
functions.

\vspace{5.mm}

\begin{parrafo} {\bf Proof of Theorem \ref{principalization}.}\label{proofprin}
Fix $I \subset \O_W$ as in Theorem \ref{principalization}, and
consider the basic object
\begin{equation}\label{unomas}
(W,(J,1),E=\emptyset),
\end{equation}
 with $J=I$, and the resolution defined by the
resolution functions. Property {\bf P2)} asserts that
$\Sing(J_{k_0},1)=\emptyset$ for some index $k_0$. It follows that
$J_{k_0}=\O_{W_{k_0}}$, namely that
$$ J \O_{W_{k_0}}= I(H_1)^{c_1}\cdot I(H_2)^{c_2}\cdots
I(H_{k_0})^{c_{k_0}}.$$ It is easy to check now that the
conditions of the Theorem are fulfilled for $W \leftarrow
W_{k_0}$.
\end{parrafo}

\begin{parrafo}{\bf Proof of Theorem \ref{classical}.}\label{proofclassical}
Let $J\subset \O_{W_0}$ be the sheaf of ideals defining $X\subset
W_0$ in Theorem \ref{classical}, and consider, as above, the
resolution of the basic object (\ref{unomas}) defined by the
resolution functions. So again $J_{k_0}=\O_{W_{k_0}}$, and hence $
J \O_{W_{k_0}}= I(H_1)^{c_1}\cdot I(H_2)^{c_2}\cdots
I(H_{k_0})^{c_{k_0}}$.

Let $V=W_0-\Sing(X)$ be the complement of the singular locus of
$X$. Note that $V$ is an open set, dense in $W_0$, and $f_0(p)=R$
for any $p \in V\cap \Sing(J,1)$. Here $X=\Sing(J,1)$, and $V\cap
\Sing(J,1)$ is dense  in $\Sing(J,1)$ since $X$ is reduced.
Furthermore, $f_0(p)=R$ for any $p\in V\cap \Sing(J,1)$ ({\bf
P4)}). So $\max f_0 \geq R$.

If $\max f_0=R$, then $\Sing(J,1)=\Max f_0$ and $X$ is smooth in
$W_0$ ({\bf P1)}).

If $\max f_0>R$, then $V$ can be identified with an open set, say
$V_1$, in $W_1$, and $f_1(p)=R$ for any $p \in V_1\cap
\Sing(J_1,1)$ ({\bf P3)}).

If $\max f_1=R$, then the strict transform of $X$ is a union of
components of $\Max f_1$, so the strict transform defines an
embedded desingularization ({\bf P1)}).

If $\max f_1>R$ then $V$ can be identified with an open subset
$V_2$ in $W_2$.

Note that that there must be an index $k$, for some $k<k_0$, so
that $\max f_k=R$. In fact this follows from {\bf P4)}, {\bf P2)},
and the fact that $\Sing(J_{k_0},1)=\emptyset$. Note that $V$ can
be identified with an open set, $V_k\subset W_k$, and that the
strict transform of $X$ in $W_k$ fulfills the conditions of the
Theorem.
\end{parrafo}


\section{ On the notion of strict transforms of ideals.}

\begin{parrafo}\label{lawtrans} The notion of strict transform of embedded schemes appears
in the very formulation of our Main Theorem \ref{classical}. A
subscheme of a given schemes is defined by a sheaf of ideals.
Given a blow-up at the scheme, there is a notion of strict
transform of ideals, corresponding to the notion of strict
transform of embedded schemes.

 A novel aspect of the proof of Theorem \ref{classical} given in
\ref{proofclassical}, as compared to the proof of Hironaka and
from previous constructive proofs (\cite{BM97}, \cite{Villa92}),
is that we do not consider, within this algorithmic procedure, the
notion of strict transform of ideals. In fact, let $J\subset \O_W$
be the sheaf of ideals defining $X\subset W_0$, and let
$$(W_0,(J,1),E_0)\leftarrow (W_1,(J_1,1),E_1)$$ be a transformation
with center $Y\subset \Sing(J,1)$. We show here that, in general,
$J_1$ is not the sheaf of ideals defining the strict transform of
$X$ in $W_1$ (i.e. is not the strict transform of $J$). Let $H
\subset W_1$ denote the exceptional locus of
$$W\leftarrow W_1$$ so that $W-Y=W_1-H$. The strict transform of
$X$ is the {\em smallest subscheme} of $W_1$ containing $ X- Y$,
via the identification $W_1-H= W-Y$. In other words, it is the
closure of $ X- Y$ in $W_1$ by this identification.

Such smallest subscheme is defined by the {\em biggest sheaf of
ideals}, say $ \widetilde{J}_{1}\subset \O_{W_{1}}$, which
coincide with $J$ when restricted to $W_1-H$. We claim that the
biggest sheaf ideal that fulfills this condition is that defined
by the increasing union of colon ideals:

$$  \widetilde{J}_{1}= \cup_k (J\O_{W_1}: I(H)^k).$$

To check this, set $U= \Spec(A)$, an open affine set of $W_1$, so
that the hypersurface $H \cap U$ is defined by an element $ a \in
A$. Let $K$ denote the ideal defined by restriction of
$\widetilde{J}_{1}$ to $U$. The localization $K\cdot A_a$ is also
a restriction of the sheaf of ideals $J$ to $U_a= \Spec(A_a)$.

Note that $ K\cdot A_a \cap A$ is the biggest ideal in $A$
defining $K\cdot A_a$ at $U_a= \Spec(A_a)$. On the other hand $
K\cdot A_a \cap A=\cup_k (K: a^k)$. Since this holds for an affine
covering of $W_1$, it turns out that $ \widetilde{J}_{1}$ is the
biggest sheaf of ideals with the previous condition.

The ideal $K$ (the restriction of $ \widetilde{J}_{1}$ to $U$), is
a finite intersection of $p$-primary ideals, called the p-primary
components. The ideal $ K\cdot A_a \cap A$ is obtained from $K$ by
neglecting, in the previous intersection, those $p$-primary
components corresponding to prime ideals containing the element
$a\in A$ (i.e. with closure of $p$ included in the exceptional
hypersurface $H$).

It is not hard to check that $$J_1\subset  \widetilde{J}_{1},$$ in
fact $J_1=(J\O_{W_1}: I(H)^1)$ according to the definition of
transformation of basic objects.

If $W_1$ arises from blowing up $W=A_k^3$ at the origin, and $J=
<Z, X^2-Y^3>\subset k[X,Y,Z]$, then $V(J_1)\cap H$ is a line,
whereas $V( \widetilde{J}_{1} )$ (the strict transform of the
curve), intersects $H$ at a unique point. So $J_1\neq
  \widetilde{J}_{1}$ in this case.
\end{parrafo}

\begin{parrafo} Resolution of singularities is defined by a proper birational
morphism, defined in a step by step procedure, each step
consisting of a suitably defined monoidal transformation. So given
equations defining the ideal $J$, and a monoidal transformation as
above, Hironaka provides equations defining the strict transform
ideal $\widetilde{J}_{1}$. This turns out being, in general, a
very difficult task. In fact a major part of the proof of Hironaka
is devoted to address this particular point; he introduces the
notions of Hilbert-Samuel functions and of normal flatness with
this purpose. An important conceptual simplification of
constructive desingularization, presented in \ref{proofclassical},
relies on the fact that it provides a proof avoiding all these
notions. In fact, we prove resolution by means of elementary
transformations (defining $J_1$), avoiding the use of the strict
transform ideal $\widetilde{J}_{1}$.
\end{parrafo}

\begin{example} \label{weakvsstrict} The following example illustrates a
situation in which both notions of transformations discussed in
\ref{lawtrans} coincide (i.e. where $J_1= \widetilde{J}_{1}$).

 Let $X \subset W$ be a closed and smooth subscheme of $W$.
Set $J=I(X)$, and note that $\Sing(J,1)=X$, and that the order of
$J$ at $\O_{W,x}$ is one at any $x \in X$. \

Let $W \longleftarrow W_1$ be the monoidal transformation with
center $Y$ which defines a transformation, say: $ (J,1) \ ; \
(J_1,1)$. In other words, assume that $Y \subset \Sing(J,1)$ (so
that $J\O_{W_1}= J_1\cdot I(H)$, where $H \subset W_1$ denotes the
exceptional locus). We claim now the following holds:

{\bf (1)} $ \Sing(J_1,1)$ ($=V(J_1)$) is the strict transform of
$X$.

 {\bf (2)} The subscheme $ X_1 \subset W_1$, defined by $J_1$, is smooth.

 \

 Note that {\bf (2)} follows from {\bf (1)}. In fact the induced
 morphism $X \leftarrow X_1$ is the blowup of $X$ at $Y$, and
 the blowup of a smooth scheme in a smooth subscheme is smooth. To
 prove {\bf 1)} note that at any point $ x\in W$, there is a regular system of
 parameters $\{ x_1,\dots,x_n \}$ such that $J_x=\langle x_1,\dots,x_r \rangle $
 and $I(Y)_x=\langle x_1,\dots,x_s \rangle $ for $r \leq s$. The fiber over
 $x \in W$ can be covered by $ Spec(\O_W[x_1/x_i,x_2,\dots, x_s/x_i,
 x_{s+1},\dots,x_n]$ for $i=1,2,\dots,s$. Finally {\bf (1)} can be
 checked directly at the charts corresponding to indices $ r+1 \leq i \leq s$.

\end{example}

\begin{parrafo}\label{anab}
There is a stronger formulation of embedded desingularization than
that in \ref{classical}, which was proved in \cite{BV2}. That
theorem proves that given $W_0$ smooth over a field $k$ of
characteristic zero, and $X_0\subset W_0$ closed and reduced,
there is a sequence of monoidal transformations, say
\begin{equation*}
\begin{array}{cccccccc}
W_{0} & \longleftarrow & (W_{1}, E_1=\{H_1\})& \longleftarrow &
(W_{2}, E_2=\{H_1,H_2\})&\cdots  &\longleftarrow &
(W_{r},E_{r}=\{H_1,H_2,..,H_r \}),\\ Y & & Y_1 & & Y_2& &  &
\end{array}
\end{equation*}
such that, in addition to the three conditions i), ii), and iii)
in \ref{classical}, it also holds that:

iv) $I(X_0)\O_{W_r}=I(H_1)^{c_1}\cdot I(H_2)^{c_2}\cdots
I(H_r)^{c_r}\cdot I(X_r)$

where $X_r$ denotes the strict transform of $X$.

\medskip

Consider the particular case in which $X$ is an irreducible
subscheme in $W_0=\Spec(k[X_1,\cdots ,X_n])$ defined by a prime
ideal $P$ of height h. In this case the theorem says that at any
point $x\in W_r$ there is a regular system of parameters $\{
Z_1,\cdots Z_n\}$ at $\O_{W_r,x}$, such that:

i)$P\cdot \O_{W_r,x}= <Z_1,\cdots , Z_h>\cdot Z_{j_1}^{a_1}\cdot
Z_{j_2}^{a_2}\cdots Z_{j_s}^{a_s}$ if $x$ is a point of the strict
transform $X_r$, and

ii) $P\cdot \O_{W_r,x}= <Z_{j_1}^{a_1}\cdot Z_{j_2}^{a_2}\cdots
Z_{j_s}^{a_s}>$ (is an ideal spanned by a monomial in these
coordinates) if $x$ is not in $X_r$.

This result does not hold, in general, for desingularizations
which make use of invariants such as Hilbert Samuel functions (
which we avoid in our proof). This algebraic formulation of
embedded desingularization is not a consequence of the theorem of
desingularization as proved by Hironaka.

\end{parrafo}


\centerline {\bf Part II}

In \ref{slip} we discussed a strong link between the set of 3-fold
points of the hypersurface $\mathbb{Y} \subset \mathbb{A}^3 $,
defined by $g= Z^3+X\cdot Y^2 \cdot Z+X^5 \in k[Z,X,Y]$, and the
smooth hypersurface $\overline{W}$ defined by $Z \in k[Z,X,Y]$.
The link showed that the reduction of 3-fold points of
$\mathbb{Y}$, by means of monoidal transformations, was equivalent
to a related problem for a suitable ideal in the smooth subscheme
$\overline{W}$ (see also \ref{tschirn1}).

This is the key for induction in resolution Theorems. In this
second Part we justify the discussion in \ref{slip} (see Example
\ref{sigueee}), and generalize this main property in Section 7 .
In section 6 we study an important preliminary: the behavior of
derivations with monoidal transformations.

\section{Derivations and monoidal transformations on smooth schemes.}

In this Section we study behavior of derivations when applying
monoidal transformations. This will be used in the next Section
\ref{scirop}, where the inductive properties discussed in
\ref{slip} will be clarified.

 Fix $W$
smooth over a field $k$, and $y \in W$ a closed point. Let
$\{x_1,\dots,x_n\}$ be a regular system of parameters at
$\O_{W,y}$.

We define an operator $ \Delta_y$ on ideals in $\O_{W,y}$ by
setting, for $J_y=<f_1,f_2,\dots,f_s>$ in $ \O_{W,y}$: $$
\Delta_y(J_y)= < f_1,f_2,\dots,f_s,  \ \frac{\partial^{}
f_j}{\partial x_{i}}
       \ / 1\leq  i \leq n ; 1 \leq j
\leq s>.$$

 Note that $\Delta_y(\Delta_y(J_y))=< f_1,f_2,\dots,f_s, \ \frac{\partial^{}
f_j}{\partial x_{i}} , \ \frac{\partial^{2} f_j}{\partial
x_{i}\partial x_j}
       \ / 1\leq  i \leq n ; 1 \leq j
\leq s> $. The whole point of restriction to fields of
characteristic zero relies on the following property:

\begin{parrafo}{\bf Characteristic zero.}  If $k$ is a field of characteristic zero and $(b \geq 1)$,
then $J_y$ has order $b$ at $ \O_{W,y}$ iff $ \Delta_y(J_y)$ has
order $b-1$.

\begin{example} Let $\O_{W,y}=k[x_1,x_2,x_3]_{<x_1,x_2,x_3>}$. $$J_y=<x_1^3+x_2^4+x_3^4>\subset
\Delta_y(J_y)=<x_1^2,x_2^3,x_3^3>\subset
\Delta_y^2(J_y)=<x_1,x_2^2,x_3^2> \subset
\Delta_y^3(J_y)=\O_{W,y}$$
 Note that, if $k$ is of characteristic zero, the orders of these ideals drop by
one :  3,2,1,0.
\end{example}
\end{parrafo}

\begin{parrafo} Further properties of the operator $ \Delta_y$ are:

{\bf i)} $J_y \subseteq  \Delta_y(J_y)\subseteq
\Delta_y(\Delta_y(J_y))= \Delta^2_y(J_y) \subseteq
\Delta^3_y(J_y)\subseteq \dots$

{\bf ii)}  $J_y\subset\O_{W,y} $ has order $b (\geq 1)$ iff
$\Delta^{b-1}_y(J_y)$ has order 1.

{\bf iii)} The order of $J_y\subset\O_{W,y} $ is $\geq s$ iff
$\Delta^{s-1}_y(J_y)$ is a proper ideal in $\O_{W,y} $.

\end{parrafo}

\begin{parrafo}
\label{PropOrdDelta} {\bf On the \(\Delta\) operator.} The locally
defined operators $\Delta_y$ can be globalized in the following
sense. Fix $W$ smooth over a field $k$, there is an operator
$\Delta$ on the class of all $\O_W$-ideals , such that: $$
J\subseteq \Delta(J)(\subset \O_W),$$ and at any closed point $y
\in W$:
$$\Delta(J)_y=\Delta_y(J_y).$$

Furthermore, the following properties hold:

\medskip

{\bf i)} $J \subseteq  \Delta(J)\subseteq \Delta^2(J)\subseteq
\dots$ ( hence $ V(J) \supseteq  V(\Delta(J))\supseteq
V(\Delta^2(J))\supset \dots$

\bigskip

{\bf ii)} $V(\Delta^{s-1}(J))=\Sing(J,s)$. In fact
$V(\Delta^{s-1}(J))$ is the closed set of points in $W$ where $J$
has order $\geq s$ (i.e. $(\Delta^{s-1}(J))_y
=\Delta_y^{s-1}(J_y)\subsetneq \O_{W,y} )$ iff the order of $J_y
\O_{W,y}$ is $\geq s$).

\bigskip

{\bf iii)} If $b$ is the biggest order of $J$,
$V(\Delta^{b}(J))=\emptyset $ and  $V(\Delta^{b-1}(J))$ is locally
included in a smooth hypersurface.

 {\bf Proof of iii)} If $b$ is the biggest order of
$J$, $\Delta^{b}(J)=\O_W $ and $\Delta^{b-1}(J)$ has order at most
1. So if $y\in V(\Delta^{b-1}(J))$, $\Delta^{b-1}(J)\O_{W,y}$ has
order 1 at $\O_{W,y}$. If $g \in \Delta^{b-1}(J)$ has order 1 at
$\O_{W,y}$, then: $$ \overline{W}=V(<g>) \supset
V(\Delta^{b-1}(J)),$$ and $\overline{W}$ is a smooth hypersurface
in a neighborhood of $y$.

\end{parrafo}
\begin{example}\label{sigue} Set $W=A_k^3=Spec(k[X,Y,Z])$, $F=Z^3+X Y^2
Z+X^5$, and $J=<F>$, as in \ref{slip}. Then:
$$ \Delta(J)=<3Z^2+XY^2, Y^2Z+5X^4, 2XYZ, F> \subset
\Delta^2(J)=<Z, XY, Y^2, X^3>  \subset \Delta^3(J)=k[X,Y,Z].$$ So,
as indicated in \ref{slip}, the 3-fold points of the hypersurface
$\mathbb{Y} \subset \mathbb{A}^3 $ defined by $V(<J>)$ are
included in smooth hypersurface $\overline{W}=V(<Z>)$.

\end{example}

\begin{parrafo}  We now address the compatibility of $\Delta$ operators with
monoidal transformations. So fix a couple $(J,b)$, and consider a
transformation
\begin{equation}\label{corta}
\begin{array}{ccccc}
 & W & \stackrel{\pi}{\longleftarrow} & W_{1} & \\ & (J,b) & & (J_1,b). &\\
\end{array}
\end{equation}
\end{parrafo}

\begin{lemma}\label{lemm1}
Given $(J,b)$ (\(J\subset \calo_W\)) and a transformation
(\ref{corta}), then:
\

1) If $b \geq 2$, (\ref{corta}) induces a transformation of
$(\Delta(J),b-1)$:

\begin{equation*}
\begin{array}{ccccc}
 & W & \stackrel{\pi}{\longleftarrow} & W_{1} & \\ & (\Delta(J),b-1) & & ((\Delta(J))_1,b-1). &\\
\end{array}
\end{equation*}

2) $(\Delta(J))_1 \subset \Delta(J_1)$.

\end{lemma}

{\em Proof:} Let $Y \subset W$ be the center of the monoidal
transformation, and let $H \subset W_1$ be the exceptional locus.
By assumption $Y \subset \Sing(J,b)$, so $ J\cdot \O_{W_1}=
I(H)^b\cdot J_1$. It follows from \ref{PropOrdDelta},{\bf ii)}
that for general $b$, $ \Sing(J,b) \subset \Sing(\Delta(J), b-1)$.
In particular $Y \subset \Sing(\Delta(J), b-1)$, which proves 1).

In order to prove 2) we first note that if $U \subset W$ is open,
a sheaf of ideals in $W$ induces a sheaf of ideals in $U$, and the
$\Delta $ operators (on $W$ and on $U$) are compatible with
restrictions. On the other hand note that the pull-back of $U$ in
$W_1$, say $U_1$, is an open set, and the induced morphism  $U
\longleftarrow U_1$ fulfills the conditions in 1) for the
restriction of $J$ to $U$.

If we can prove that 2) holds over $U$ (at $U \longleftarrow
U_1$), for all $U$ in an open covering of $W$, then it is clear
that 2) holds. Therefore we may argue locally.

 Let \(\xi\in W\) be a closed point and choose a regular system
of parameters \(\{x_1,\dots x_n\}\subset \calo_{W,\xi}\) so that
the center of the monoidal transformation is locally defined by
\(\langle x_1,\dots,x_s\rangle\). Now consider an affine
neighborhood \( U \) of \( \xi \) such that \( x_{1},\ldots,x_{s}
\) are global sections of \( \mathcal{O}_{U} \), and such that \(
J \) is generated by global sections, say  \( f_1,\dots,f_r \). We
may also assume that $ \left\{\frac{\partial\ }{\partial
x_{1}},\ldots,
        \frac{\partial\ }{\partial x_{n}}\right
\}$ are global derivations, and that $\Delta(J)$ is generated by
the global sections $ \left\{f_{k}\right\}_{k=1}^{r}\cup
\left\{\frac{\partial
      f_{k}}{\partial
      x_{j}}\right\}_{\genfrac{}{}{0pt}{}{k=1,\ldots,r}{j=1,\ldots,n}}$.

By the previous discussion we may assume that \( U=W \). The
scheme \(W_1\) is defined by patching the affine rings
$$A_i=\calo_{W}[x_1/x_i,\dots,x_s/x_i],\mbox{ } \mbox{ }\mbox{ }
\mbox{ }i\in \{1,\ldots,s\},$$ and \( I(H)=\langle x_i\rangle \)
at \(A_i\). For each index $ k \in \{1,\dots,r\}$ there is a
factorization \(f_k=x_i^b {g_i}^{(k)}\), and
$\{{g_i}^{(1)},{g_i}^{(2)},\dots,{g_i}^{(r)}\}$ generate the
restriction of $J_1$ to $\Spec(A_i)$, say $J_1^{(i)}$. In order to
prove 2) we must show that, for each index $ k \in \{1,\dots,r\}$:

\

a) $\frac{f_k
      }{x_i^{b-1}} \in \Delta(J_1^{(i)})$, and

\

b)$\frac{(\frac{\partial
      f_{k}}{\partial
      x_{j}})
      }{x_i^{b-1}} \in \Delta(J_1^{(i)}).$
\

The assertion in a) is clear since $\frac{f_k
      }{x_i^{b-1}}=x_ig_i^{(k)} \in J_1^{(i)}\subset  \Delta(J_1^{(i)})$.
      We now address b). In what follows we fix an index $ k \in \{1,\dots,r\}$ and set
      $f=f_k$. We also fix  an index $ j \in \{1,\dots,n\}$ and set
      $\delta=\frac{\partial
      }{\partial
      x_{j}}$ which is a global derivation on
      $U$.

      Note that
\begin{equation*}
       \delta\left(\frac{x_{j}}{x_{i}}\right)=
\frac{\delta(x_{j})}{x_{i}}-
\frac{x_{j}}{x_{i}}\frac{\delta(x_{i})}{x_{i}},
\end{equation*}
and that $$I(H)\cdot \delta|_{\mbox{Spec}(A_i)}=x_i\cdot \delta:
A_i\to A_i,$$ and hence \(\cali(H)\cdot \delta\) is an invertible
sheaf of derivations on \(W_1\).

Now in  $A_i$ consider the factorization \(f=x_i^b g_i\), so \(
g_i\in J_1^{(i)} \subset A_i\), and \(x_i\cdot\delta \) is a
derivation on  \(A_i\). Finally check that
\begin{equation*}
\frac{\delta(f)}{x_i^{b-1}}=
       \frac{x_i\delta(x_i^b\cdot
g_i)}{x_i^{b}}=
       \frac{x_i\delta(x_i^b)}{x_i^b}g_i+
x_i^b\frac{(x_i\delta)(g_i)}{x_i^b})=
       b\cdot\delta(x_{i})\cdot g_i+
(x_i \delta)(g_i).
\end{equation*}
This already proves b) since the right hand side is in
$\Delta(J_1^{(i)})$. \qed

 Our argument also shows that this equality is stable by any {\em
sequence} of transformations (see \ref{luminy}).

\begin{remark}\label{remk1}
Fix $ K\subset J $ two ideals in $\O_W$, and couples $(J,b)$ and
$(K,b)$. Then clearly:

\

a) $\Sing(J,b) \subset \Sing(K,b)$.

\

b) Any transformation, as in (\ref{corta}), of $(J,b)$, induces
the transformation

\begin{equation*}
\begin{array}{ccccc}
 & W & \stackrel{\pi}{\longleftarrow} & W_{1} & \\ & (K,b) & & (K_1,b) &\\
\end{array}
\end{equation*}
 and $K_1\subset J_1.$

\end{remark}

\begin{parrafo}\label{luminy}  We finally extend the previous result to study the behavior of $\Delta$
operators with an arbitrary sequence of transformations.
\end{parrafo}

\begin{corollary}\label{cor1}
Fix a couple $(J,b)$ (\(J\subset \calo_W\)) and a sequence of
transformations
\begin{equation}\label{larga}
\begin{array}{cccccccc}
 & W & \stackrel{\pi_1}{\longleftarrow} & W_{1} & \stackrel{\pi_2}{\longleftarrow}
 &\ldots& \stackrel{\pi_r}{\longleftarrow}&W_r\\ & (J,b) & & (J_1,b) &&&&(J_k,b).\\
\end{array}
\end{equation} \

1) If $b\geq 2$, then (\ref{larga}) induces a sequence of
transformations
\begin{equation*}
\begin{array}{cccccccc}
 & W & \stackrel{\pi_1}{\longleftarrow} & W_{1} & \stackrel{\pi_2}{\longleftarrow}
 &\ldots& \stackrel{\pi_r}{\longleftarrow}&W_r,\\ & (\Delta(J),b-1) & & ((\Delta(J))_1,b-1) &&&&((\Delta(J))_r,b-1),\\
\end{array}
\end{equation*}
and

2) $(\Delta(J))_r \subset \Delta(J_r)$.

\end{corollary}
\begin{proof} The case when $r=1$ is in \ref{lemm1}. Consider now the case $r=2$,
namely
\begin{equation*}
\begin{array}{cccccc}
 & W & \stackrel{\pi_1}{\longleftarrow} & W_{1} & \stackrel{\pi_2}{\longleftarrow}
 &W_3\\ & (J,b) & & (J_1,b) &&(J_2,b).\\
\end{array}
\end{equation*}
Then \ref{lemm1} asserts that $\pi_1$ defines a transform of $
(\Delta(J),b-1)$, say $((\Delta(J))_1,b-1)$, and that
$(\Delta(J))_1 \subset \Delta(J_1)$. The same result says that
$\pi_2$ defines a transform of $ (\Delta(J_1),b-1)$, say
$((\Delta(J_1))_1,b-1)$, and that $(\Delta(J_1))_1 \subset
\Delta(J_2)$. The statement follows in this case from \ref{remk1}.

The general case $r\geq 2$ follows similarly, by induction.
\end{proof}
\begin{corollary}\label{cor2}
Fix a couple $(J,b)$ (\(J\subset \calo_W\)) and, as before, a
sequence of transformations (\ref{larga}). Assume that $b \geq 2$.
Then, for each index $ 1 \leq j \leq b-1$:
\

1) The sequence (\ref{larga}) induces a sequence of
transformations
 $((\Delta^{(j)}(J)),b-1-(j-1))$, say
\tiny
\begin{equation*}
\begin{array}{cccccccc}
 & W & \stackrel{\pi_1}{\longleftarrow} & W_{1} & \stackrel{\pi_2}{\longleftarrow}
 &\ldots& \stackrel{\pi_r}{\longleftarrow}&W_r\\ & (\Delta^{(j)}(J),b-1-(j-1)) & & ((\Delta^{(j)}(J))_1,b-1-(j-1)) &&&&
 ((\Delta^{(j)}(J))_r,b-1-(j-1))\\
\end{array}
\end{equation*}
\normalsize

 and

2) $(\Delta^{(j)}(J))_r \subset \Delta^{(j)}(J_r)$.

\begin{proof}
Note that for $j=1$, $\Delta^{(j)}=\Delta$ and we obtain the
previous corollary. So we prove now the statement for $j$ assuming
that it holds $j-1$. Set $J^{*}=\Delta^{(j-1)}(J)$ and
$b^*=b-1-(j-2)$. By induction:

i) The sequence of transformations (\ref{larga}) induces
transformations of $(J^{*},b^{*})$, say:

\begin{equation*}
\begin{array}{cccccccc}
 & W & \stackrel{\pi_1}{\longleftarrow} & W_{1} & \stackrel{\pi_2}{\longleftarrow}
 &\ldots& \stackrel{\pi_r}{\longleftarrow}&W_r,\\ & (J^{*},b^{*}) & & (J^{*}_1,b^{*}) &&&&(J^{*}_r,b^{*})\\
\end{array}
\end{equation*}
and

ii) $J^*_r \subset \Delta^{(j-1)}(J_r)$.

\

Applying our previous Corollary \ref{cor1} to i), we get:

\

i') The sequence in i) induces transformations of
$(\Delta(J^{*}),b^{*}-1)$:

\begin{equation*}
\begin{array}{cccccccc}
 & W & \stackrel{\pi_1}{\longleftarrow} & W_{1} & \stackrel{\pi_2}{\longleftarrow}
 &\ldots& \stackrel{\pi_r}{\longleftarrow}&W_r,\\ & (\Delta(J^{*}),b^{*}-1)  & &((\Delta(J^{*}))_1,b^{*}-1) &&&&
 ((\Delta(J^{*}))_r,b^{*}-1)\\
\end{array}
\end{equation*}
and

ii') $(\Delta(J^{*}))_r \subset \Delta(J^*_r) $.

\

Here $\Delta(J^{*})=\Delta^{(j)}(J)$ and i') is statement 1). On
the other hand, applying $\Delta$ to ii) we get $$\Delta(J^*_r)
\subset \Delta^{(j)}(J_r),$$ which together with ii') proves 2).

\end{proof}
\end{corollary}

In the next Section we shall apply Corollary \ref{cor2}, basically
in the case $j=b-1$. The reader might want to look into Example
\ref{sigueee} to get have an overview of this application of the
Corollary.


\section{Simple couples and a form of induction on resolution
problems.}\label{scirop}

\begin{parrafo}\label{simplecouples} The purpose of this Section is the study of
{\em simple couples} $(J,b)$ (\(J\subset \calo_W\)). Examples of
simple couples appear already in Section \ref{tchirn}. They will
play a central role in our inductive arguments (induction on the
dimension of the ambient space). The main results of this Section
are Theorem \ref{theorem72} and Proposition \ref{proptch}, where
the notion of stability of induction discussed in \ref{slip} is
formalized.
\end{parrafo}
\begin{parrafo} Fix $J\subset \calo_W $, assume that $J_x\neq 0 (\subset
\calo_{W,x})$ for any $ x \in W$, and define a function
\begin{equation}
ord_J: W \to \mathbb{N},
\end{equation}
 where $ord_J(x)$ denotes the order of $J_x$ in the local ring $ \calo_{W,x}$.

 Note that $ord_J$ is upper-semi-continuous (\ref{semicon}). In fact, for any positive
 integer $s$: $$\{ x \in W / ord_J(x) \geq
 s\}=V(\Delta^{s-1}(J))\mbox{ (see \ref{PropOrdDelta})}.$$
 \begin{remark}
 The following conditions are equivalent:

 1) $\max-ord_J=b$ (where, as in \ref{semicon}, $\max-ord_J$ denotes the maximum value achieved).

 2) $V(\Delta^{b-1}(J))\neq \emptyset \mbox{ and }
 V(\Delta^{b}(J))=
 \emptyset$.

3) $\max-ord_{\Delta^{b-1}(J)}=1$.
 \end{remark}

 The equivalence follows from the properties of the $\Delta$ operator discussed in \ref{PropOrdDelta}.
\end{parrafo}

 \begin{definition}\label{nomeacuerdo}  We say that $(J,b)$ is a {\em simple
 couple} if the previous conditions hold for $J$ and $b$.
 \end{definition}

The following theorem is a central result in this section.
\begin{theorem}\label{theorem72}
If  $(J,b)$ (\(J\subset \calo_W\)) is a simple couple, and
\begin{equation*}
\begin{array}{ccccccc}
 & W & \stackrel{\pi}{\longleftarrow} & W_{1} & & &\\ & (J,b) & & (J_1,b) &&&\\
\end{array}
\end{equation*}
is a transformation, then either $\Sing(J_1,b)=\emptyset$ or
$(J_1,b)$ is a simple couple.
\end{theorem}
The case $b=1$ will be proved in Proposition \ref{caso1}, and the
case $b\geq 2$ in Proposition \ref{casob}.

We shall first draw attention to the case of simple couples of the
form $(J,1)$.

 \begin{remark}\label{cubrii}
 The following conditions are equivalent:

 1) $\max-ord_J=1$.

 2) $V(J))\neq \emptyset \mbox{ and }
 V(\Delta^{}(J))=
 \emptyset$.

3) There is an open covering $\{ U_{\lambda}\}_{\lambda \in
\Lambda}$ of $W$, and for each $\lambda$ a smooth hypersurface
$\overline{W}_{\lambda}$ in $U_{\lambda}$ such that
$I(\overline{W}_{\lambda})\subset J_{\lambda}$, where
$J_{\lambda}$ denotes the restriction of $J$ to $U_{\lambda}$.
 \end{remark}
 For the proof of 3), note that an ideal of order one in a local
 regular ring $\O_{W,x}$ contains an element of order one; and that element defines a smooth hypersurface in some open neighborhood of
 $x\in W$.

\begin{remark}\label{cubri}
Fix, as before, an open covering of $W$, say $\{
U_{\lambda}\}_{\lambda \in \Lambda}$, and a monoidal
transformation with center $Y\subset W$, say $W \longleftarrow
W_1$. For each index $\lambda $ set $U_{\lambda}^{(1)}\subset W_1$
as the pull-back of $U_{\lambda}$. In this way we get
$$U_{\lambda} \longleftarrow U_{\lambda}^{(1)}$$ which is either a
monoidal transformation (in case $Y \cap U_{\lambda} \neq
\emptyset$), or the identity map (if $Y \cap U_{\lambda} =
\emptyset$). Note also that $\{ U_{\lambda}^{(1)}\}_{\lambda \in
\Lambda}$ is an open cover of $W_1$.

\end{remark}

\begin{proposition}\label{caso1}
Fix $J \subset \O_{W}$ with maximum order 1, and a sequence of
transformations
\begin{equation}\label{caso22}
\begin{array}{cccccccc}
 & W & \stackrel{\pi_1}{\longleftarrow} & W_{1} & \stackrel{\pi_2}{\longleftarrow}
 &\ldots& \stackrel{\pi_r}{\longleftarrow}&W_r\\ & (J,1)  & &(J_1,1) &&&&
 (J_r,1)\\
\end{array}
\end{equation}
then the maximum order of $J_r$ is either 1 or 0 (i.e. $J_r =
\O_{W_r}$ in the last case).
\end{proposition}

\begin{proof}
 Define an open covering $\{ U_{\lambda}\}_{\lambda \in
\Lambda}$ of $W$, and inclusions
\begin{equation}\label{eqdelcaso1}
 I(\overline{W}_{\lambda})
\subset J_{\lambda},
\end{equation}
where $\overline{W}_{\lambda}$ is a smooth hypersurface in
$U_{\lambda}$, as indicated in Remark \ref{cubrii},3).

The sequence (\ref{caso22}) defines, for each index $\lambda$, a
sequence of transformations:
\begin{equation*}
\begin{array}{cccccccc}
 & U_{\lambda}& \stackrel{\pi_1}{\longleftarrow} & U_{\lambda}^{(1)} & \stackrel{\pi_2}{\longleftarrow}
 &\ldots& \stackrel{\pi_r}{\longleftarrow}& U_{\lambda}^{(r)}\\ & (J_{\lambda},1)  & &((J_{\lambda})_1,1) &&&&
 ((J_{\lambda})_r,1),\\
\end{array}
\end{equation*}
and also
\begin{equation*}
\begin{array}{cccccccc}
 & U_{\lambda}& \stackrel{\pi_1}{\longleftarrow} & U_{\lambda}^{(1)} & \stackrel{\pi_2}{\longleftarrow}
 &\ldots& \stackrel{\pi_r}{\longleftarrow}& U_{\lambda}^{_{\lambda}(r)}\\ & (I(\overline{W}_{\lambda}),1)   & &((I(\overline{W}_{\lambda}))_1,1)
 &&&&
 ((I(\overline{W}_{\lambda}))_r,1).\\
\end{array}
\end{equation*}

Furthermore $$(I(\overline{W}_{\lambda}))_r \subset
(J_{\lambda})_r$$ by Remark \ref{remk1}. Let
$\overline{W}_{\lambda}^{(r)} \subset U^r_{\lambda}$ denote the
strict transform of $\overline{W}_{\lambda}$. Since
$\overline{W}_{\lambda}$ is smooth in $U_{\lambda}$, Example
\ref{weakvsstrict} asserts that $\overline{W}_{\lambda}^{(r)}$ is
smooth, and defined by the ideal $(I(\overline{W}_{\lambda}))_r$.
In particular $(I(\overline{W}_{\lambda}))_r$ has maximum order at
most one, and hence the same holds for $(J_{\lambda})_r$. Since
the open sets $(U_{\lambda})^{(r)}$ cover $W_r$ it follows that
$J_r$ has order at most 1.
\end{proof}

\begin{proposition}\label{casob}
Fix $J \subset \O_{W}$ with maximum order $b\geq 2$, and consider
a sequence of transformations
\begin{equation}
\begin{array}{cccccccc}\label{retu}
 & W & \stackrel{\pi_1}{\longleftarrow} & W_{1} & \stackrel{\pi_2}{\longleftarrow}
 &\ldots& \stackrel{\pi_r}{\longleftarrow}&W_r\\ & (J,b) & &(J_1,b)  &&&&
 (J_r,b). \\
\end{array}
\end{equation}
 Then then the maximum order of $J_r
(\subset \O_{W_r})$ is at most $b$.
\end{proposition}

\begin{proof}
From \ref{PropOrdDelta} we conclude that the maximum order of
$\Delta^{b-1}(J)(\subset \O_{W})$ is 1. Corollary \ref{cor2}
applied for $j=b-1$ says that (\ref{retu}) defines the sequence of
transformations

\begin{equation}\label{ecuacasob}
\begin{array}{cccccccc}
 & W & \stackrel{\pi_1}{\longleftarrow} & W_{1} & \stackrel{\pi_2}{\longleftarrow}
 &\ldots& \stackrel{\pi_r}{\longleftarrow}&W_r\\ & (\Delta^{b-1}(J),1)  & &((\Delta^{b-1}(J))_1, 1)  &&&&
 ((\Delta^{b-1}(J))_r, 1), \\
\end{array}
\end{equation}
 and that
$(\Delta^{b-1}(J))_r \subset \Delta^{b-1}(J_r)$. On the other hand
Proposition \ref{caso1} asserts that $(\Delta(J))_r$ has order at
most 1, and hence $ \Delta^{b-1}(J_r)$ has order at most one. From
this and \ref{PropOrdDelta} we conclude that $J_r$ has order at
most b.
\end{proof}

\begin{remark}\label{rk710} There is a stronger outcome that follows from the proof of Proposition
\ref{casob} that relates to induction in the dimension of the
ambient space. Note that $J$ has highest order $b$, so
$\Delta^{b-1}(J)$ has highest order one. We can argue as in the
proof of Proposition \ref{caso1}, and define an open cover $\{
U_{\lambda}\}_{\lambda \in \Lambda}$ of $W$, and for each index
$\lambda$, a smooth hypersurface $ \overline{W}_{\lambda}\subset
U_{\lambda}$, defined by
\begin{equation}\label{eq710}
I(\overline{W}_{\lambda}) \subset (\Delta^{b-1}(J))_{\lambda}.
\end{equation}

 Now use the compatibility of the $\Delta$ operator with
restriction to open sets and check that
$(\Delta^{b-1}(J))_{\lambda}=(\Delta^{b-1}(J_{\lambda}))$. Note
also that $\Sing(J,b)\cap U_{\lambda}\subset
\overline{W}_{\lambda}.$ Recall that (\ref{ecuacasob}) defines,
for each index $\lambda$, a sequence of transformations of
$((\Delta^{b-1}(J))_{\lambda},1)$, say:
\begin{equation*}
\begin{array}{cccccccc}
 & U_{\lambda}& \stackrel{\pi_1}{\longleftarrow} & U_{\lambda}^{(1)} & \stackrel{\pi_2}{\longleftarrow}
 &\ldots& \stackrel{\pi_r}{\longleftarrow}& U_{\lambda}^{(r)}\\ & ((\Delta^{b-1}(J))_{\lambda},1)  & &(((\Delta^{b-1}(J))_{\lambda})_1, 1) &&&&
 (((\Delta^{b-1}(J))_{\lambda})_r, 1),\\
\end{array}
\end{equation*}
and also
\begin{equation*}
\begin{array}{cccccccc}
 & U_{\lambda}& \stackrel{\pi_1}{\longleftarrow} & U_{\lambda}^{(1)} & \stackrel{\pi_2}{\longleftarrow}
 &\ldots& \stackrel{\pi_r}{\longleftarrow}& U_{\lambda}^{_{\lambda}(r)}\\ & (I(\overline{W}_{\lambda}),1)   & &((I(\overline{W}_{\lambda}))_1,1)
 &&&&
 ((I(\overline{W}_{\lambda}))_r,1).\\
\end{array}
\end{equation*}

Furthermore, $(I(\overline{W}_{\lambda}) )_r \subset
((\Delta^{b-1}(J))_{\lambda})_r$, and $(I(\overline{W}_{\lambda})
)_r$ defines a smooth hypersurface $\overline{W}_{\lambda}^{(r)}
\subset U_{\lambda}^{(r)}$ which is the strict transform of
$\overline{W}_{\lambda}$. We finally note that $\{
U_{\lambda}^{(r)}\}_{\lambda \in \Lambda}$ is a cover of
$W^{(r)}$, and taking restriction of the inclusion
$(\Delta^{b-1}(J))_r \subset \Delta^{b-1}(J_r)$, we get that:
$$((\Delta^{b-1}(J))_{\lambda})_r=((\Delta^{b-1}(J))_r)_{\lambda}
\subset (\Delta^{b-1}(J_r))_{\lambda},$$ and hence
$(I(\overline{W}_{\lambda}) )_r \subset
(\Delta^{b-1}(J_r))_{\lambda}$. In particular $$
(\Sing((J)_r,b)\cap U_{\lambda}^{(r)}=) \Sing((J_{\lambda})_r,b)
\subset \overline{W}_{\lambda}^{(r)}.$$

\end{remark}

\begin{lemma}\label{lematch}
Fix $J \subset \O_{W}$ with maximum order $b$. There is an open
covering, say $\{ U_{\lambda}\}_{\lambda \in \Lambda}$ of $W$, and
for each index $\lambda$ a smooth hypersurface $
\overline{W}_{\lambda}\subset U_{\lambda}$, such that the
following properties hold:

{\bf P1)} $\Sing(J_{\lambda},b) \subset \overline{W}_{\lambda}$.

{\bf P2)} For {\em any} sequence
\begin{equation}\label{ulambda}
\begin{array}{cccccccc}
 & W & \stackrel{\pi_1}{\longleftarrow} & W_{1} & \stackrel{\pi_2}{\longleftarrow}
 &\ldots& \stackrel{\pi_r}{\longleftarrow}&W_r\\ & (J,b) & &(J_1,b)  &&&&
 (J_r,b) \\
\end{array}
\end{equation}
 and setting by restriction, for each $\lambda$, say:
\begin{equation}\label{ulambdaa}
\begin{array}{cccccccc}
 & U_{\lambda}& \stackrel{\pi_1}{\longleftarrow} & U_{\lambda}^{(1)} & \stackrel{\pi_2}{\longleftarrow}
 &\ldots& \stackrel{\pi_r}{\longleftarrow}& U_{\lambda}^{(r)},\\ & (J_{\lambda},b)   &
 &((J_{\lambda})_1,b)
 &&&& ((J_{\lambda})_r,b)\\
\end{array}
\end{equation}
then $\{ U^{(r)}_{\lambda}\}_{\lambda \in \Lambda}$ is an open
covering of $W_r$, and
\begin{equation}\label{bulambda1}
\Sing(J_r,b)\cap U^{(r)}_{\lambda}= \Sing((J_{\lambda})_r,b)
\subset \overline{W}_{\lambda}^{(r)},
\end{equation}
where $\overline{W}_{\lambda}^{(r)}$ is the smooth hypersurface
defined by the strict transform of $\overline{W}_{\lambda}$.
\end{lemma}
\begin{proof} The case $b=1$ (in which $\Sing(J,1)=V(J)$) is in the proof of Proposition \ref{caso1}.
The case $b\geq 2$ is in Remark \ref{rk710}, and relies entirely
on the inclusion (\ref{eq710}).

\end{proof}

\begin{parrafo}\label{comentario}
Let $\overline{W}_{\lambda}^{(i)}$ denote the strict transform of
$\overline{W}_{\lambda}^{(0)}$ in $ U_{\lambda}^{(i)}$(see
(\ref{ulambdaa})). A consequence of (\ref{bulambda1}) is that all
the centers of monoidal transformations involved in
(\ref{ulambdaa}) are included in $\overline{W}_{\lambda}^{(i)}$;
hence (\ref{ulambdaa}) defines a sequence of monoidal
transformations
\begin{equation}\label{ulambda3}
\overline{W}_{\lambda} \longleftarrow
\overline{W}_{\lambda}^{(1)}\longleftarrow \cdots \longleftarrow
\overline{W}_{\lambda}^{(r)}.
\end{equation}
\end{parrafo}

\begin{proposition}\label{proptch}
Fix $J \subset \O_{W}$ with maximum order $b$. There is an open
covering, say $\{ U_{\lambda}\}_{\lambda \in \Lambda}$ of $W$, and
for each index $\lambda$ a closed and smooth hypersurface $
\overline{W}_{\lambda}\subset U_{\lambda}$, and a couple
$(K_{\lambda}^{(0)}, b!)$ with $K_{\lambda}^{(0)} \subset
\O_{\overline{W}_{\lambda}}$, such that, in addition to {\bf P1)}
and {\bf P2)} (\ref{lematch}), the following property holds:

\

{\bf P3)} The sequence (\ref{ulambda3}) defined by (\ref{ulambda})
as above, induces a sequence of transformations
\begin{equation}\label{ulambda1}
\begin{array}{cccccccc}
 & \overline{W}_{\lambda}  & \stackrel{\pi_1}{\longleftarrow} & \overline{W}_{\lambda}^{(1)} & \stackrel{\pi_2}{\longleftarrow}
 &\ldots& \stackrel{\pi_r}{\longleftarrow}&\overline{W}_{\lambda}^{(r)}\\ & (K_{\lambda},b!) & & ((K_{\lambda})_1,b!)  &&&&
 ((K_{\lambda})_r,b!) \\
\end{array}
\end{equation}
and
\begin{equation}\label{ulambdar}
 \Sing((J_{\lambda})_r,b) = \Sing((K_{\lambda})_r,b!)(\subset
\overline{W}_{\lambda}^{(r)}).
\end{equation}
\end{proposition}
\begin{remark}\label{rkdelconv}{\em On the converse.} Set
$W=U_{\lambda}$ so that $(J,b)=(J_{\lambda},b)$. The equality in
(\ref{ulambdar}) asserts, by induction on $r$, that any sequence
\ref{ulambda1} induces a sequence (\ref{ulambda}). And
furthermore, if \ref{ulambda1} is a resolution, so is
(\ref{ulambda}).

We are interested mainly in this converse, since we will argue by
increasing induction on the dimension of the ambient space. If we
accept, by induction, that there is a resolution \ref{ulambda1}
for each index $\lambda$, then we will have defined a resolution
(\ref{ulambdaa}) for each $\lambda$. We will define these
resolutions so that they patch to a resolution (\ref{ulambda}).

Full details of the proof of Proposition \ref{proptch} will be
given in Part IV, however the following example illustrates the
basic idea of the proof.
\end{remark}

\begin{example}\label{sigueee} In Example \ref{sigue} we considered the case $W=A_k^3=Spec(k[X,Y,Z])$, and $$J=<Z^3+X Y^2
Z+X^5>, $$ an ideal of maximum order $b=3$. In such example we
noted that $Z\in \Delta^2(J)=<Z, XY, Y^2, X^3>$, and we considered
the smooth hypersurface $\overline{W}=V(<Z>)$. This is a
particular example of Lemma \ref{lematch}, where there is no need
to consider the open covering $\{ U_{\lambda}\}_{\lambda \in
\Lambda}$ of $W$. In fact here the Lemma applies globally in $W$.
In this example $b!=6$, and Proposition \ref{proptch} applies by
setting $K=J^* $ as in (\ref{f1}).

A similar situations holds, more generally, in \ref{tschirn1}, for
$K=J^*=\langle c_i^{\frac{b!}{i}}, i=2,3,\dots,b \rangle. $
\end{example}

\begin{remark}\label{72local} The compatibility of the $\Delta$ operator with
open restrictions has played an important role in the proofs in
this section. If the transformation in Theorem \ref{theorem72} is
defined with center $Y\subset W$, and if $H\subset W_1$ denotes
the exceptional locus, then $J\O_{W_1}=I(H)^b\cdot J_1$, and $J_1$
has at most order $b$. Suppose now that the highest order of $J$
along points in $W$ is bigger than $b$, but that we simply know
that the order of $J$ is constant and equal to $b$ along any point
of the center $Y$. Since the order of $J$ along points in $W$
defines an upper-semi-continuous function on $W$, then there is an
open neighborhood, say $U \subset W$ of $Y$, so that $b$ is the
highest order of the restriction $J_U$. In particular there is an
open neighborhood $U_1$ of $H$ in $W_1$ so that the restriction
$(J_1)_{U_1}$ has highest order $\leq b$.

\end{remark}

\begin{remark}\label{rk715} The compatibility of the $\Delta$ operator with
open restrictions will also play a role in our proof of
Proposition \ref{proptch}, and this will allow us to present the
ideals $K_{\lambda}$ so that they are also compatible with a
restriction of $W$ to an open set $U$, at least if the restricted
ideal $J_U$ is again of highest order $b$.

There is yet another context in which there is a natural
compatibility of the operator $\Delta$, which are not open
restrictions, but will also play a role in the proof of
Proposition \ref{proptch}. In fact, set $W\longleftarrow
W_1=W\times A_k^1$ where $A_k^1$ denotes the affine line and the
map is the projection on the first coordinate. Note that if $J$ is
an ideal in $\O_W$, and if $\Delta_1$ denotes the operator on the
smooth scheme $W_1$, then
$$\Delta_1(J\O_{W_1})=\Delta(J)\O_{W_1}.$$ Note that a covering
$\{ U_{\lambda}\}_{\lambda \in \Lambda}$ of $W$ induces by
pull-back, a covering of $W_1$. The setting of Proposition
\ref{caso1} and the inclusions (\ref{eqdelcaso1}) are compatible
with pull-backs; and so are the setting of Proposition \ref{casob}
and the inclusions (\ref{eq710}). This will guarantee the
compatibility of all our development for this particular kind of
projection.
\end{remark}


\bigskip

\centerline {\bf Part III}

\bigskip

\section{On how the algorithm works. Examples.}
We finally sketch the main ideas and invariants involved in our
definition of Resolution Functions in \ref{resfunct}, which lead
us to the simple proofs of the Main Theorems in \ref{proofprin}
and \ref{proofclassical}. Recall the notion of permissible
sequence of transformations of pairs, say
\begin{equation*}
\label{resol}
\begin{array}{ccccccc}
W_{0} & \longleftarrow & (W_{1}, E_1=\{H_1\})& \longleftarrow &
\cdots & \longleftarrow & (W_{k},E_{k}=\{H_1,H_2,\dots,H_k \}),\\
Y & & Y_1 & & & &
\end{array}
\end{equation*}
in which we require that each monoidal transformation $W_i
\leftarrow W_{i+1}$ be defined so that all exceptional
hypersurfaces introduced have normal crossings (Prop.
\ref{hacer}).

Given $J\subset \O_W$, there is an expression of the total
transform (Def. \ref{deftransid}), say
$$J\O_{W_k}=I(H_1)^{a_1}I(H_2)^{a_2}\cdots I(H_k)^{a_k}\cdot \mathcal{A}_k.$$
This factorization is unique if we require the $a_i$ to be the
highest possible exponents in any such expression. In
\ref{proofprin} we want to achieve $\mathcal{A}_k=\O_{W_k}$ with
the conditions stated in Theorem \ref{principalization}. We will
argue in steps to achieve the proof of that theorem, each step
will introduce an exceptional hypersurface, and this will lead us
to consider a pair $(W,E=\{H_1,\dots,H_r\})$, rather then simply
$W$, and also permissible transformations of pairs
\begin{equation}\label{cerezas}
(W,E) \leftarrow (W_1,E_1)\leftarrow \cdots \leftarrow (W_k,E_k);
\end{equation}
always in the conditions of Prop. \ref{hacer}.

In \ref{normalc1} we have defined a basic object as a couple
$(J,b)$, $J\subset \O_W$, together with a pair $(W,E)$. A sequence
of transformations, say
\begin{equation}
\label{resol0}
\begin{array}{cccccccc}
(W,(J,b),E)& \leftarrow &(W_1,(J_1,b),E_1) & \leftarrow &\cdots
&\leftarrow & (W_k,(J_k,b),E_k),&
\\
\end{array}
\end{equation}
is a sequence of transformations of couples, say
\begin{equation*}
\begin{array}{cccccccc}
W& \leftarrow &W_1 & \leftarrow &\cdots &\leftarrow & W_k,&
\\(J,b)&&(J_1,b)& &  &  & (J_k,b) ,&
\\
\end{array}
\end{equation*}
(see (\ref{sectransp})), which also defines a sequence of
transformations of pairs, as in (\ref{cerezas}).

We shall say that (\ref{resol0}) is a resolution of $(W,(J,b),E)$
if $ V(\Delta^{b-1}(J_k))= \emptyset$. Note that
$V(\Delta^{b-1}(J_k))= \emptyset$ is equivalent to
$Sing(J_k,b)=\emptyset$, and also to the condition
$\max-ord_{J_k}<b$.

So the resolution would provide an expression of the form:
$$J\O_{W_k}=I(H_1)^{a_1}I(H_2)^{a_2}\cdots I(H_k)^{a_k}\cdot J_k
\mbox{,  and } \max-ord_{J_k}=b' < b.$$ If $b'=0$ we have achieved
what is stated in Theorem \ref{principalization}. If not we repeat
the argument, and try to produce a resolution of $(J_k,b')$ and
$(W_k,E_k)$. It is clear that ultimately we come to the case
$b'=0$.

Our task is to produce a resolution of $(J,b)$ and $(W,E)$, in
some explicit manner, in which centers of monoidal transformations
are defined by an upper-semi-continuous function. In some
particular cases this will be clear from the data involved (see
\ref{monomialcase} ). But, in general, the strategy will be to
reduce to the case in which $b=\max-ord_J$, namely to the case of
simple couples (\ref{nomeacuerdo}).

In case $b=\max-ord_J$, Theorem \ref{proptch} says that there is
$\overline{W}\subset W$, at least locally, and that (\ref{resol0})
induces
\begin{equation}
\label{resol2}
\begin{array}{cccccccc}
(\overline{W},(K,d),\overline{E}=\emptyset)& \leftarrow &
(\overline{W}_1,(K_1,d),\overline{E}_1) & \leftarrow &\cdots
&\leftarrow & (\overline{W}_k, (K_k,d),\overline{E}_k)&
\\
\end{array}
\end{equation}
such that $V(\Delta^{d-1}(K_k))= \emptyset$. It is important to
point out here that we will argue by induction, and hence we would
like to reverse the argument; namely, to define (\ref{resol0}) in
terms of (\ref{resol2}). We now indicate the difficulties to
overcome.

\bigskip

{\bf The three difficulties for an inductive argument:}

\bigskip

{\bf D1)} $(K,d)$ encodes information of $(J,b)$, but not of the
set of hypersurfaces  $E$ in $W$. Theorem \ref{proptch} asserts
that, after restriction to an open subset of $W$, (\ref{resol2})
will define a sequence of transformations of couples, say
\begin{equation*}
\begin{array}{cccccccc}
W& \leftarrow &W_1& \leftarrow &\cdots &\leftarrow & W_k&
\\(J,b)&&(J_1,b)& &  &  & (J_k,b), & \mbox{ }
\\
\end{array}
\end{equation*}
 such that $V(\Delta^{b-1}(J_k))=
\emptyset$. However this sequence might not define a sequence
(\ref{resol0}). In fact, it might not be permissible over $(W,E)$
because of the presence of hypersurfaces of $E$.

This is an important point to overcome. As indicated above, since
we will argue in steps, we introduce hypersurfaces with normal
crossings (those in $E$), and we want to preserve this condition
of normal crossings in all exceptional hypersurfaces to be
introduced in forthcoming steps.

\bigskip

{\bf D2)} The couple $(K,d)$ might not be a simple couple (might
not be such that $ d= \max-ord_K$). Take for example the case
$J=\langle z^3-x^2\cdot y^2\rangle$ and the couple $(J,3)$ in the
affine 3-space. Clearly $3=\max-ord_J$ so the couple is simple.
Since $z\in \Delta^2(J)$, we may take $\overline{W}$ as the affine
plane, and $(K,d)=(\langle x^2\cdot y^2 \rangle,3)$. Note that
$\max-ord_K=4,$ so that $(K,d)$ is not a simple couple
(\ref{nomeacuerdo}).

\bigskip

{\bf D3)} If $(J,b)$ is a simple couple (i.e. if $\max-ord_j=b$),
then $\overline{W}$ is defined by choosing, locally at a point $x
\in V(\Delta^{b-1}(J))$, an element of order one in
$\Delta^{b-1}(J)_x$. In general this choice is not unique, and the
definition of $(K,d)$ ($K\subset \calo_{\overline{W}})$ is local
at $x$. Our form of induction should provide a resolution
\ref{resol0}, with independence of open restrictions and of
choices of $\overline{W}$.

\bigskip

\begin{parrafo}
Set $J\subset \O_W$ and $(W,E)$ as before. Assume, in accordance
with D2), that $b\geq \max-ord_J$. So here $(J,b)$ might not be
simple. Consider a sequence of transformations, say:
\begin{equation}
\label{resol3}
\begin{array}{cccccccc}
(W,(J,b),E)& \leftarrow &(W_1,(J_1,b),E_1) & \leftarrow & \dots
&\leftarrow & (W_s,(J_s,b),E_s).&
\\
\end{array}
\end{equation}
We claim that this provides a factorization of $J_s$, say
$$J_s=I(H_1)^{b_1}I(H_2)^{b_2}\cdots I(H_s)^{b_s}\cdot
\overline{J}_s$$ so that $\overline{J}_s$ does not vanish along
$H_i$, $1\leq i \leq s$. In this manner we may consider $(J_s,b)$,
together with this factorization of $J_s$. This extra structure on
$(J_s,b)$ will allow us to overcome D2), namely to reduce the
general case to the case of simple couples.
\end{parrafo}

\medskip

\begin{example} Set $J=<x_1, x_2^2>^4$, $W=\mathbb{A}^2_k$
\begin{equation*}
\begin{array}{cccccc}
(W,(J,3),E=\emptyset)& \leftarrow
&(W_1,(J_1,3),E_1=\{H_1\})&\leftarrow&&(W_2,(J_2,3),E_2=\{
H_1,H_2\})
\\&&J\O_{W_1}=I(H_1)^4\cdot \mathcal{M}_p^4&&&J_1 \O_{W_2}=I(H_1)^2\cdot I(H_2)^6\\
&&J_1=I(H_1)^2 \mathcal{M}_p^4&&&J_2=I(H_1)^2\cdot I(H_2)^3.\\
\end{array}
\end{equation*}
\end{example}
Here $W \leftarrow W_1$ is the blow-up at $0\in \mathbb{A}^2_k$,
$p\in W_1$ is a point in the exceptional line $H_1$,
$\mathcal{M}_p$ is the sheaf of functions that vanish at $p$, and
finally $W_1 \leftarrow W_2$ is the blow-up at $p$.

\begin{remark}\label{monomialcase}  If $\overline{J}_s=\O_{W_s}$, we say that $(J_s,b)$ is within the {\it monomial case}.
In this case it is easy to extend (\ref{resol3}) to a resolution;
namely, to define for some $k\geq s$:
\begin{equation}
\label{resol4}
\begin{array}{ccccccccc}
(W,(J,b),E)&\leftarrow \cdots &\leftarrow &(W_s,(J_s,b),E_s) &
\leftarrow & \cdots &\leftarrow & (W_k,(J_k,b),E_k)&
\\
\end{array}
\end{equation}
so that $V(\Delta^{b-1}(J_k))= \emptyset$. The following example
illustrates this fact. Note that in the previous example
$\overline{J}_2=\O_{W_r}$.
\end{remark}

\begin{example}Consider transformations with centers $Y_j$:
\begin{equation*}
\begin{array}{cccccc}
(W_2,(J_2,3),E_2=\{H_1,H_2\})&  \stackrel{id}{\leftarrow}
&(W_3,(J_3,3),E_3=\{H_1,H_2\})&\leftarrow&&(W_4,(J_4,3),E_4=\{
H_1,H_2,H_4\})
\\J_2=I(H_1)^2\cdot I(H_2)^3&&J_3=I(H_1)^2\cdot I(H_2)&&&J_3=I(H_1)^2\cdot I(H_4)^0 \cdot I(H_2)\\
Y_2=H_2&&Y_3=H_1\cap H_2&&&V(\Delta^2(J_4))=\emptyset.\\
\end{array}
\end{equation*}
The first transformation is defined with center at the
hypersurface $H_2$. So the morphism is the identity map, but we
take here $H_2\in E_2$ to be the exceptional locus. Note that
$J_3$ is not $J_2$.
\end{example}

\bigskip

\begin{parrafo}{\bf On the function $\vord$.}

\bigskip

Given a sequence of transformations of basic objects, say
(\ref{resol3}), we have defined an expression:
$$J_s=I(H_1)^{b_1}I(H_2)^{b_2}\cdots I(H_s)^{b_s}\cdot
\overline{J}_s$$ so that $\overline{J}_s$ does not vanish along
$H_i$, $1\leq i \leq s$. Define now:

$$ \vord_s: V(\Delta^{b-1}(J_s)) \to {\mathbb N}$$ $$ \vord_s(x)=\nu_x( \overline{J}_s),\mbox{
(the order of } (\overline{J}_s)_x \mbox{ at } \O_{W_s,x}).$$

\bigskip
Note that:

 1) The function is
upper-semi-continuous. In particular $\Max \vord$ is closed.
\medskip

2) For any index $ i \leq s$, there is an expression
$$J_i=I(H_1)^{b_1}I(H_2)^{b_2}\cdots I(H_i)^{b_i}\cdot
\overline{J}_i,$$ and hence a function $\vord_i:
V(\Delta^{b-1}(J_s)) \to {\mathbb N}$ can be defined.

\medskip

Another property of these functions is:

3)  If each step $(W_i,E_i) \leftarrow (W_{i+1},E_{i+1})$ in
(\ref{resol3}) is defined with center $Y_i \subset \Max \vord_i $,
then $$ \max \vord \geq \max \vord_1\geq \dots \geq \max
\vord_s.$$
\end{parrafo}
\medskip

Property 3) follows from the fact that, if $\max \vord_s=b'$, then
$\Max \vord_s=V(\Delta^{b'-1}(\overline{J_s}))$ (the closed set of
$(\overline{J_s},b'))$, where $(\overline{J_s},b')$ is, by
definition, a simple couple.

\begin{example}\label{curvita} Set $(J,1)$; $J=<x^2-y^5>$, and $W=\mathbb{A}^2_k$. Let $C$
denote the curve defined by $J$. \tiny
\begin{equation*}
\begin{array}{ccccccccc}
(W,(J,1),E=\emptyset)& \leftarrow
&(W_1,(J_1,1),E_1=\{H_1\})&\leftarrow&&(W_2,(J_2,1),E_2=\{
H_1,H_2\})&\leftarrow && (W_3,(J_3,1),E_3)
\\&&J_1=I(H_1)\cdot I(C')&&&J_2=I(H_1)^1\cdot I(H_2)^1 \cdot I(C'')&&&J_3=I(H_1)^1\cdot I(H_2)^1\cdot I(H_3)^2 \cdot I(C''')\\
Y=0\in \mathbb{A}^2_k&&Y_1=C'\cap H_1 &&&Y_2=H_1 \cap H_2.&&&\\
\end{array}
\end{equation*}
\end{example}

Here the $Y_i$ are the centers of the monoidal transformations,
and $C'$, $C''$, and $C'''$ are strict transforms of $C$. In this
example $$ \max \vord_J=2 \ ; \ \max \vord_{J_1}=1 \ ; \max
\vord_{J_2}=1;\max \vord_{J_3}=1;$$ and the sequence is defined by
setting $$Y=\Max \vord_J=Y ; Y_1=\Max \vord_{J_1} \mbox{, and }
Y_2=\Max \vord_{J_2}.$$

\medskip

\begin{parrafo} \label{poift} {\bf On the inductive function t.}

Consider, as before, a sequence
\begin{equation}
\label{resol6}
\begin{array}{cccccccc}
(W,(J,b),E)& \leftarrow &(W_1,(J_1,b),E_1) & \leftarrow &\cdots
&\leftarrow & (W_s,(J_s,b),E_s),&
\\
\end{array}
\end{equation}
where each $W_i\leftarrow W_{i+1}$ is defined with center $Y_i
\subset \Max \vord_i$, so that: $$ \max \vord \geq \max
\vord_1\geq \dots \geq \max \vord_s.$$ Set $s_0 \leq s$ such that
$$ \max \vord \geq \dots \geq \max \vord_{s_{0-1}}> \max
\vord_{s_0}= \max \vord_{s_{0+1}}= \dots= \max \vord_{s},$$ and $$
E_s=E_s^+ \sqcup E_s^- \mbox{ (disjoint union)},$$ where $E_s^-$
are the strict transform of hypersurfaces in $E_{s_0}$. Define
$$ t_s: V(\Delta^{b-1}(J_s)) \to {\mathbb N} \times {\mathbb N}  \mbox{ (ordered
lexicographically). }$$
$$t_s(x)=(\vord_s(x), n_s(x))$$ $$ n_s(x)= \sharp \{H_i\in E_s^-, /
x\in H_i\}$$

\medskip

One can check that:

1) the function is upper-semi-continuous. In particular $\Max t_s$
is closed.

2) There is a function $t_i$ for any index $ i \leq s$.

\bigskip

Example \ref{exabc} illustrates the following properties which
also hold for this function:

$\bullet$ If each $(W_i,E_i) \leftarrow (W_{i+1},E_{i+1})$ in
\ref{resol6} is defined with center $Y_i \subset \Max t_i $, then
$$ \max t \geq \max t_1\geq \dots \geq \max t_s.$$

$\bullet$ If $\max t_s=( b', r)$ (here $\max \vord_s=b'$) then
$\Max t_s \subset \Max \vord_s$.

$\bullet$ If $\Max t_s$ has codimension 1 in $W_s$, then it is
smooth. Moreover, in such case $Y_s=\Max t_s$ is a permissible
center, defining
\begin{equation*}
\begin{array}{cccccccc}
(W_s,(J_s,b),E_s)& \leftarrow &(W_{s+1},(J_{s+1},b),E_{s+1}), &  &
& &&
\\
\end{array}
\end{equation*}
 and $ \max t_s > \max t_{s+1}$ (hence $\max \vord_s \geq
\max \vord_{s+1}).$
\begin{example}\label{exabc}

\

{\bf 0)} Consider $(J,1)$; $J=<x^2-y^3>$ defining a curve
$C\subset W=\mathbb{A}^2_k$.

Here $t(x)=(1,0) $ at any $x\in C$ except at $\overline{0}\in
\mathbb{A}^2_k$, $t(\overline{0})=(2,0)$. So
 $$\max t=(2,0) \mbox{ and } \Max t=\overline{0}\in
\mathbb{A}^2_k. \mbox{ Let now }$$
\begin{equation*}
\begin{array}{cccc}
(W,(J,1),E=\emptyset)& \leftarrow &(W_{1},(J_1,1),E_{1}=\{H_1\}) &
\\
\end{array}
\end{equation*}
be the quadratic transformation at $\overline{0}\in
\mathbb{A}^2_k$.

{\bf 1)} Let $C'\subset W_1$ denote the strict transform of $C$.
Here $$J_1=I(H_1)\cdot \overline{J}_1$$ where
$\overline{J}_1=I(C')$, and $t_1(x)=(1,0) $ at any $x\in C'$
except for $p=C'\cap H_1$, where $t_1(p)=(1,1)$. So
 $$\max t_1=(1,1) \mbox{ and } \Max t_1= p.$$ Set
\begin{equation*}
\begin{array}{cccc}
(W_1,(J_{1},1),E_1)& \leftarrow
&(W_{2},(J_{2},1),E_{2}=\{H_1,H_2\}) &
\\
\end{array}
\end{equation*}
with center at $p\in W_1$.
\bigskip

{\bf 2)} If $C''\subset W_2$ denotes the strict transform of $C$,
$$J_2=I(H_1)\cdot I(H_2)\cdot \overline{J}_2$$ where
$\overline{J}_2=I(C'')$.

Now $t_2(x)=(1,0) $ at any $x\in C''$ except for $q=C''\cap
H_1\cap H_2 $, where $t_2(q)=(1,1)$. So $$\max t_2=(1,1) \mbox{
and } \Max t_2= q.$$ Set
\begin{equation*}
\begin{array}{cccc}
(W_{2},(J_{2},1),E_{2}=\{H_1,H_2\})& \leftarrow
&(W_{3},(J_{3},1),E_{3}=\{H_1,H_2, H_3\}) &
\\
\end{array}
\end{equation*}
with center at $q\in W_2$.
\bigskip

{\bf 3)} Now $$J_3=I(H_1)\cdot I(H_2)\cdot I(H_3)^2\cdot
\overline{J}_3$$ where $\overline{J}_3=I(C''')$ (ideal of the
strict transform). Finally check that $t_3(x)=(1,0) $ at any $x\in
C'''$. So
$$\max t_3=(1,0) \mbox{ and } \Max t_3= C'''. $$

 This is a case in which $\Max t$ has codimension 1. Note
 that $\Max t_3$ is a smooth hypersurface, and the blow-up at $\Max t_3$ defines a
 permissible transformation (the identity map):
\begin{equation*}
\begin{array}{cccc}
(W_{3},(J_3,1),E_{3}=\{H_1,H_2,H_3\})& \leftarrow

&(W_{3},(J_4,1),E_{3}=\{H_1,H_2, H_3\}) &
\\
\end{array}
\end{equation*}
with $J_4=I(H_1)\cdot I(H_2)\cdot I(H_3)^2$.
\end{example}
\end{parrafo}

\parrafo\label{overcom}{\bf Overcoming difficulties D1) and D2)}

We finally indicate a further property of the function $t_s$,
which leads to constructive desingularization by induction. To
this end set:
\begin{equation}
\label{resol7}
\begin{array}{cccccccc}
(W,(J,b),E)& \leftarrow &(W_1,(J_1,b),E_1) & \leftarrow &\cdots
&\leftarrow & (W_s,(J_s,b),E_s)&
\\
\end{array}
\end{equation}
so that $$ \max \vord \geq  \max \vord_1\geq \dots \geq \max
\vord_s.$$ And define, as before, the function $$ t_s:
V(\Delta^{b-1}(J_s)) \to {\mathbb N}\times {\mathbb N}.$$
\bigskip
This last property can be stated as follows:

There is a couple $(J''_s, b'')$ with the following properties:
\bigskip

$\bullet$ $V(\Delta^{b''-1}(J''))= \Max t_s$, and $\max
ord_{J''_s}=b''$ (i.e. the couple is a simple couple).

\bigskip

$\bullet$ Let $\overline{W}_s$ be a smooth hypersurface containing
$V(\Delta^{b''-1}(J''))$, and set $(K,d)$ ($K\subset
\O_{\overline{W}_s})$ as in Proposition \ref{proptch}. Then any
resolution, say:
\begin{equation}
\label{presol8}
\begin{array}{cccccccc}
(\overline{W}, (K,d),\overline{E}_s=\emptyset)& \leftarrow &
(\overline{W}_1 ,(K_1,d), \overline{E}_{s+1})& \leftarrow &\cdots
&\leftarrow & (\overline{W}_k, (K_k,d),\overline{E}_{s+k})&
\\
\end{array}
\end{equation}
($ V(\Delta^{d-1}(K_k))= \emptyset$), induces an extension of
(\ref{resol7}), say:
\begin{equation}
\label{resol9}
\begin{array}{cccccccc}
(W_s,(J_{s},b),E_s)& \leftarrow &(W_{s+1},(J_{s+1},b),E_{s+1}) &
\leftarrow &\cdots &\leftarrow & (W_{s+k},(J_{s+k},b),E_{s+k}),&
\\
\end{array}
\end{equation}
\bigskip
such that $$\max t_s= \max t_{s+1}=\dots = \max t_{s+k-1}> \max
t_{s+k}.$$ Furthermore $$ \Max t_{s+i}=V(\Delta^{d-1}(K_i)) (
=\Sing(K_i,d)) $$ for $ i=0,\cdots, k-1$.

\begin{parrafo}{\bf Example of constructive resolution.}

\begin{example} The curve $C$ defined by  $J=<x^2-y^5>$ in
$W=\mathbb{A}^2_k$ is irreducible, in particular reduced. We
attach to it the basic object
\begin{equation*}
\begin{array}{cccc}
 &(W,(J,1),E=\emptyset),&&\\
\end{array}
\end{equation*}
and the function $$ t: V(\Delta^0(J))=V(J) \to {\mathbb N}\times
{\mathbb N}.$$ Here $t(x)=(1,0)$ except at the origin
$\overline{0}\in \mathbb{A}^2_k$, $t(\overline{0})=(2,0)$.

Note that in Example \ref{curvita}:

$\bullet$ $\max t=(2,0)$ and $Y=\Max t=\overline{0}$;

$\bullet$ $\max t_1=(1,1)$ and $Y_1= \Max t_1$;

$\bullet$ $\max t_2=(1,1)$ and $Y_2=\Max t_2$

$\bullet$ $\max t_3= (1,0)$, $\Max t_3=C'''$, is a smooth
hypersurface (see \ref{overcom}). Thus, this defines an embedded
desingularization.

 Compare with the proof of Theorem \ref{classical}.

\end{example}
\begin{example}
The hypersurface $Z^2+X^2+ Y^3=0$ is irreducible with an isolated
singularity at $\overline{0}\in \mathbb{A}^3_k$. Set
$W=\mathbb{A}^3_k$ , $J=<Z^2+X^2+ Y^3>$. According to the proof of
Theorem \ref{classical} in \ref{proofclassical},
desingularization is achieved at some intermediate step of the
resolution of the basic object:
\begin{equation*}
\begin{array}{cccc}
 &(W,(J,1),E=\emptyset).&&\\
\end{array}
\end{equation*}
The function $t: V(J) \to {\mathbb N}\times {\mathbb N}$ takes
value $t(x)=(1,0)$ except at the singular point,
$t(\overline{0})=(2,0)$. In this case, and following the notation
in \ref{overcom}:

$\bullet$ $\max t=(2,0)$.

$\bullet$ $(J'',b'')$ can be defined as $(J,2)$.

$\bullet$ $\overline{W}=V(<Z>)$ (in fact $Z\in \Delta^1(J)$).

$\bullet$ $(K,d)$ can be defined by $(<X^2+ Y^3>,2)$.

Here $\overline{W}=\mathbb{A}^2_k$, and the blow-up at
$\overline{0}\in \mathbb{A}^2_k$ defines a resolution, namely
\begin{equation*}
\begin{array}{ccccc}
&(\overline{W},(K,2),E=\emptyset)&\longleftarrow &(\overline{W}_1,(K_1,2),E_1=\{\overline{H}_1\})&\\
\end{array}
\end{equation*}
 and $ V(\Delta(K_1))=\emptyset$.
 According to \ref{overcom}, this
defines
\begin{equation*}
\begin{array}{ccccc}
&(W,(J,1),E=\emptyset)&\longleftarrow &(W_1,(J_1,1),E_1=\{H_1\})&\mbox{ and } \max t > \max t_1.\\
\end{array}
\end{equation*}

In fact $\max t_1=(1,1)$. So again, we argue as in \ref{overcom},
and attach a couple $(J'',b'')$ to the value $\max t_1=(1,1)$.
Moreover, a smooth hypersurface $\overline{W}$ and a couple
$(K,d)$ can be defined so that a resolution, say:
\begin{equation}
\label{resol8}
\begin{array}{cccccccc}
(\overline{W}, (K,d),\overline{E}_s=\emptyset)& \leftarrow &
(\overline{W}_1 ,(K_1,d) \overline{E}_{s+1})& \leftarrow &\cdots
&\leftarrow & (\overline{W}_k, (K_k,d),\overline{E}_{s+k})&
\\
\end{array}
\end{equation}
(such that $ V(\Delta^{d-1}(K_k))= \emptyset$), induces:
\begin{equation}
\label{resol19}
\begin{array}{cccccccc}
(W_1,(J_{1},1),E_1)& \leftarrow &(W_{2},(J_{2},1),E_{2}) &
\leftarrow &\cdots &\leftarrow & (W_{s},(J_{s},1),E_{s})&
\\
\end{array}
\end{equation}
such that $$(1,1)=\max t_1= \max t_{2}=\dots = \max t_{s-1}> \max
t_{s}=(1,0).$$ Note that $J_s$ is the sheaf of ideals of the
strict transform of the hypersurface, that $\Max t_s=V(J_s)$. So
$\Max t_s$ is a hypersurface, and the last property in \ref{poift}
says that this is an embedded desingularization.
\end{example}

\end{parrafo}



\centerline {\bf Part IV}

In this Part we will address constructive resolution in detail.
Part III was devoted to give an overview of the invariants
involved, and examples of constructive resolution. This last Part
IV can be read independently of the previous one, so we will
introduce all invariants, and prove resolution theorems in full
generality.

\begin{section} {Tchirnhausen revisited.}\label{Tchirnhausen revisited}
The objective of this Section is to prove Proposition
\ref{proptch} (see also \ref{demproptch}), which is the form of
induction that leads to resolution. This form of induction is that
suggested by the examples in Section \ref{tchirn}.

 In Example \ref{example2} we treated a case of a simple basic object where $W=A^3$, and
$b=3$. There the covering $\{ U_{\lambda}\}_{\lambda \in \Lambda}$
is trivial (i.e. $\{ U_{\lambda}= W \}$), and $Z \in
\Delta^{(2)}(J)$ defines a smooth hypersurface
$\overline{W}=V(<Z>) (\subset W)$. Moreover, in that example the
couple $(J^*, 6)$ ( $J^* \subset \O_{\overline{W}}$ in (\ref{f1}))
plays the role of $ (K_{\lambda},b!)$ with property P3) in
Proposition \ref{proptch}, to be defined in \ref{demproptch}.

\begin{remark} \label{RemCoeffComplet}
     We will assume here that the setting of Remark \ref{rk710} holds for $U_{\lambda}= W $,
but in a more general form, where the role of the smooth
hypersurface $\overline{W}$ is played now by an arbitrary smooth
subscheme, say \( Z\subset W \). In other words, assume that $b$
is the highest order of $J \subset \O_W$, and that for any
sequence of transformations of couples, say
\begin{equation}\label{ecW}
\begin{array}{cccccccc}
W& \leftarrow &W^{(1)}& \leftarrow &\cdots &\leftarrow & W^{(r)}&
\\(J,b)&&(J_1,b)& &  &  & (J_r,b), & \mbox{ }
\\
\end{array}
\end{equation}
then
\begin{equation}\label{ecCc}
\Sing((J)_r,b) \subset Z^{(r)},
\end{equation}
where $Z^{(r)}$ is the smooth subscheme in $W^{(r)}$ defined by
the strict transform of $Z$.

Note, in particular, that (\ref{ecW}) induces a sequence of
monoidal transformations
\begin{equation}\label{ecCcc}
 Z \longleftarrow
Z^{(1)}\longleftarrow \cdots \longleftarrow Z^{(r)}.
\end{equation}

 Let \( \xi\in
        Z \) be a
closed point, and let $\{
        z_{1},\ldots,z_{r},x_{1},\ldots,x_{n}
\}$ be a regular system of
        parameters in \( \mathcal{O}_{W,\xi} \)
such that \(
        \mathcal{I}(Z)_{\xi}=(z_{1},\ldots,z_{r}) \).
Consider the isomorphisms
        \begin{equation*}
\hat{\mathcal{O}}_{W,\xi}\cong
k(\xi)[[z_{1},\ldots,z_{r},x_{1},\ldots,x_{n}]],
      \qquad
\hat{\mathcal{O}}_{Z,\xi}\cong k(\xi)[[x_{1},\ldots,x_{n}]],
\end{equation*}
where the right hand sides are rings of formal series. Given \(
f\in\mathcal{O}_{W,\xi} \), let \( \hat{f} \) denote the image at
\( \hat{\mathcal{O}}_{W,\xi} \), say:
\begin{equation*}
      \hat{f}=\sum_{i_{1},\ldots,i_{r}=0}^{\infty}
a_{i_{1},\ldots,i_{r}}z_{1}^{i_{1}}\cdots z_{r}^{i_{r}},
\end{equation*}
       where each \( a_{i_{1},\ldots,i_{r}}\in
k(\xi)[[x_{1},\ldots,x_{n}]] \).  Note that
\begin{equation}
\label{EqParcCoeff}
      (i_{1}!\cdots i_{r}!)a_{i_{1},\ldots,i_{r}}=
\varphi\left(
      \frac{\partial^{i_{1}+\cdots+i_{r}} f}{\partial
z_{1}^{i_{1}}\cdots
      \partial z_{r}^{i_{r}}}
      \right),
\end{equation}
where \(
\varphi:k(\xi)[[z_{1},\ldots,z_{r},x_{1},\ldots,x_{n}]]\rightarrow
k(\xi)[[x_{1},\ldots,x_{n}]] \) is the quotient map induced by the
inclusion \( Z\subset W \) at \( \xi \). Note also that, for a
fixed integer \( b \),
        \begin{equation} \label{EqCoeffOrd}
\nu_{\xi}(f)\geq b
      \Longleftrightarrow
\nu_{\xi}(a_{i_{1},\ldots,i_{r}})\geq b-(i_{1}+\cdots+i_{r}),
\end{equation}
for all $ i_{1},\ldots,i_{r}$ with $0\leq i_{1}+\cdots+i_{r}<b$
(here the left hand side is the order at \(\calo_{W, \xi}\), and
the right hand side is the order at \(\calo_{Z, \xi}\)). Set now
\begin{equation}\label{eccoef}
I(f,b)=\langle
(a_{i_{1},\ldots,i_{r}})^{\frac{b!}{b-(i_{1}+\cdots+i_{r})}} /
0\leq i_{1}+\cdots+i_{r}<b \rangle
\end{equation}
and reformulate (\ref{EqCoeffOrd}) by means of the equivalence
 \begin{equation} \label{EqCoeffOrd1}
\nu_{\xi}(f)\geq b
      \Longleftrightarrow
\nu_{\xi}(I(f,b))\geq b!.
\end{equation}
\end{remark}

\begin{lemma}\label{lemita}

 Assume now that:

 1) $f$ and $\{z_{1},\ldots,z_{r},y_{1},\ldots,y_{n} \}$ are global sections
of $\O_W$ and $J=<f>$,

2) the sheaf of differentials $\Omega_W$ is free with basis $\{ d(
z_{1}),\dots ,d( z_{r}),d( y_{1}), \dots ,d( y_{n})\}$,

3) $I(Z)=<z_{1},\ldots,z_{r}>$ and $Z$ fulfills the property
expressed in (\ref{ecCc}) for $(<f>,b)$.

Then there is a couple $(J^*, d)$, with $J^* \subset \O_{ Z}$,
such that any sequence of transformations of $(J,b)$, say
(\ref{ecW}), induces a sequence of transformations for the couple
$(J^*, d)$, say

\begin{equation}\label{ecC1}
\begin{array}{cccccccc}
Z& \leftarrow &Z^{(1)}& \leftarrow &\cdots &\leftarrow & Z^{(r)}&
\\(J^*,d)&&(J^*_1,d)& &  &  & (J^*_r,d), & \mbox{ }
\\
\end{array}
\end{equation}
 and $$ \Sing((J)_r,b)=\Sing((J^*)_r,d) (\subset
Z^{(r)}).$$ Conversely, any sequence (\ref{ecC1}) induces a
sequence (\ref{ecW}).

\end{lemma}
\begin{proof} For the converse stated in the last line see \ref{rkdelconv}.

If $g$ is a global section in $\O_W$, let $\overline{g}$ denote
the class in $\O_Z$. Set
  $$J^*=\langle \left( \frac{1}{i_{1}!\cdots i_{r}!}
\overline{\left( \frac{\partial^{i_{1}+\cdots+i_{r}} f}{\partial
z_{1}^{i_{1}}\cdots
      \partial z_{r}^{i_{r}}}\right)} \right)^{\frac{b!}{b-(i_{1}+\cdots+i_{r})}} / 0\leq
i_{1}+\cdots+i_{r}<b \  \rangle.$$

Fix a closed point  \( \xi\in Z \) with residue field $ k(\xi)$,
and let $\alpha_i \in k(\xi)$ denote the class of $y_i$ at the
point. Set $\{z_{1},\ldots,z_{r},x_{1},\ldots,x_{n} \} \subset
\hat{\mathcal{O}}_{W,\xi}$, where $x_i=y_i-\alpha_i  $.  Note also
that, despite the change of coordinates, the global derivations
$\frac{\partial}{\partial z_i}$ (defined in terms of
$\{z_{1},\ldots,z_{r},y_{1},\ldots,y_{n} \}$) induces the
derivation $\frac{\partial}{\partial z_i}$ on
$\hat{\mathcal{O}}_{W,\xi}$ defined in terms of
$\{z_{1},\ldots,z_{r},x_{1},\ldots,x_{n} \}$.

It follows now from (\ref{EqParcCoeff}) and (\ref{EqCoeffOrd1})
that, setting $d=b!$, $\Sing(J,b)=\Sing(J^*,d) (\subset Z)$. \

 Fix a closed point $\xi_r \in Z^{(r)}$, and set \( \xi_{k} \)  as the image of \( \xi_{r}
\) in  \( Z^{(k)} \). In particular \( \xi_{0}\in Z^{(0)}=Z \). We
may assume, by induction, that:

1) there is a regular system of parameters
$\{z_{k-1,1},\ldots,z_{k-1,r},x_{k-1,1},\ldots,x_{k-1,n}\}$ at
$\hat{\mathcal{O}}_{W_{{k-1}},\xi_{k-1}}\cong R_{k-1}=
k(\xi_{k-1})[[z_{k-1,1},\ldots,z_{k-1,r},x_{k-1,1},\ldots,x_{k-1,n}]],$

2) $\hat{I(Z^{(k-1)})}=<z_{k-1,1},\ldots,z_{k-1,r}>$, and

3) there is a generator $\hat{f}_{k-1}$ of $(J)_{k-1}$, together
with an expression:
\begin{equation*}
         \hat{f}_{k-1}^{}=
\sum_{i_{1},\ldots,i_{r}=0}^{\infty} a_{k-1,i_{1},\ldots,i_{r}}^{}
         z_{k-1,1}^{i_{1}}\cdots
z_{k-1,r}^{i_{r}},
    \end{equation*}
with $a_{k-1,i_{1},\ldots,i_{r}} \in
k(\xi_{k-1})[[x_{k-1,1},\ldots,x_{k-1,n}]](\subset
\hat{\mathcal{O}}_{W_{{k-1}},\xi_{k-1}}) $.

Note that there is a natural identification of the subring
$k(\xi_{k-1})[[x_{k-1,1},\ldots,x_{k-1,n}]]$ with the quotient
ring $\hat{\mathcal{O}}_{Z_{k-1},\xi_{k-1}}$; and using 1), 2),
and 3), define $I(\hat{f}_{k-1},b)\subset
\hat{\mathcal{O}}_{Z_{k-1},\xi_{k-1}}$ as in (\ref{eccoef}).

A change of coordinates in the subring $R_{k-1}$, extends to a
change of coordinates at $\hat{\mathcal{O}}_{W_{{k-1}},\xi_{k-1}}$
by fixing $\{z_{k-1,1},\ldots,z_{k-1,r}\}$.

This particular kind of change of coordinates in
$\hat{\mathcal{O}}_{W_{{k-1}},\xi_{k-1}}$ fixes the ideal in 2),
and modifies the expression in 3) by changing each coefficient
$a_{k-1,i_{1},\ldots,i_{r}} \in
k(\xi_{k-1})[[x_{k-1,1},\ldots,x_{k-1,n}]]$.

The induced change of coordinates in the quotient ring
$\hat{\mathcal{O}}_{Z_{k-1},\xi_{k-1}}$ is compatible with our
definition of the ideal $I(\hat{f}_{k-1},b)$, defined in terms of
expression 3). The point is that, after enlarging $k(\xi_{k-1})$
to $k(\xi_{k})$, and taking a suitable change of coordinates as
before, we may choose

1') coordinates $
\{z_{k,1},\ldots,z_{k,r},x_{k,1},\ldots,x_{k,n}\}$ in
$\hat{\mathcal{O}}_{W_{k},\xi_{k}}$ with
\begin{equation*}
     z_{k,i}=\frac{z_{k-1,i}}{x_{k-1,1}},
\quad i=1,\ldots,r, \qquad x_{k,1}=x_{k-1,1},\quad
x_{k,i}=\frac{x_{k-1,i}}{x_{k-1,1}}, \qquad i=1,\ldots,n,
\end{equation*} so that

2') \( \mathcal{I}(Z_{k})_{\xi_{k}}= \left\langle
z_{k,1},\ldots,z_{k,r}\right\rangle \).  Set the expression

3') $ \hat{f}_{k}^{}= \frac{\hat{f}_{k-1}^{}}{x_{k-1,1}^{b}}=
\sum_{i_{1},\ldots,i_{r}=0}^{\infty} a_{k,i_{1},\ldots,i_{r}}^{}
z_{k,1}^{i_{1}}\cdots z_{k,r}^{i_{r}}, \mbox{ } \mbox{ }   \mbox{
where } \mbox{ } \mbox{ } a_{k,i_{1},\ldots,i_{r}}^{}=
\frac{a_{k-1,i_{1},\ldots,i_{r}}^{}}{x_{k-1,1}^{b-
(i_{1}+\cdots+i_{r})}}.$

 Note here that $\hat{f}_{k}$ is a generator of $(J)_k$. Furthermore, since
      $$
     \left(
a_{k,i_{1},\ldots,i_{r}}^{}
\right)^{\frac{b!}{b-(i_{1}+\cdots+i_{r})}}=\frac{\left(
a_{k-1,i_{1},\ldots,i_{r}}^{}
\right)^{\frac{b!}{b-(i_{1}+\cdots+i_{r})}}}{x_{k-1,1}^{b!}}
     $$
 it follows that the transform of the couple $(I(\hat{f}_{k-1},b), b!)$
 ($I(\hat{f}_{k-1},b)\subset \hat{\mathcal{O}}_{Z_{k-1},\xi_{k-1}}$) is
 $I(\hat{f}_{k},b)\subset \hat{\mathcal{O}}_{Z_{k},\xi_{k}}$ and the Lemma follows now by (\ref{eccoef}).
\end{proof}
\begin{parrafo}\label{demproptch} {\em Proof of \ref{proptch}.} We first consider a covering
$\{ U_{\lambda}\}_{\lambda \in \Lambda}$ of $W$, so that there is
a closed and smooth hypersurface $ \overline{W}_{\lambda}\subset
U_{\lambda}$, and $I(\overline{W}_{\lambda}) \subset
(\Delta^{b-1}(J))_{\lambda}$ as in \ref{rk710}. After suitable
refinement we may assume that, for each $\lambda$, the conditions
of Lemma \ref{lemita} hold for $Z=\overline{W}_{\lambda}$, and
$J=\langle f_j \rangle$ in \ref{rk710}, where $\{f_1,\dots
,f_j,\dots,f_l\}$ are global sections in $\O_{U_{\lambda}}$ that
span $J_{U_{\lambda}}=\langle f_1,\dots ,f_l \rangle$.

Finally, one can check that a couple $(K_{\lambda}^{(0)}, b!)$,
with property {\bf P3)} in Proposition \ref{proptch}, is defined
by setting:
 $$K_{\lambda}^{(0)}=\langle \left( \frac{1}{i_{1}!\cdots i_{r}!}
\overline{\left( \frac{\partial^{i_{1}+\cdots+i_{r}} f_j}{\partial
z_{1}^{i_{1}}\cdots
      \partial z_{r}^{i_{r}}}\right)} \right)^{\frac{b!}{b-(i_{1}+\cdots+i_{r})}} / 0\leq
i_{1}+\cdots+i_{r}<b \ ; \  j=1,\dots ,l \rangle.$$

\end{parrafo}
\end{section}

\section{On resolution functions I.}\label{On resolution functions.}

\begin{parrafo} In this, and in the next Section \ref{Monomial
Case},
we
show that resolution of basic objects can be achieved once we know
how to define resolution for a simple class of basic objects.
\end{parrafo}
\begin{definition}\label{simplebo}
We will say that a basic object $(W,(J,b),E)$ is a {\em simple}
basic object, if  $(J,b)$ is a simple couple (\ref{nomeacuerdo}),
$J\neq 0(\subset \O_W)$, and $E=\emptyset $ (or, more generally,
if $H_i\cap \Sing(J,b)=\emptyset$ for any $H_i\in E$).

\bigskip

The following result indicates the relevance of simple basic
objects for inductive arguments.
\end{definition}
\begin{proposition}\label{sobreR(1)} Fix a simple basic
object $(W,(J,b),E=\emptyset)$ . Set $\Sing(J,b)=\cup_{1\leq i
\leq s}Z_i$, where each $Z_i$ denotes an irreducible component of
this closed set, and let $R(1)(\Sing(J,b))$ be the union of those
$Z_i$ of codimension one in $W$. Then

a) $R(1)(\Sing(J,b))$ is open and closed in $\Spec(J,b)$ (i.e. a
union of connected components), and it is a closed and smooth
hypersurface in $W$. Moreover, no other component of $\Sing(J,b)$
meets $R(1)(\Sing(J,b))$.

b) If $(W,(J,b),E=\emptyset)\longleftarrow (W_1,(J_1,b),E_1)$ is
defined with center $R(1)(\Sing(J,b)$, then $W_1=W$ and
$\Sing(J_1,b)=\Sing(J,b)-R(1)(\Sing(J,b))$.

In particular $(W_1,(J_1,b),E_1)$ is a simple basic object and
$R(1)(\Sing(J_1,b))=\emptyset$.
\end{proposition}
\proof a) If $Z_1$ is of codimension one, and if $x\in Z_1 \cap
Z_i$ for some other component $Z_i$, then $\Spec(J,b)$ cannot be
included in a smooth hypersuface locally at $x$, in contradiction
with property {\bf P1)} in \ref{lematch}. The same property
insures that $R(1)(\Sing(J,b))$ is regular.

b) The blow-up on a hypersurface is the identity map, so $W_1=W$.
The second assertion follows from property {\bf P2)} in
\ref{lematch}. In fact, locally at a point $x\in R(1)(\Sing(J,b))$
there is a smooth hypersurface $\overline{W}$, such that locally
at $x\in W$, $\Sing(J,b)=\overline{W}$. Moreover, locally at $x\in
W_1=W$, $\Sing(J_1,b)\subset \overline{W}_1$, where
$\overline{W}_1$ is the strict transform of $\overline{W}$ by
blowing up at the center $\overline{W}$. So $\overline{W}_1$, and
hence $\Sing(J_1,b)$, are empty locally at $x$.

\begin{parrafo}{\bf Resolution functions and
the principle of Patching.}\label{patchingfunct} If
$(W,(J,b),E=\emptyset)$ is a  simple basic object, Proposition
\ref{sobreR(1)} says that, after blowing-up at the center
$R(1)\Sing(J,b)$, we may assume that $R(1)\Sing(J,b)=\emptyset$;
and point is that in this setting we can profit of the form of
induction on the dimension $d$ in Proposition \ref{proptch}. In
fact, if the simple basic object is such that
$R(1)\Sing(J,b)=\emptyset$, then there is a covering $\{
U_{\lambda}\}_{\lambda \in \Lambda}$ of $W$, and for each index
$\lambda$ a $d-1$ dimensional basic object
$\overline{B}^{d-1}_{\lambda}=(\overline{W}_{\lambda},(K_{\lambda}^{(0)},
b!), \emptyset)$, such that $\Sing(J,b) \cap
U_{\lambda}=\Sing(K_{\lambda}^{(0)}, b!)$. Note that
$R(1)\Sing(J,b)=\emptyset$ asserts that the ideals
$K_{\lambda}^{(0)}$ are non-zero; a condition required in our
definition of basic object.

Assume that $(I^{d-1},\geq)$ has been defined, together with the
functions defining, as in \ref{resfunct}, resolutions of $d-1$
dimensional basic objects. {\em We will require} that
$$f_{\overline{B}^{d-1}_{\lambda}}=f_{\overline{B}^{d-1}_{\beta}}$$
along points in $\Sing(J,b) \cap U_{\lambda} \cap U_{\beta}$
(condition of patching). In such case we can define a function
$$f^{d-1}_B:\Sing(J,b)\to I^{d-1}$$ simply by patching the
functions $f_{\overline{B}^{d-1}_{\beta}}$. The function
$f^{d-1}_B$ is upper-semi-continuous, and $U_{\delta} \cap \Max
f_B^{d-1} =\Max f_{\overline{B}^{d-1}_{\delta}}$ whenever
$U_{\delta} \cap \Max f_B^{d-1} \neq \emptyset$. Therefore $\Max
f_B^{d-1}$ is a center defining $$(W,(J,b), E=\emptyset )
\longleftarrow (W_1,(J_1,b), E_1 ).$$
 Assume, for simplicity, that
$U_{\delta} \cap \Max f_B^{d-1} \neq \emptyset$ for each index
$\delta$, then $W_1$ can be covered by $\{
U^{(1)}_{\lambda}\}_{\lambda \in \Lambda}$ (notation as in
\ref{cubri}), and for each $\lambda$ we obtain
$$\overline{B}^{d-1}_{\lambda}=(\overline{W}_{\lambda},(K_{\lambda}^{(0)},
b!), \emptyset)\longleftarrow
(\overline{B}^{d-1}_{\lambda})_1=(\overline{W}_{\lambda}^{(1)},(K_{\lambda}^{(1)},
b!), \overline{E}_1)$$ with $\overline{W}_{\lambda}^{(1)} \subset
U^{(1)}_{\lambda}$. Furthermore, the closed set
$\Sing(J_1,b)\subset W_1$ is such that $\Sing(J_1,b) \cap
U_{\lambda}^{(1)}=\Sing(K_{\lambda}^{(1)}, b!)$. {\em We will
require that}
$$f_{(\overline{B}^{d-1}_{\lambda})_1}=f_{(\overline{B}^{d-1}_{\beta})_1}$$
along points in $\Sing(J_1,b) \cap U_{\lambda}^{(1)} \cap
U_{\beta}^{(1)}$ so as to define $f^{d-1}_{B_1}:\Sing(J_1,b)\to
I^{d-1}$(requirement of patching); and in such case $\Max
f^{d-1}_{B_1}$ defines a transformation of $(W_1,(J_1,b),E_1)$.

The point is that if all these requirements of patching hold again
and again, the resolutions of the different basic objects
$\overline{B}^{d-1}_{\lambda}=(\overline{W}_{\lambda},(K_{\lambda}^{(0)},
b!), \emptyset)$, defined in terms of the functions on
$(I^{d-1},\geq)$, patch so as to define a resolution of
$(W,(J,b),\emptyset)$. This would provide resolution of simple
basic objects of dimension $d$.

{\bf Conclusion:} Resolution of simple basic objects
$(W,(J,b),\emptyset)$  can be achieved by blowing up successively
at $\Max f^{d-1}_{B_i}$, for $f^{d-1}_{B_i}:\Sing(J_i,b)\to
I^{d-1}$ defined as above, if the condition of patching holds.

{\bf General strategy for resolution of basic objects:}

{\bf 1)} Define the functions so that the patching principle
holds.

{\bf 2}) Reduce the problem of resolution of a basic object to
that of simple basic objects (\ref{simplebo}).
\end{parrafo}

\begin{parrafo}\label{defdeword} Fix a basic object and a sequence of
transformations
\begin{equation}\label{trabo}
B_0=(W,(J,b),E) \longleftarrow B_1=(W_1,(J_1,b), E_1) \ldots
\longleftarrow B_k=(W_k,(J_k,b), E_k)
\end{equation}
where $E=\{H_1,\dots,H_r\}$ and $E_k=\{H_1,\dots ,H_r,
\dots,H_{r+k}\}$. There is an expression relating $J_k$ with the
total transform, say:
\begin{equation}\label{expressbo}
J \O_{W_k }= I(H_{r+1})^{c_1}I(H_{r+2})^{c_2}\cdots
I(H_{r+k})^{c_k}. J_k
\end{equation}
as in (\ref{express}). Note here that $J_k$ might vanish along
some of the exceptional hypersurfaces $\{H_{r+1},\ldots
,H_{r+k}\}$; and, as indicated in the example in
(\ref{ecuabarra}), there is also a well defined expression:
\begin{equation}\label{expressbo1}
J_k = I(H_{r+1})^{a_{r+1}}I(H_{r+2})^{a_{r+2}}\cdots
I(H_{r+k})^{a_{r+k}}. \overline{J}_k
\end{equation}
at $\O_{W_k }$, so that $\overline{J}_k$ does not vanish along any
exceptional hypersurface. Set $d= dim \ W_k$, and
\begin{equation*}\label{eqord}
\begin{array}{rccl}
\ord^{d}_k: & \Sing(J_k,b) & \to & {\mathbb Q}\\ & x &
\longrightarrow &

\frac{\nu_{J_k}(x)}{b},
\end{array}
\end{equation*}
\begin{equation*}
\begin{array}{rccl}
\word^{d}_k: & \Sing(J_k,b) & \to & {\mathbb Q}\\ & x &
\longrightarrow &

\frac{\nu_{\overline{J}_k}(x)}{b},
\end{array}
\end{equation*}
where \(\nu_{J_k}(x)\) ($\nu_{\overline{J}_k}(x)$) denotes the
order of the ideal \(J{\mathcal O}_{W_k,x}\) (
\(\overline{J}_k{\mathcal O}_{W_k,x}\) ) at \({\mathcal
O}_{W_k,x}\). Both functions are upper semi-continuous and, since
$J_0=\overline{J}_0=J$, they coincide for $k=0$.

We will also define an upper semi-continuous function by setting
$$\overline{a}_{r+i}:  \Sing(J_k,b)  \to  {\mathbb Q}$$

where $\overline{a}_{r+i}(x)=\frac{a_{r+i}(x)}{b}$ and
$a_{r+i}(x)=a_{r+i}$ in (\ref{expressbo1}) if $x\in H_{r+i}$, and
 $a_{r+i}(x)=0$ otherwise.

The role of the denominator $b$ is of no use for the moment and
the reader might want to ignore it. We will justify the presence
of $b$ in \ref{invexp}.

\end{parrafo}
\begin{remark}\label{rk96}
Assume that $\max \word^d_k = \frac{b'}{b}$, and that $
(W_k,(J_k,b), E_k)\longleftarrow (W_{k+1},(J_{k+1},b), E_{k+1})$
is defined with center, say $Y_k\subset \Max \word^d_k$. So the
function $\word^d_k$ takes only the value $\frac{b'}{b}$ along
points of $Y_k$. Then the expression (\ref{expressbo1}),
corresponding now to $J_{k+1}$, is:
\begin{equation}\label{eq77}
J_{k+1} = I(H_{r+1})^{a_1}I(H_{r+2})^{a_2}\cdots I(H_{r+k})^{a_k}
\cdot  I(H_{r+k+1})^{a_{k+1}}\cdot \overline{J}_{k+1}.
\end{equation}
Furthermore:

1) $(\overline{J}_k, b')$ is a simple couple
(\ref{simplecouples}), and $\Max \word_k=\Sing (\overline{J}_k,
b')\cap \Sing(J_k,b)$.

2)  $ (\overline{J}_{k+1}, b')$ is the transform of the simple
couple $(\overline{J}_k, b')$, and $\max \word_{k+1} \leq
\frac{b'}{b}$.

In fact, if $\overline{J}_{k+1}$ is not to vanish along
$H_{r+k+1}$ it must be defined as the proper transform of
$\overline{J}_{k}$ (\ref{deftransid}). The first assertion follows
from our choice of center and the second from Theorem
\ref{theorem72}.
\end{remark}

\begin{parrafo} \label{functk} We will impose conditions on a sequence
(\ref{trabo}). Set $\frac{b'_i}{b}=\max \word_{i}$. Assume all
centers $$Y_i\subset \Max \word_i,\mbox{ for } i=0,\dots, k;$$ and
hence, that
 $$\max \word_{k+1}\geq
\max \word_{k+1}\geq \ldots \geq \max \word_{k+1}$$ (namely
$\frac{b'_0}{b}\geq \frac{b'_1}{b}\geq\ldots \geq
\frac{b'_{k-1}}{b}\geq \frac{b'_k}{b}$) by the previous remark.
Let $k_0$ be the smallest index such that
$$\frac{b'_{k_0}}{b}=\frac{b'_{k_0+1}}{b}= \ldots
=\frac{b'_{k}}{b}.$$ For each index $ k_0 \leq j \leq k$ define a
partition on the set of hypersurfaces in $E_j$, say $ E_j=E_j^-
\cup E_j^+$, where $E_j^-=\{H_1,\ldots, H_r, \ldots H_{r+k_0} \}$
and $E_j^+=\{H_{r+k_0+1},\ldots H_{r+j} \}$. So for $j=k_0$
$E_{k_0}=E_{k_0}^-$; and for $j>k_0$, $E_{j}^-$ consists of the
strict transforms of hypersurfaces in $E_{k_0}$. We finally order
${\mathbb Q}\times {\mathbb N}$ lexicographically, and set
\begin{equation*}
\begin{array}{rccl}
t^d_k: & \Sing(J_k,b) & \to & {\mathbb Q}\times {\mathbb N}\\ & x
& \longrightarrow &

(\word^{d}_k(x),n^d_k(x))
\end{array}
\end{equation*}

\begin{equation*}
       n_{k}^d(\xi)=\left\{
\begin{array}{lll}
     \#\{H\in E_{k}\mid \xi\in H\} & {\rm if} &

\word_{k}^d(\xi)<\frac{b'_k}{b} \\
     \#\{H\in E_{k}^{-}\mid \xi\in H\}
& {\rm if} &
     \word_{k}^d(\xi)=\frac{b'_k}{b}.
      \end{array}
\right.
\end{equation*}

To see that $t^d_k$ is upper-semi-continuous we argue
coordinate-wise: fix integers $m$, $n$, and note that
$$G_{(m,n)}=\{\xi \in \Sing(J_k,b) / t^d_k(\xi)\geq
(\frac{m}{b},n)\}$$ is closed. The statement is clear if
$\frac{m}{b}=\frac{b'_k}{b}$. If not, set $G_{(m,n)}= F_1 \cup
F_2$, for $$F_1=\{\xi  / \word_k(\xi)\geq \frac{m+1}{b}\} \mbox{,
and } F_2=\{\xi / \word_k(\xi)\geq \frac{m}{b};\#\{H\in
E^-_{k}\mid \xi\in H\} \geq n \}.$$

\end{parrafo}

\begin{remark}\label{sobreord}
Because of the lexicographic order on $ {\mathbb Q}\times {\mathbb
N}$ we see that $\max t^d_k=(\max \word_k^d, a)$ for some integer
$ 0\leq a\leq \mbox{ dim } W=d$; and hence that $\Max t_k^d\subset
\Max \word_k^d$. With notation as in \ref{functk}, at most $a$
hypersurface of $E_k^-$ cut at a point of $\Max \word_k^d$, and
$\Max t_k^d$ are the points with this condition.

If a transformation $(W_k,(J_k,b),E_k)\leftarrow
(W_{k+1},(J_{k+1},b),E_{k+1})$, is defined with center $Y_k\subset
\Max t^d_k$, then $ \max \word_k \geq \max \word_{k+1}$
(\ref{rk96},2).

If $ \max \word_k > \max \word_{k+1}$, then $\max t_k^d > \max
t_{k+1}^d$. On the other hand, if  $\max \word_k = \max
\word_{k+1}$, then $E_{k+1}^-$ will consist on the strict
transform of hypersurfaces in $E_k^-$. It is clear that in such
case $\max t^d_k \geq \max t^d_{k+1}$.

\end{remark}

\begin{parrafo}\label{proyecciondebo}{\bf Projections of Basic
Objects.}
 So far we have only considered transformations on basic objects
 defined by monoidal transformations. Set $W_{k+1}=W_k\times A^1$ (the affine line), $W_k \leftarrow
W_{k+1}$ the projection, and define
\begin{equation}\label{defprodebo}
 (W_k,(J_k,b),E_k)\leftarrow
(W_{k+1},(J_{k+1},b),E_{k+1})
\end{equation}
 where $E_{k+1}$ is
the pull-back of hypersurfaces in $E_k$, and $J_{k+1}=J_k\cdot
\O_{W_{k+1}}$. We call this a {\em projection} of basic objects.
Projections will play a key role when proving the patching
conditions discussed in \ref{patchingfunct}.

Note that if a point $x_{k+1}\in W_{k+1}$ maps to a point
$x_{k}\in W_{k}$, then the order of $J_k$ at $\O_{W_k,x_k}$ is the
same as the order of $J_{k+1}$ at $\O_{W_{k+1},x_{k+1}}$.

Here that dimension of $W_{k+1}=dim\ W_k+1$, but ignoring
superscripts, the functions $\word_k$, $ n_k$ and $t_k$ can also
be extended to functions  $\word_{k+1}$, $ n_{k+1}$ and $t_{k+1}$
at $\Sing(J_{k+1},b)$ (pull-back of $\Sing(J_k,b)$). Furthermore,
if $ x_{k+1}\in \Sing(J_{k+1},b)$ maps to  $ x_{k}\in
\Sing(J_{k},b)$, then $\word_{k+1}(x_{k+1})=\word_{k}(x_{k})$,
$n_{k+1}(x_{k+1})=n_{k}(x_{k})$, and
$t_{k+1}(x_{k+1})=t_{k}(x_{k})$.

In other words, given a sequence
\begin{equation}\label{traboconproy}
B_0=(W,(J,b),E) \longleftarrow B_1=(W_1,(J_1,b), E_1) \ldots
\longleftarrow B_{k_0}=(W_{k_0},(J_{k_0},b), E_{k_0})
\end{equation}

where each index $i<k$ $(W_i,(J_i,b),E_i) \leftarrow
(W_{i+1},(J_{i+1},b),E_{i+1})$ is defined by

1) a monoidal transformation with center $Y_i \subset \Max
\word_i$, or by

2) a projection of basic objects;

we have that $$\max \word_0 \geq \max \word_1 \geq \dots \max
\word_k.$$ In particular, the partitions of $E_i=E_i^+\cup
E_{i+1}^-$, and the functions $t^d_i: \Sing(J_k,b)  \to {\mathbb
Q}\times {\mathbb N}$ defined in \ref{functk} can also be defined
in this setting. Furthermore, if all centers $Y_i\subset \Max t_i
(\subset \Max \word_i)$ then
 $$\max t_0 \geq \max t_1 \geq \dots \max t_{k_0}.$$

\end{parrafo}
\begin{proposition}\label{prop7.9}
Consider a sequence (\ref{traboconproy}) (of transformations and
projections), and assume that either $k_0=0$ or that $ \max t^d_0
\geq \max t^d_1\ldots \geq \max t^d_{k_0-1}> \max t^d_{k_0}$. Then
there is a simple basic object $(W_{k_0},(K_{k_0},c),\emptyset)$
(Def. \ref{simplebo}) such that any sequence
\begin{equation}\label{tratkk}
(W_{k_0},(K_{k_0},c),\emptyset) \longleftarrow
(W_{k_0},(K_{k_0+1},c),E'_1) \longleftarrow \ldots
(W_{k_0},(K_{k_0+m},c),E'_{m})
\end{equation}
(of transformations and projections) induces a sequence over
$(W_{k_0},(J_{k_0},b),E_{k_0})$, say:
\begin{equation}\label{tratkkk}
(W,(J,b),E) \ldots \longleftarrow
(W_{k_0},(J_{k_0},b),E_{k_0})\longleftarrow \ldots \longleftarrow
(W_{k_0+m},(J_{k_0+m},b),E_{k_0+m}).
\end{equation}
Furthermore, (\ref{tratkk}) and (\ref{tratkkk}) are related by the
following properties:

{\bf P1)} $\max t^d_{k_0}=\ldots=\max t^d_{k_0+m-1}$ in
(\ref{tratkkk}).

{\bf P2)} $\Max t^d_{k_0+j}=\Sing(K_{k_0+j},c)$ for $ 0\leq j \leq
m-1$.

{\bf P3)} If $\Sing(K_{k_0+m},c)=\emptyset$, then
$\Sing(J_{k_0+m},b)=\emptyset$ or $\max t^d_{k_0+m-1} >\max
t^d_{k_0+m}.$

{\bf P4)} If $\Sing(K_{k_0+m}, c)\neq \emptyset$, then
 $\max t^d_{k_0+m-1} =\max
t^d_{k_0+m}$ and $\Sing(K_{k_0+m},c)=\Max t_{k_0+m}$.
\end{proposition}
We begin by the following two remarks, needed to sketch a proof of
this Proposition.

\begin{remark}\label{invexp}
Fix a basic object $B=(W,(J,b),E)$ and a positive integer $m \geq
1$. Set $J^m \subset \O_W$ the $m$-th power of $J$, and consider
$B_m=(W,(J^m,m\cdot b),E)$. Note that $$\Sing(J,b)=\Sing(J^m,m
\cdot b).$$ In particular, a smooth center $Y$ defines a
transformation of one basic object iff it defines a transformation
of both, say:
$$(W,(J,b),E)\longleftarrow (W_1,(J_1,b),E_1) \mbox{ and
}(W,(J^m,m\cdot b),E)\longleftarrow (W_1,((J^m)_1,m\cdot
b),E_1)\mbox{ respectively}.
$$ Since the total transform of $J^m$, namely $J^m.\O_{W_1}$, is
the $m$-th power of the total transform of $J$, it follows that
$(J^m)_1$ is the $m$-th power of $J_1$ (i.e. $(J^m)_1= J_1^m $).
The same holds after any sequence of transformations. Therefore a
resolution of $B$ induces a resolution of $B_m$, and the other way
around. It will turn out that the resolution of $B$, defined by
the resolution functions in \ref{resfunct}, will coincide with the
resolution of $B_m$ defined by the resolution functions. For the
time being note that at a point $\xi \in \Sing(J,b)=\Sing(J^m,m
\cdot b)$, $\frac{\nu_{{J}}(\xi)}{b}=\frac{\nu_{{J}^m}(\xi)}{m
\cdot b}$, where $\nu_J(x)$ denotes the order of $J$ at
$\O_{W,x}$.
\end{remark}
\begin{remark}\label{intersec}
Given two basic objects $(W,(J,b), E)$ and $(W,(K,c),E)$ (same
(W,E)), note that $$\Sing(J,b) \cap \Sing(K,c)=\Sing(J^c,c\cdot
b)\cap \Sing(K^b,b\cdot c)=\Sing(J^c+ K^b,b \cdot c).$$ If
$Y\subset \Sing(J^c+ K^b,b\cdot c)$ defines $(W,(J^c+ K^b,b\cdot
c) ,E) \leftarrow (W_1,((J^c+ K^b)_1,b\cdot c) ,E_1)$, then $Y$
also defines $(W,(J,b), E)\leftarrow (W_1,(J_1,b), E_1)$, and
$(W,(K,c), E)\leftarrow (W_1,(K_1,c), E_1)$. Check finally that
$(J^c+ K^b)_1=J_1^c+K_1^b$, and hence that
$$\Sing(J_1,b) \cap \Sing(K_1,c)=\Sing((J^c+ K^b)_1,c\cdot b).$$
So, by induction, a sequence of transformations of $(W,(J^c+
K^b,b\cdot c) ,E)$,  induces sequences of transformations of
$(W,(J,b), E)$ and of $(W,(K,c), E)$, and for each index $i$,
$$\Sing(J_i,b) \cap \Sing(K_i,c)=\Sing((J^c+ K^b)_i,c\cdot b).$$
Because of this property we will set formally:
\begin{equation}\label{forinter}
(W,(J,b), E)\cap (W,(K,c), E)=(W,(J^c+ K^b,b \cdot c) ,E)
\end{equation}
\end{remark}
\begin{example} Let $X_0 $ be an
 curve in a smooth surface $W_0$, analytically irreducible at a closed
 point $\xi_0 \in W_0$. These data allow us to define, for
each integer $k$, a sequence of $k$ quadratic transformations over
$W_0$. In fact, if $W_0 \leftarrow W_1$ is defined with center
$\xi_0$, the strict transform $X_1$ intersects the exceptional
locus $H_1$ at a unique point, say $\xi_1$. Set $W_1 \leftarrow
W_{2}$ with center $\xi_1$. By iteration we get $W_i \leftarrow
W_{i+1}$, with exceptional hypersurface $H_{i+1}$, and $\xi_{i+1}=
H_{i+1}\cap X_{i+1} \in W_{i+1}$.

Set $X_0=V(<x^4- y^5>)\subset W_0=\Spec(k[x,y])$. For any $k$, the
sequence of length $k$ defined by this curve induces a sequence of
transformations of  $(W_0,J=<x^4- y^5>,1), E_0=\emptyset)$. For
$k=1$: $$ (W_0,J=<x^4- y^5>,1), E_0=\emptyset)\leftarrow
(W_1,(J_1,1), E_1).$$ Check first that
$J_1=I(H_1)^3\overline{J}_1$, and that $ \max t^2_0=(4,0)> \max
t^2_1= (1,1)$. So for $k=1$, $k_0=1$ in the setting of Proposition
\ref{prop7.9}. Show that $(W_1,(\overline{J}_1,1),\emptyset)\cap
(W_1,(I(H_1),1),\emptyset)$ (see (\ref{forinter})) plays the role
of $(W_{k_0},(K_{k_0},c),\emptyset)$ in the Proposition. Note
finally that for the sequence of five quadratic transformations
defined by the curve (i.e. for $k=5$): $$ \max t^2_1= \max t^2_2=
\max t^2_3= \max t^2_4> \max t^2_5.$$
\end{example}
\begin{proof}{\em Of Prop \ref{prop7.9}.} We define
$(W_{k_0},(K_{k_0},c),\emptyset)$ with those properties, and we do
so by taking suitable intersections (\ref{intersec}).

If $\max t^d_{k_0}=(\frac{b_{k_0}}{b},n_{k_0})$, $b_{k_0}$ is the
highest order of $\overline{J}_{k_0}$ along points in
$\Sing(J_{k_0})$. Set
$$(W_{k_0},(A_{k_0},c),\emptyset)=(W_{k_0},(J_{k_0},b),\emptyset)\cap
(W_{k_0},(\overline{J}_{k_0},b_{k_0}),\emptyset)$$ so that
$\Sing(A_{k_0},c)$ is $\Max \word_{k_0}$. By assumption
$E_{k_0}=E^{-}_{k_0}$, and at most $n_{k_0}$ hypersurfaces of
$E^{-}_{k_0}$ can come together at a point of $\Max \word_{k_0}$.
For a subset $S \subset E^{-}_{k_0}$, with $n_{k_0}$
hypersurfaces, set $F_S=\Max \word \cap (\cap_{H_i \in S} H_i)$.
Given two such subsets $S_1\neq S_2$ note that $F_{S_1}\cap
F_{S_2}=\emptyset$ since $n_{k_0}$ is a maximum. Furthermore $\Max
t_{k_0}=\cup F_S$ for all $S$ as before. Recall that each $H_i \in
E_{k_0}$ is a smooth hypersurface, and set
$(W_{k_0},(I(H_{i}),1),\emptyset)\cap(W_{k_0},(I(H_{j}),1),\emptyset)(=(W_{k_0},(I(H_{i})+I(H_j),1),\emptyset)).$
Finally set $B_{k_0}=\sum_{S} \sum_{H_{j_i}\in S}I(H_{j_i})$, and
$$(W_{k_0},(K_{k_0},c),\emptyset)=(W_{k_0},(A_{k_0},c),\emptyset)\cap
(W_{k_0},(B_{k_0},1),\emptyset),$$ and check that
$\Sing(K_{k_0},c)=\Max t^d_{k_0}$.

If $Y_{k_0} \subset \Sing(K_{k_0},c)$ is a center of
transformation for this basic object, then, for any $H_i \in
E_{k_o}$, either $Y_{k_0} \subset H_i$ or $Y_{k_0} \cap H_i=
\emptyset$. In particular $Y_{k_0}$ has normal crossing with
$E_{k_0}$ and defines a transformation of
$(W_{k_0},(J_{k_0},b),E_{k_0})$. Furthermore, using
\ref{weakvsstrict} and the previous Remarks, we conclude that
either $\max t^d_{k_0}>\max t^d_{k_0+1}$, in which case
$\Sing(K_{k_0+1},c)=\emptyset,$ or $\max t^d_{k_0}=\max
t^d_{k_0+1}$, in which case $\Sing(K_{k_0+1},c)=\Max t^d_{k_0+1}$
(notation as in (\ref{tratkk}) and (\ref{tratkkk})). \

In this last case $E'_1$ (in  (\ref{tratkk})) is $E_{k_0+1}^+$ in
(\ref{tratkkk}). If $Y_{k_0+1}$ is a center that defines a
transformation of $(W_{k_0+1},(K_{k_0+1},c),E^{'}_1)$, then
$Y_{k_0+1}$ must have normal crossing with $E_{k_0+1}^+$, and on
the other hand, for any hypersurface $H_i \in E_{k_0+1}^-$ either
$Y_{k_0+1}\subset H_i$ or $Y_{k_0+1}\cap H_i=\emptyset$. This
insures that $Y_{k_0+1}$ has normal crossing with $E_{k_0+1}$, and
defines a transformation of $(W_{k_0+1},(J_{k_0+1},b),E_{k_0+1})$.

All properties in  Proposition \ref{prop7.9} follow by iteration
of this argument.

We end this proof by showing that
$(W_{k_0},(K_{k_0},c),\emptyset)$ is a simple basic object. To
check this note that if $J$ has highest order $b$, then $J^c+ K^b$
has highest order $b \cdot c$ in (\ref{forinter}). So it suffices
to check that $B_{k_0}$ has highest order 1, which is clear.
\end{proof}



\section{On resolution functions II: the Monomial Case}\label{Monomial Case}
\begin{parrafo}\label{com101}
Proposition \ref{prop7.9} asserts that if we knew how to define
resolutions of simple basic objects, then we could define an
extension (\ref{tratkkk}), so that $\max t^d_{k_0}=\ldots=\max
t^d_{k_0+m-1}>\max t^d_{k_0+m}.$ Since (\ref{tratkkk}) is a
sequence of transformations of $(W,(J,b), E)$, the first
coordinate of $\max t^d$ is of the form $\frac{l}{b}$ for a
positive integer $l$; and the second coordinates is at most the
dimension of $W$. So by iteration of resolutions of the simple
basic objects in Proposition \ref{prop7.9}, we could force the
last value $\max t^d$ to drop again and again. Ultimately, we come
to the case in which either (\ref{tratkkk}) is a resolution, or
the first coordinate of $\max t^d_{k_0+m}$ is zero. This last case
is called {\it the monomial case}. In Proposition
\ref{prop7.9casoword0} we will provide resolution for this case.
Note also that in the monomial case $\overline{J}_k=\O_{W_k}$ in
(\ref{expressbo1}), locally at any point in $\Sing(J_k,b)$, so
$J_k$ is locally spanned by a monomial.
\end{parrafo}
\begin{example} \label{exmonomial} Let $W=\Spec(k[X_1,X_2,X_3])$
 and $(W,(J,5), \emptyset)$, $J=<X_1^6\cdot X_2^7 \cdot X_3^4>$. The singular
locus is a union of two hypersurfaces $V(<X_1>)\cup V(<X_2>)$.
Blowing up $V(<X_2>)$ we get $W_1=W$ and $(J_1,5)$, $J_1=<X_1^6
\cdot X_2^2 \cdot X_3^4>$. The singular locus is a union of a
hypersurface with a line. Blowing up at the hypersurface we get
$W_2=W$ and $J_2=<X_1^1\cdot X_2^2 \cdot X_3^4>$ where the
singular locus is a line. A resolution is finally achieved by
blowing up such line.

It is simple to establish a general strategy, in the monomial
case, so that, as in this example, resolution is achieved by
blowing up at maximal dimension components of the singular locus.
\end{example}

Note that, for a monomial basic object, the closed set \(
\Sing(J,b) \) is the  union of some of the irreducible components
of intersections of the hypersurfaces \( H_{i} \).  In fact,
consider the functions \( a_{i_{1}},\ldots, a_{i_{p}} \) defined
in \ref{defdeword}, and an irreducible component \( C \) of the
intersection \( H_{i_{1}}\cap\cdots\cap H_{i_{p}} \); then the
functions \( a_{i_{1}},\ldots, a_{i_{p}} \) are constant on \( C
\), and \( C \) is included in \( \Sing(J,b) \) if and only \(
a_{i_{1}}+\cdots+a_{i_{p}}\geq b \) along \( C \).

\begin{definition} \label{DefGamma}
Let \( (W,(J,b),E) \) be a monomial basic object. Define the
function:
\begin{gather*}
h:\Sing(J,b)\longrightarrow
\Gamma=\mathbb{Z}\times\mathbb{Q}\times\mathbb{Z}^\mathbb{N} \\
h(\xi)=
     (-p(\xi),\omega(\xi),\ell(\xi)).
\end{gather*}
where, if $\xi \in \Sing(J,b)$, the values $p(\xi)$, $\omega(\xi)$
and $\ell(\xi)$ are defined as follows:

\begin{equation} \label{DefGam1}
p(\xi)=
       \min\left\{
        q\mid\exists i_1,\ldots,i_q,\
\begin{array}{l}
            a_{i_1}(\xi)+\cdots+a_{i_q}(\xi)\geq b  \\
\xi\in H_{i_1}\cap\cdots\cap H_{i_q}
\end{array}\right\},
\end{equation}

\begin{multline} \label{DefGam2}
\omega(\xi)=
       \max\left\{
\frac{a_{i_1}(\xi)+\cdots+a_{i_q}(\xi)}{b}\mid \right.\\
       \left.
q=p(\xi), \mbox{ } a_{i_1}(\xi)+\cdots+a_{i_q}(\xi)\geq b,\mbox{ }
\xi\in H_{i_1}\cap\cdots\cap H_{i_q}\right\},
\end{multline}
and
\begin{multline} \label{DefGam3}
\ell(\xi)=
       \max\left\{(i_1,\ldots,i_q,0,\ldots)\in \mathbb{Z}^{\mathbb{N}}\mid\right. \\
\left. q=p(\xi), \mbox{ }
\frac{a_{i_1}(\xi)+\cdots+a_{i_q}(\xi)}{b}=\omega(\xi),
       \mbox{
}\xi\in H_{i_1}\cap\cdots\cap H_{i_q}\right\}.
\end{multline}
In the last formula we consider the lexicographical order in \(
\mathbb{Z}^{\mathbb{N}} \).
\end{definition}

Fix a point \( \xi\in\Sing(J,b) \) and let \( C_{1},\ldots,C_{s}
\) be the irreducible components of \( \Sing(J,b) \) at \( \xi \).
\bigskip

$\bullet$ The first coordinate of \( h(\xi) \) is \( -p(\xi) \),
where \( p(\xi) \) is the minimal codimension of \(
C_{1},\ldots,C_{s} \).

Denote by \( C'_{1},\ldots,C'_{s'} \) the components with minimum
codimension \( p(\xi) \) (i.e. of highest dimension at the point
$\xi$).

$\bullet$ The second coordinate of \( h(\xi) \) is \(
\omega(\xi)=\dfrac{b'}{b} \), where \( b' \) is the maximum order
of \( J \) along the components \( C'_{1},\ldots,C'_{s'} \).

Denote by \( C''_{1},\ldots,C''_{s''} \) the components with
maximum order.

$\bullet$ The last coordinate of \( h(\xi) \), \( \ell(\xi) \),
corresponds to one \( C''_{j} \), for some index \( j \).

So for a fixed point \( \xi \), with \( p(\xi) \) we have selected
the irreducible components of \( \Sing(J,b) \), at \( \xi \), of
highest dimension.  With \( \omega(\xi) \) we have select, among
the previous components, those where the order of \( J \) is
maximum. Finally with \( \ell(\xi) \) we select a unique component
containing \( \xi \).

\begin{parrafo}\label{monoparrafo1}
Now one can check that the function \( h \) is
upper-semi-continuous, and that the closed set \( \Max{h} \) is
regular. In fact if \(
\max{h}=(-p_{0},w_{0},(i_{1},\ldots,i_{p_{0}},0,\ldots)) \), then
\( \Max{h} \) is a union of connected components of the regular
scheme \( H_{i_{1}}\cap\cdots\cap H_{i_{p_{0}}} \).

It is clear that \( \Max{h} \) is a permissible center for the
basic object \( (W,(J,b),E) \). Let
\begin{equation*}
(W,(J,b),E)\stackrel{\Pi}{\longleftarrow}(W_{1},(J_{1},b),E_{1})
\end{equation*}
be the transformation with center \( \Max{h} \), and let \(
E_{1}=\{H_{1},\ldots,H_{r},H_{r+1}\} \), where, by abuse of
notation, \( H_{i} \) denotes the strict transform of \( H_{i} \),
for \( i=1,\ldots,r,\) and \( H_{r+1} \) is the exceptional
divisor of \( \Pi \). The basic object \( (W_{1},(J_{1},b),E_{1})
\) is also monomial, in fact for \( \xi\in\Sing(J_{1},b) \) we
have
\begin{equation}
       J_{\xi}=
\mathcal{I}(H_{1})^{a'_{1}(\xi)}_{\xi}\cdots
\mathcal{I}(H_{r})^{a'_{r}(\xi)}_{\xi}
\mathcal{I}(H_{r+1})^{a'_{r+1}(\xi)}_{\xi}, \label{eq:ExpJ1}
\end{equation}
where the functions \( a'_{i} \) are given by:
\begin{equation}
       \begin{array}{lll}

a'_{i}(\xi)=a_{i}(\Pi(\xi)) & \forall \xi\in H_{i} &

\text{and}\ i=1,\ldots,r; \\

a'_{r+1}(\xi)=a_{i_{1}}(\Pi(\xi))+\cdots+a_{i_{p_{0}}}(\Pi(\xi))-b
&

\forall \xi\in H_{r+1}. \\
        \end{array}
\label{eq:Expaprim}
 \end{equation}

As in Definition \ref{DefGamma},  a function \( h_{1} \) has been
associated to the basic object \( (W_{1},(J_{1},b),E_{1}) \), and
one can check that the maximum value has dropped:
 \begin{equation*}
 \max{h}>\max h_{1}.
\end{equation*}

In fact, for any point \( \xi\in\Sing(J_{1},b) \):
\begin{equation}
       \begin{array}{lll}

h_{1}(\xi)=h(\Pi(\xi)) & \text{if} &
       \Pi(\xi)\not\in\Max h  \\

h_{1}(\xi)<h(\Pi(\xi)) & \text{if} & \Pi(\xi)\in\Max h.

\end{array}
         \label{eq:IgualMenor}
 \end{equation}
 It is not hard
to check now that this function  \( h \) defines
 a resolution in the monomial case:
 \end{parrafo}
 \begin{proposition}\label{prop7.9casoword0}
Consider a sequence (\ref{traboconproy}) (of transformations and
projections), and assume that $ \max t^d_0 \geq \max t^d_1\ldots
\geq \max t^d_{k_0-1}> \max t^d_{k_0}$, and that $\max
\word_{k_0}=0$. A resolution
\begin{equation}\label{tratkkk00}
 (W_{k_0},(J_{k_0},b),E_{k_0})\longleftarrow \ldots
\longleftarrow (W_{k_0+m},(J_{k_0+m},b),E_{k_0+m}).
\end{equation}
is defined by the functions $h_i:\Sing(J_{k_0+i},b) \to \Gamma$,
by blowing up successively at $Y_{k_0+i}=\Max h_i$.
\end{proposition}



\section{General basic objects and resolution functions.} \label{sectiongbo} In \ref{patchingfunct}
we already discussed the need to generalize the notion of basic
object in order to profit from a form of induction on the
dimension of basic objects, which would enable us to achieve
resolutions of basic objects. This leads us to the notion of {\em
general basic
  objects} which will be
developed in this section. Recall that in the setting of
\ref{patchingfunct}, namely the case of a simple basic object
$(W,(J,b),E=\emptyset)$ (in which dim $W=d$, and where
$R(1)(\Sing(J,b))=\emptyset$), there is a form of induction on the
dimension $d$. In fact, in such case there is a covering $\{
U_{\lambda}\}_{\lambda \in \Lambda}$ of $W$, and for each index
$\lambda$ a $d-1$ dimensional basic object
$\overline{B}^{d-1}_{\lambda}=(\overline{W}_{\lambda},(K_{\lambda}^{(0)},
b!), \emptyset)$, such that $\Sing(J,b) \cap
U_{\lambda}=\Sing(K_{\lambda}^{(0)}, b!)$. The outcome of the
previous sections \ref{On resolution functions.} and \ref{Monomial
Case} is to show that resolutions of simple basic objects implies
resolutions of arbitrary basic objects. However in doing so, we
expect to argue inductively by defining the functions $\word$, $n$
(see (\ref{functk})), and $h$ (see (\ref{DefGamma})), for these
{\it locally defined} basic objects
$\overline{B}^{d-1}_{\lambda}$. In this Section we provide a
precise formulation of these locally defined basic objects. The
key point, that will ultimately allow us to define the functions
$\word$, $n$ , and $h$ in this more ample context, is the fact
that the singular loci of these d-1 dimensional basic objects,
namely the sets $\Sing(K_{\lambda}^{(0)}, b!) $, patch and define
the closed set $\Sing(J,b)$. In fact $\Sing(J,b) \cap
U_{\lambda}=\Sing(K_{\lambda}^{(0)}, b!)$. Furthermore, this form
of patching will also hold after transformations; a concept that
will be made precise in the following definition. The covering $\{
U_{\lambda}\}_{\lambda \in \Lambda}$ of $W$ and the $d-1$
dimensional basic object
$\overline{B}^{d-1}_{\lambda}=(\overline{W}_{\lambda},(K_{\lambda}^{(0)},
b!), \emptyset)$ will define, in the sense of the following
definition, a $d-1$ dimensional {\em general basic object}.

\begin{definition}
\label{defbo} A {\em \(d\)-dimensional general basic object over a
  pair \((W,E)\)} ($W$ smooth, $E=\{H_1,\dots ,H_s\}$ as in (\ref{normalc1})),
   consists of an open covering of \( W
\), say \(\{U_{\alpha}\}_{\alpha \in \Lambda} \); and setting
\((U_{\alpha},E_{\alpha})\) as the restriction of \((W,E)\) to
\(U_{\alpha}\), there is:
\begin{enumerate}
\item[(i)] {\bf A collection of basic objects.} For
  every \(\alpha \in \Lambda\) there is a
closed and smooth
  \(d\)-dimensional subscheme \(\widetilde{W}_{\alpha}
\subset
  U_{\alpha}\),
 which intersects transversally all hypersurfaces \(E_{\alpha}\), in the sense that
 $H_i\cap \widetilde{W}_{\alpha}=(H_{\alpha})_i$ is either empty
  or a smooth hypersurface
of $\widetilde{W}_{\alpha}$, defining a pair
\((\widetilde{W}_{\alpha} , \widetilde{E}_{\alpha})=\{
(H_{\alpha})_1,\dots ,(H_{\alpha})_s\}\). And, for each $\alpha$
there is a basic object
\begin{equation*}
(\widetilde{W}_{\alpha},
(B_{\alpha},d_{\alpha}),\widetilde{E}_{\alpha}).
\end{equation*}
Obviously, for each index \(\alpha\) the closed set
$\Sing(B_{\alpha},d_{\alpha})\subset U_{\alpha}$ is locally closed
in  \(W\).

\item[(ii)] {\bf A patching condition.} There is a closed subset
\(F \subset W\) such that
\begin{equation*}
F \cap U_{\alpha}=\Sing(B_{\alpha},d_{\alpha})
\end{equation*}
for every \(\alpha\in \Lambda\).

\item[(iii)] {\bf Stability of patching (I).} Let
\begin{equation*}
(W,E) \longleftarrow (W_1,E_1)
\end{equation*}
be  a permissible transformation with center \(Y\subset F\)
(\ref{normalc1}), let \(\{U_{\alpha,_1}\}\) be  the pullback of
\(\{U_{\alpha}\}_{\alpha\in \Lambda}\)  to  \(W_1\), and for each
\(\alpha\in  \Lambda\) let
\begin{equation*}
(\widetilde{W}_{\alpha}, (B_{\alpha},d_{\alpha}),
\widetilde{E}_{\alpha})\longleftarrow (\widetilde{W}_{\alpha,1},
(B_{\alpha,1},d_{\alpha}), \widetilde{E}_{\alpha,1}).
\end{equation*}
be the corresponding transformation of basic objects. Then there
is a closed set \(F_1 \subset W_1\) so that
\begin{equation*}
F_1 \cap U_{\alpha,1}=\Sing(B_{\alpha,1},d_{\alpha})
\end{equation*}
for each index \(\alpha \in
 \Lambda\).

\item[(iv)] {\bf Stability of patching (II).} Let \(W
\longleftarrow W_1=W \times A^1\) be the projection and let
\begin{equation*}
(W,E) \longleftarrow (W_1,E_1)
\end{equation*}
where $E_1$ is defined as the set of pull-backs of hypersurfaces
in $E$. Let \(\{U_{\alpha,1}\}\) be the pullback of
\(\{U_{\alpha}\}_{\alpha\in I}\) to  \(W_1\), and for each
\(\alpha\in \Lambda\) set
\begin{equation*}
(\widetilde{W}_{\alpha}, (B_{\alpha},d_{\alpha}),
\widetilde{E}_{\alpha})\longleftarrow (\widetilde{W}_{\alpha,1},
(B_{\alpha,1},d_{\alpha}), \widetilde{E}_{\alpha,1}),
\end{equation*}
where $\widetilde{W}_{\alpha,1}=\widetilde{W}_{\alpha}\times A^1$,
$\widetilde{E}_{\alpha,1}$ is the pull-back of hypersurfaces in
$\widetilde{E}_{\alpha}$, and
$B_{\alpha,1}=B_{\alpha}\O_{\widetilde{E}_{\alpha,1}}$. Then there
is a closed set \(F_1 \subset W_1\) such that, for each index
\(\alpha \in  \Lambda\)
\begin{equation*}
F_1 \cap U_{\alpha,1}=\Sing(B_{\alpha,1},d_{\alpha}).
\end{equation*}

\item[(v)] {\bf Stability of patching (III). } The   patching
condition defined in (iii) and (iv)   holds after any sequence of
transformations: Given a sequence of transformations of pairs, $$
\begin{array}{ccccccccc}
(W_0,E_0) &\longleftarrow & (W_1,E_1) & \longleftarrow & \ldots &
\longleftarrow & (W_r,E_r) & \longleftarrow  & (W_{r+1},E_{r+1})
\\ \cup & & \cup & & & & \cup & & \\ F_0 & & F_1 & & & & F_r &
\end{array}
$$ where for $i=0,1,\ldots,r$,  $W_{i+1}\to W_i$ is defined either
by:
\begin{enumerate}
\item[(1)] blowing up  at centers \(Y_i \), permissible for the
pair \((W_i,E_i)\), and \(Y_i\) included in the inductively
defined closed sets \(F_i \subset W_i\), or \item[(2)] a
projection \(p: W_{i+1} \to W_i\),
\end{enumerate}

\noindent there is an open covering \(\{U_{\alpha,r+1}\}\) of
\(W_{r+1}\) (the pull back of \(\{U_{\alpha}\}\)), a sequence of
transformations of basic objects,
\begin{multline} \label{transalpha}
(\widetilde{W}_{\alpha}, (B_{\alpha},d_{\alpha}),
\widetilde{E}_{\alpha})\longleftarrow (\widetilde{W}_{\alpha,1},
(B_{\alpha,1},d_{\alpha}),\widetilde{E}_{\alpha,1 })
\longleftarrow\cdots \\ \cdots\longleftarrow
(\widetilde{W}_{\alpha,r+1},(B_{\alpha,r+1},d_{\alpha}) ,
\widetilde{E}_{\alpha,r+1}),
\end{multline}
   and a closed
set $F_{r+1}\subset W_{r+1}$,  such that for each $\alpha\in
\Lambda$,

\begin{equation*}
F_{r+1} \cap U_{\alpha,r+1}=\Sing(B_{\alpha,r+1},d_{\alpha}).
\end{equation*}

\item[(vi) ] \textbf{Restriction to open sets.} If \( V\subset W
\) is an open set, consider the restriction of all data to \( V
\): the open covering \( \{U_{\alpha}\cap V\}_{\alpha\in\Lambda}
\), the basic objects \( (\widetilde{W}_{\alpha},
(B_{\alpha},d_{\alpha}),\widetilde{E}_{\alpha})_ {V} \),

 and the closed set \( F_{V}=F\cap V \). Then
we require that all properties (i), (ii), (iii) (iv) and (v) hold
for the restriction.
\end{enumerate}
Last condition (vi) could be avoided if we assume
desingularization. In fact, if $Y \subset F\cap V$ is a smooth
center, the closure of $Y$ in $W$ might be singular. If we assume
desingularization we may assume that the closure is regular, and
that the transformation over $V$ is a restriction of a
transformation over $W$. However we want to prove
desingularization, so we impose condition (vi).

 A general basic object will be denoted by
\(({\mathcal F},(W,E))\), the restriction to an open set \( V \)
will be
 denoted by $(\mathcal{F}_{V},(V,E_{V})).$ Note that we have defined two
notions of transformations of general basic objects: one as in
\ref{defbo}(iii) (by a monoidal transformations), and another one
as in \ref{defbo}(iv), by a projection. This last transformation
increases the dimension by one.

We denote a sequence (of transformations and projections) as
\begin{equation}
\label{trangbo}
\begin{array}{ccccccc}
({{\mathcal F}}_0,(W_0,E_0)) & \longleftarrow & \ldots &
\longleftarrow & ({{\mathcal F}}_{r},(W_{r},E_{r})) &
\longleftarrow & ({{\mathcal F}}_{r+1},(W_{r+1},E_{r+1}))\\ \cup &
& & & \cup & & \cup\\ F_0 & & & & F_{r} & & F_{r+1}.
\end{array}
\end{equation}
\end{definition}

\begin{remark} \label{notaopen}
If \( (\mathcal{F},(W,E)) \) is \( d \)-dimensional, then  \( d \)
can be different from \( \dim{W} \).

1) A basic object \( (W,(J,b),E) \) defines a $d$-dimensional
general basic object \((\mathcal{F},(W,E)) \), with the trivial
open covering and \( d=\dim{W} \).

2) A simple basic object $(W,(J,b),E=\emptyset)$, with dim $W=d$
and $R(1)(\Sing(J,b))=\emptyset $, also defines a general basic
object \((\mathcal{F},(W,E)) \) of dimension $d-1$.

This follows from Proposition \ref{proptch}, as was indicated in
\ref{patchingfunct}; see also \ref{rk715} for the case of
projections.
\end{remark}

\begin{remark} \label{onsgbo}

A general basic object can be described  by giving two different
open coverings. What is important here are the closed sets \(F\)
that it defines. That is why in the notation for general basic
objects \(({\mathcal F},(W,E))\) there is no reference to the open
covering which appears in the definition.

A general basic object \(({{\mathcal F}},(W,E))\), defined
 in terms of an open cover $\{U_{\alpha}\}$ of $W$ and basic
 objects \(
(\widetilde{W}_{\alpha},(B_{\alpha},d_{\alpha}),\widetilde{E}_{\alpha})\),
is said to be a {\em simple} general basic object, when all the
basic objects \(
(\widetilde{W}_{\alpha},(B_{\alpha},d_{\alpha}),\widetilde{E}_{\alpha})\)
are simple (\ref{simplebo}).

We now extend the result in Proposition \ref{sobreR(1)} to the
case of general basic objects.

Let $R(1)(F)$ be the union of $d-1$ dimensional components of $F$
(so that $R(1)(F)\cap
U_{\alpha}=R(1)(\Sing((B_{\alpha},d_{\alpha})))).$

a) $R(1)(F)$ is open and closed in $F$ (i.e. a union of connected
components), and smooth in $W$.

b) Setting $({{\mathcal F}},(W,E)) \longleftarrow ({{\mathcal
F}}_1,(W_1,E_1)) $ with center $R(1)(F)$, then $W_1=W$ and $F_1= F
-R(1)(F)$; in particular:

 c) $ ({{\mathcal F}}_1,(W_1,E_1))$  is simple and
$R(1)(F_1)=\emptyset$.

Finally, one can generalize  \ref{notaopen}, 2), to show that if
c) holds, then $ ({{\mathcal F}}_1,(W_1,E_1))$ has a structure of
$d-1$ dimensional general basic object (where $d=\mbox{ dimension
of } ({{\mathcal F}},(W,E))$).

\end{remark}

\begin{definition}\label{resgbo}
A resolution of a general basic object \(({{\mathcal
F}}_0,(W_0,E_0))\) is a sequence of transformations as in
(\ref{trangbo}) which fulfills the following two conditions:
\begin{enumerate}
\item[(i)] The sequence involves only monoidal transformations
(\ref{defbo},(iii)). \item[(ii)] The closed set \(F_{r+1}\) is
empty.
\end{enumerate}

Note that if $\{U_{\alpha}\}$ is an open covering of $W$ as in
Definition \ref{defbo}, then for any \( \alpha \)  we obtain a
resolution of the basic object \(
(\widetilde{W}_{\alpha},(B_{\alpha},d_{\alpha}),\widetilde{E}_{\alpha})\)
as defined in \ref{normalc1}.
\end{definition}

\begin{parrafo} \label{compopen}
We will assign, to each general basic object \(({\mathcal
F},(W,E))\), an upper semi-continuous function $f_{\mathcal{F}}:
F\to (T, \geq) $ (on the closed set $F\subset W$ as in
\ref{defbo}, (ii)). Such functions will be defined so that they
are compatible with open restrictions. In other words, if $V$ is
an open subset of $W$, the closed set of the restriction \(
(\mathcal{F}_{V},(V,E_{V})) \) is $F\cap V$, and we require that
the restriction of $f_{\mathcal{F}}$ to $F\cap V$ be
$f_{\mathcal{F}_V}$. The following is an example.
\end{parrafo}
\begin{lemma}\label{sobretruco}
Let \(({\mathcal F},(W,E))\) be a general basic object, let
$\{U_{\alpha}\}_{\alpha\in \Lambda}$ be the corresponding open
covering of $W$, and let
$(\widetilde{W}_{\alpha},(B_{\alpha},d_{\alpha}),\widetilde{E}_{\alpha})
$ be the collection of $d$-dimensional basic objects  associated
to \(({\mathcal F},(W,E))\). Then the functions
\begin{equation*}
\ord^{d}_{\alpha}: (F \cap
U_{\alpha}=)\Sing(B_{\alpha},d_{\alpha}) \to {\mathbb Q}
\end{equation*}
patch so as to define
\begin{equation*}
\ord^{d}_{\mathcal{F}}:F \to {\mathbb Q}.
\end{equation*}
\end{lemma}

The proof of this Lemma will be developed in Section
\ref{Hironaka}. It is an example of the principle of patching of
functions (\ref{patchingfunct}). Indeed it is the main example,
and the proof in Section \ref{Hironaka} will clarify why
projections were considered in \ref{defbo}, (iv) (and in
\ref{proyecciondebo}).

\begin{parrafo}
 Define \(({\mathcal F}_0,(W_0,E_0))\)
as before, by the covering $\{U_{\alpha}\}_{\alpha\in \Lambda}$
and basic objects
$(\widetilde{W}_{\alpha},(B_{\alpha},d_{\alpha}),\widetilde{E}_{\alpha})
$. Recall that a sequence of transformations
\begin{equation}
\label{trangbo1}
\begin{array}{ccccccc}
({{\mathcal F}}_0,(W_0,E_0)) & \longleftarrow & \ldots &
\longleftarrow & ({{\mathcal F}}_{r},(W_{r},E_{r})) &
\longleftarrow & ({{\mathcal F}}_{r+1},(W_{r+1},E_{r+1}))\\ \cup &
& & & \cup & & \cup\\ F_0 & & & & F_{r} & & F_{r+1},
\end{array}
\end{equation}
induces, for each index \(\alpha\),  a sequence of transformations
of basic objects
\begin{multline}
\label{transalpha1} (\widetilde{W}_{\alpha},
(B_{\alpha},d_{\alpha}), \widetilde{E}_{\alpha})\longleftarrow
((\widetilde{W}_{\alpha})_1, ((B_{\alpha})_1,d_{\alpha}),
(\widetilde{E}_{\alpha})_1)\longleftarrow\cdots
\\
\cdots((\widetilde{W}_{\alpha})_{r},
((B_{\alpha})_{r},d_{\alpha}),
(\widetilde{E}_{\alpha})_{r})\longleftarrow
((\widetilde{W}_{\alpha})_{r+1}, ((B_{\alpha})_{r+1},d_{\alpha}),
(\widetilde{E}_{\alpha})_{r+1});
\end{multline} and for each index $k$, set
\begin{equation}\label{jbarragbo}
(B_{\alpha})_{k} = I(H_{\alpha , 1})^{a_{\alpha ,1 }}\cdot
I(H_{\alpha , 2})^{a_{\alpha ,2}}\cdots I(H_{\alpha ,
k})^{a_{\alpha , k}}\cdot  \overline{(B_{\alpha})}_k
\end{equation}
as in \ref{expressbo1}.

\end{parrafo}

\begin{lemma}\label{gendetagbo}
Assume that sequence (\ref{trangbo1}) is such that, for each index
$0\leq k \leq r$:

1) The functions $$\word^{d}_{\alpha,k}: (F_k\cap (U_{\alpha})_k=)
\Sing((B_{\alpha})_{k},d_{\alpha})
 \to  {\mathbb Q} \mbox{ and } \overline{a}_{\alpha ,i}: \Sing(((B_{\alpha})_{k},d_{\alpha})) \to  {\mathbb Q}$$
patch to define functions $$\word_k:F_k\to {\mathbb
 Q}  \mbox{ and } \overline{a}_{i}: F_k \to  {\mathbb Q}.$$

2) If  \(({\mathcal F}_k,(W_k,E_k))\leftarrow ({\mathcal
F}_{k+1},(W_{k+1},E_{k+1})) \) is defined by a monoidal
transformation, assume that the center $Y_k$ is such that: $$Y_k
\subset \Max \word_k (\subset W_k).$$ Then, under assumptions 1)
and 2), the functions defined in terms the expression
(\ref{jbarragbo}) for the index $r+1$, namely the functions
$$\word^{d}_{\alpha,r+1}: \Sing((B_{\alpha})_{r+1},d_{\alpha})
 \to  {\mathbb Q} \mbox{ and } \overline{a}_{\alpha ,i}: \Sing(((B_{\alpha})_{r+1},d_{\alpha})) \to  {\mathbb Q}$$
( see \ref{defdeword}), patch, and define functions
$$\word_{r+1}:F_{r+1}\to {\mathbb
 Q}  \mbox{ and } \overline{a}_{i}: F_{r+1} \to  {\mathbb Q}.$$
\end{lemma}
\proof

A) Assume that $(\mathcal{F}_r,(W_r,E_r))\leftarrow
(\mathcal{F}_{r+1},(W_{r+1},E_{r+1}))$ is defined by a projection
$W_r \leftarrow W_{r+1}=W_r\times A^1$. Then for each $\alpha$,
the expression (\ref{jbarragbo}) for the index $r+1$ is the
pull-back of the expression for index $r$. In such case the
patching of functions with index $r+1$ follows from the case of
index $r$.

B) Assume that $(\mathcal{F}_r,(W_r,E_r))\leftarrow
(\mathcal{F}_{r+1},(W_{r+1},E_{r+1}))$ is defined by a center $Y_r
\subset F_r$ and let $H_{r+1}\subset W_{r+1}$ denote the
exceptional locus. Choose a point $x \in F_{r+1}(\subset W_{r+1})$
and assume that $x \in U_{\alpha_1, r+1}\cap U_{\alpha_2, r+1}$.
Consider the two expressions:
\begin{equation}
(B_{\alpha_j})_{k} = I(H_{\alpha_j , 1})^{a_{\alpha_j ,1 }}\cdot
I(H_{\alpha_j , 2})^{a_{\alpha_j ,2}}\cdots I(H_{\alpha_j ,
k})^{a_{\alpha_j , k}}\cdot  \overline{(B_{\alpha_j})}_k \mbox{
for } j=1,2.
\end{equation}
We want to prove that
$\word^{d}_{\alpha_1,r+1}(x)=\word^{d}_{\alpha_2,r+1}(x)$, and
that $ \overline{a}_{\alpha_1 ,i}(x)= \overline{a}_{\alpha_2
,i}(x)$ for each $H_i\in E_{r+1}$. If $x \in F_{r+1}-H_{r+1}$ then
$x$ can be identified with a point in $F_{r}$, and the equalities
follow by assumption. So assume that $x \in F_{r+1}\cap H_{r+1}$.
Note that
$$\ord_{\alpha_j}(x)= \sum \overline{a}_{\alpha_j
,i}(x)+\word_{\alpha_j,r+1}(x).$$ Lemma \ref{sobretruco} asserts
that $\ord_{\alpha_1}(x)=\ord_{\alpha_2}(x)$; and by assumption,
we also know that $\overline{a}_{\alpha_1
,i}(x)=\overline{a}_{\alpha_2 ,i}(x)$ for each hypersurface $H_i$
with index $i< r+1$. So it suffices to prove that
\begin{equation}\label{poiu}
\overline{a}_{\alpha_1 ,r+1}(x)=\overline{a}_{\alpha_2 ,r+1}(x).
\end{equation}
Note that $x \in H_{r+1}\subset W_{r+1}$ maps to a point $x'\in
Y_{r}\subset W_r$. Let $y$ be the generic point of the irreducible
component of $Y_r$ containing $x'$. To settle (\ref{poiu}), note
that
$$\overline{a}_{\alpha_j ,r+1}(x)=\sum_{H_t\in E_r \mbox{ and } Y_r
\subset H_t}\overline{a}_{\alpha_j ,r}(y)+ \word_{\alpha_j
,r}(y),$$ and, by assumption, all terms are independent of $j$.

\begin{remark}\label{deftgbo}

Lemma \ref{gendetagbo} was proved under some assumptions on
(\ref{trangbo1}) (for each index $0\leq k \leq r$). Note that such
assumptions hold for $r=0$. As in \ref{functk}, we see that
 $$\max \word_{k+1}\geq
\max \word_{k+1}\geq \ldots \geq \max \word_{k+1}$$ since
$Y_i\subset \Max \word_i\subset F_i$.  Set $\frac{b'_i}{b}=\max
\word_{i}$, so $\frac{b'_0}{b}\geq \frac{b'_1}{b}\geq\ldots \geq
\frac{b'_{k-1}}{b}\geq \frac{b'_k}{b}$, and let $k_0$ be the
smallest index such that $\frac{b'_{k_0}}{b}=\frac{b'_{k_0+1}}{b}=
\ldots =\frac{b'_{k}}{b}.$

For each index $ k_0 \leq j \leq k$ define a partition on the set
of hypersurfaces in $E_j$, say $ E_j=E_j^- \cup E_j^+$, where
$E_j^-=\{H_1,\ldots, H_r, \ldots H_{r+k_0} \}$ and
$E_j^+=\{H_{r+k_0+1},\ldots H_{r+j} \}$.

For $j=k_0$ $E_{k_0}=E_{k_0}^-$, and for $j>k_0$, $E_{j}^-$ are
the strict transforms of hypersurfaces in $E_{k_0}$. Order
${\mathbb Q}\times {\mathbb N}$ lexicographically, and set
\begin{equation*}
\begin{array}{rccl}
t^d_k: & F_k & \to & {\mathbb Q}\times {\mathbb N}\\ & x &
\longrightarrow &

(\word^{d}_k(x),n^d_k(x))
\end{array}
\end{equation*}
\begin{equation*}
       n_{k}^d(\xi)=\left\{
\begin{array}{lll}
     \#\{H\in E_{k}\mid \xi\in H\} & {\rm if} &

\word_{k}^d(\xi)<\frac{b'_k}{b} \\
     \#\{H\in E_{k}^{-}\mid \xi\in H\}
& {\rm if} &
     \word_{k}^d(\xi)=\frac{b'_k}{b}.
      \end{array}
\right.
\end{equation*}
One can check, as for the case of basic objects, that this
function is upper-semi-continuous.
\end{remark}

 \vspace{5mm} We now extend \ref{prop7.9} to the
setting of general basic objects.

\begin{proposition}\label{prop119}  Consider a sequence of transformations
\begin{equation}
\label{trangbo11}
\begin{array}{ccccccc}
({{\mathcal F}}_0,(W_0,E_0)) & \longleftarrow & \ldots &
\longleftarrow & ({{\mathcal F}}_{k_0-1},(W_{k_0-1},E_{k_0-1})) &
\longleftarrow & ({{\mathcal F}}_{k_0},(W_{k_0},E_{k_0}))\\ \cup &
& & & \cup & & \cup\\ F_0 & & & & F_{k_0-1} & & F_{k_0},
\end{array}
\end{equation}
where each transformation is either a projection, or
transformation with centers $Y_i\subset \Max t_i$; and assume that
$max t_{i-1} \geq max t_i$. Assume, in addition, that $\max
t^d_{k_0-1}> \max t^d_{k_0}$, or that $k_0=0$. Then, there is a
simple general basic object $({{\mathcal
G}}_{k_0},(W_{k_0},E_{k_0}'))$, such that any resolution
\begin{equation}
\label{tratkkk1}\tiny
\begin{array}{cccccc}
({{\mathcal G}}_{k_0},(W_{k_0},E_{k_0}'))& \longleftarrow &
({{\mathcal G}}_{k_0+1},(W_{k_0+1},E_{k_0+1}')) & \longleftarrow &
\ldots & ({{\mathcal
G}}_{k_0+m},(W_{k_0+m},E_{k_0+m}')))\\G_{k_0}&&G_{k_0+1}&&&G_{k_0+m}=\emptyset
\end{array}
\end{equation}
induces a sequence of transformations
\begin{equation}
\label{trangbo11kk}\tiny
\begin{array}{cccccccc}
({{\mathcal F}}_0,(W_0,E_0)) & \longleftarrow & \ldots &
\longleftarrow & ({{\mathcal F}}_{k_0},(W_{k_0},E_{k_0})) &
\longleftarrow & \cdots &({{\mathcal
F}}_{k_0+m},(W_{k_0+m},E_{k_0+m}))\\ \cup & & & & \cup & & &
\cup\\ F_0 & & & & F_{k_0}& & & F_{k_0+m}.
\end{array}
\end{equation}

Furthermore, this sequence has the following two properties:

{\bf P1)} $\max t^d_{k_0}=\ldots=\max t^d_{k_0+m-1}$; and
$F_{k_0+m}=\emptyset$ or $\max t^d_{k_0+m-1}>\max t^d_{k_0+m}$.

{\bf P2)} $\Max t^d_{k_0+j}=G_{k_0+j}$ for $ 0\leq j \leq m-1$.

\end{proposition}
\proof

Note that the properties in Proposition \ref{prop7.9} assert that
$({{\mathcal G}}_{k_0},(W_{k_0},E_{k_0}'))$ is indeed a general
basic object.

\begin{proposition}\label{prop119000}
Assume that (\ref{trangbo11}) is such that $\max t_{k_0-1}>\max
t_{k_0}$ and that $\max \word_{k_0}=0$. Then, there are
upper-semi-continuous functions $h_i:F_{k_0+i} \to \Gamma$, and a
resolution
$$({{\mathcal F}}_{k_0},(W_{k_0},E_{k_0}))\leftarrow ({{\mathcal
F}}_{k_0+1},(W_{k_0+1},E_{k_0+1}))\leftarrow \cdots \leftarrow
({{\mathcal F}}_{k_0+m},(W_{k_0+m},E_{k_0+m})).$$  The resolution
defined by blowing up successively on $\Max h_i$.
\end{proposition}
\proof This is an extension of \ref{prop7.9casoword0} to the case
of general basic object. The fact that $h_i$ are well defined as
functions on $F_{k_0+i}$ follows from Lemma \ref{gendetagbo}.

\begin{theorem}
\label{Theoremd} {\em (Theorem (d))} Fix a positive integer $d$.
There is a totally ordered set $I^d$, and for each \( d
\)-dimensional general basic object \((\mathcal{F}_0,(W_0,E_0))\),
a function $f_{\mathcal{F}_0}:F_0 \to I^d$. The functions defined
so that:

i)$f_{\mathcal{F}_0}$ is upper-semi-continuous, and $\Max
f_{\mathcal{F}_0}\subset F_0$ is a smooth permissible center for
\((\mathcal{F}_0,(W_0,E_0))\).

ii) For each \((\mathcal{F}_0,(W_0,E_0))\), there is a resolution
\(R_{\mathcal{F}_0}\) ( Def \ref{resgbo}), say
\begin{equation}
\label{trangbo2}
\begin{array}{ccccccc}
({{\mathcal F}}_0,(W_0,E_0)) & \longleftarrow & \ldots &
\longleftarrow & ({{\mathcal F}}_{r},(W_{r},E_{r})) &
\longleftarrow & ({{\mathcal F}}_{r+1},(W_{r+1},E_{r+1}))\\ \cup &
& & & \cup & & \cup\\ F_0 & & & & F_{r} & & F_{r+1}=\emptyset,
\end{array}
\end{equation}
obtained by blowing up successively at $\Max f_{\mathcal{F}_i}$ (
$f_{\mathcal{F}_i}:F_i \to I^d$).
\end{theorem}

\proof The proof is based on inductive argument, so we first show
why Theorem \ref{Theoremd} holds for $0-$dimensional general basic
objects: Note that in such case, each
$(\widetilde{W}_{\alpha},(B_{\alpha},d_{\alpha}),E_{\alpha})$ is
zero dimensional, so we can assume that each
$\widetilde{W}_{\alpha}$ is a point, and hence, each $B_{\alpha}$
is a non-zero ideal in a field. Therefore,
$\Sing(B_{\alpha},d_{\alpha})=\emptyset$, and hence,
$F_0=\emptyset$. Here we can take $I^0$ to be a point; it plays no
role in any case.

Set \(T^d=\{\infty\}\sqcup ({\mathbb Q}\times {\mathbb Z}) \sqcup
\Gamma\) with $ \Gamma $ as in \ref{DefGamma}. This disjoint union
is totally ordered by setting \(\infty\) as the biggest element,
and \(\alpha < \beta \) if \(\beta \in ({\mathbb Q}\times {\mathbb
Z})\) and \(\alpha \in \Gamma \). We now set \( I^{d}=T^d \times
I^{d-1} \) ordered lexicographically. In our proof,
upper-semi-continuous functions \( f_i^d: F_i \to I_{d}\) will be
defined with the property stated in the theorem; namely, that a
resolution (\ref{trangbo2}) will be achieved by taking successive
monoidal transformations with centers $Y_i=\underline{Max} f_i^d
(\subset F_i)$.

For the index $i=0$ we know that the locally defined functions
$\ord^{d}_{\alpha}: \Sing(B_{\alpha},d_{\alpha}) \to {\mathbb Q}$
 and $ n^{d}_{\alpha}: \Sing(B_{\alpha},d_{\alpha}) \to {\mathbb
Z}$ patch so as to define functions
\begin{equation*}
\ord^{d}_{\mathcal{F}}:F \to {\mathbb Q} \qquad\text{and}\qquad
n^{d}_{\mathcal{F}}: F \to{\mathbb Z},
\end{equation*}
and also an upper-semi-continuous function $$t_0=(\ord^{d},n^d):
F_0 \to {\mathbb Q} \times {\mathbb Z} $$ as in \ref{deftgbo}
(recall that $\word_0=\ord_0$).

We attach to the value $max \ t^d_0$ the simple general basic
object $({{\mathcal G}}_{0},(W_{0},E_{0}'))$ as in \ref{prop119},
so that $G_0=\Max t^d_0$. We now define a resolution of this
simple general basic object:
\begin{equation}
\begin{array}{cccccc}\label{jjjjj}
({{\mathcal G}}_{0},(W_{0},E_{0}'))& \longleftarrow & ({{\mathcal
G}}_{1},(W_{1},E_{1}')) & \longleftarrow & \ldots & ({{\mathcal
G}}_{m},(W_{m},E_{m}'))\\G_{0}&&G_{1}&&&G_{m}=\emptyset
\end{array}
\end{equation}
and, in order to define this resolution, we first apply the
transformation with center $R(1)(G_0)$, if not empty, so as to
obtain a $d-1$ general basic object (see \ref{onsgbo}). We then
proceed to define the resolution (\ref{jjjjj}) by induction (i.e.
by blowing up successively at $\Max f_i^{d-1}\subset G_i$). This
defines
\begin{equation}
\begin{array}{cccccc}\label{jjjjjj}
({{\mathcal F}}_{0},(W_{0},E_{0}))& \longleftarrow & ({{\mathcal
F}}_{1},(W_{1},E_{1})) & \longleftarrow & \ldots & ({{\mathcal
F}}_{m},(W_{m},E_{m})))\\F_{0}&&F_{1}&&&F_{m}
\end{array}
\end{equation}
and functions $t^{d}_i:F_i\to{\mathbb Q} \times {\mathbb Z}$ for
$i=0,\dots, m$. Furthermore, $G_i=\Max t^{d}_i$ for $i=0,\dots,
m-1$, and
$$\max t^{d}_0=\max t^{d}_1= \dots =\max t^{d}_{m-1}\mbox{; and either }
F_m=\emptyset \mbox{ or } \max t^{d}_{m-1}>\max t^{d}_{m}.$$
\begin{parrafo}\label{ultimo}
We now define functions $f^{d}_i$, however, for the time being,
only at points of $G_i =\Max t^{d}_i(\subset F_i)$:

i) $ f_0^d(x)=(\max t^d,\infty)\in I^{d}=T^d \times I^{d-1}\mbox{
if } x\in R(1)(G_0)$;

ii) $ f_0^d(x)=(\max t^d,f_0^{d-1}(x))\in I^{d}=T^d \times
I^{d-1}\mbox{ if } x\notin R(1)(G_0)$;

iii) $f_i^d(x)=(\max t^d,f_0^{d-1}(x))\in I^{d}=T^d \times
I^{d-1}\mbox{ for } i=1, \dots , m-1.$

Note now that sequence (\ref{jjjjjj}), induced by (\ref{jjjjj}),
is also defined by blowing up successively at $\Max f^d_i$, by the
way such functions are defined.

The condition $ \max t^{d}_{m-1}>\max t^{d}_{m}$ asserts that
$\max \word_{m-1} \geq \max \word_m$. In fact $ \max
t^{d}_{i}=(\max \word_i,a_i)$ where $a_i\leq \mbox{ dim } W_i=
\mbox{ dim } W$. So either $\max \word_{m-1}>\max \word_{m}$; or
$\max \word_{m-1}=\max \word_{m}$ and $ 0\leq a_m<a_{m-1}\leq dim
W$. If $F_m\neq \emptyset$ in (\ref{jjjjjj}) we define
\begin{equation}
\begin{array}{cccccc}\label{marisma}
({{\mathcal F}}_{m},(W_{m},E_{m}))& \longleftarrow & ({{\mathcal
F}}_{m+1},(W_{m+1},E_{m+1})) & \longleftarrow & \ldots &
({{\mathcal F}}_{m_1},(W_{m_1},E_{m_1}))\\F_{m}&&F_{m+1}&&&F_{m_1}
\end{array}
\end{equation}
distinguishing two different cases:

A) If $\max \word_{m}>0$, argue as above, and define
(\ref{marisma}) as the resolution of a simple general basic object
$({{\mathcal G}}_{m},(W_{m},E_{m}))$ (\ref{prop119}). Finally
define functions $f_{m+i}$ on $G_{m+i}$, as done before, so that
(\ref{marisma}) is obtained by blowing up $\Max f_{m+i}^d$.

B) If $\max \word_{m}=0$ set (\ref{marisma}) as in
\ref{prop119000}. Hence, a resolution $({{\mathcal
F}}_{m},(W_{m},E_{m}))$, and thus of $({{\mathcal
F}}_{0},(W_{0},E_{0})))$, is obtained by blowing up at $\Max
h_{m+i}\subset F_{m+i}$. Finally set: $$f_{m+i}: F_{m+i}\to I^d \
, \ f^d_{m+i}(x)=(h_{m+i},\infty)\in I^{d}=T^d \times I^{d-1};$$
and note that (\ref{marisma}) is defined by blowing up at $\Max
f^{d}_{m+i}\subset F_{m+i}$.

Finally note that case A) will occur only finitely many time, and
either a resolution is achieved or case B) occurs, and thus we
always achieve a resolution. In fact, $\max \word$ can take only
finitely many values, and the second coordinate of $\max t^d$ is a
positive integer $\leq dim \ W$. This finiteness was discussed in
\ref{com101} for the case of one basic object. Note that a general
basic object can be covered finitely many affine sets.

In this way we have defined functions $f^d_i$, and a resolution
(\ref{trangbo2}), obtained by blowing up at $\Max f^d_i$. However
our functions $f^d_i$ are only defined in $G_i (\subset F_i)$ (set
formally $G_i=F_i$ in Case A)). So we now extend the definition of
the functions to all $F_i$, and we must do so in a way that does
not modify the sets $\Max f_i^d$ already considered.

Since $\Max  f_i^d \subset G_i \subset F_i$, a point $ x \in
F_i-G_i$ can be identified with a point, say $x \in F_{i+1}$.
Furthermore, since (\ref{trangbo2}) is a resolution, there is
smallest index $i_0 > i$ such that $x$  can be identified with a
point, say again $x \in G_{i_0} (\subset F_{i_0})$. Define
$$f_i(x)=f_{i_0}(x).$$ Note that $t^d_i(x)=t^d_{i_0}(x)$ (if
$\word_i(x)>0$), that $h_i(x)= h_{i_0}(x)$ (if $\word_i(x)=0$);
and that an open neighborhood of $x$ in $F_i$ can be identified
with a neighborhood of $x$ in $G_{i_0}$. Finally argue
coordinate-wise (as in \ref{functk}) to show that the extended
functions $f_i^d: F_i \to I^d=T^d\times I^{d-1}$ are in fact
upper-semi-continuous.

The compatibility with open restrictions of the functions $t_i^d$
and $h_i$, and also that of $f_i^{d-1}$ (by induction), insure
that the same property holds for the functions $f_i^d$ (see
\ref{compopen}).
\end{parrafo}

\begin{parrafo}\label{resfunct2} {\em On Resolution functions and Proof of \ref{resfunct}}.
Recall that basic objects are, in particular, general basic
objects (\ref{notaopen}), so Theorem \ref{Theoremd} provides, for
each dimension $d$, resolution functions as in \ref{resfunct}.
Here $$I^{d}=T^d \times I^{d-1}=T^d \times T^{d-1} \times \cdots
\times
 T^0$$ and $f_i^d(x)$ can be expressed with $d+1$ coordinates. For
 instance case i) in \ref{ultimo} is:

i) $ f_0^d(x)=(\max t^d,\infty,\infty, \dots ,\infty)\in I^{d}
\mbox{ if } x\in R(1)(G_0)$.

If $W$ is smooth of dimension $d$, and $X\subset W$ is a smooth
hypersurface, then the basic object $(W,(J,1),E=\emptyset)$
defines a $d$-dimensional general basic object. In this case
$t^d_0(x)=(1,0)$, and $ f_0^d(x)=((1,0),\infty,\infty, \dots
,\infty) $ for any $x\in \Sing(J,1)$.

If $X\subset W$ is smooth of codimension two, then $
f_0^d(x)=((1,0),(1,0),\infty, \dots,\infty) $ for any $x\in
\Sing(J,1)$.

 If $X\subset W$ is smooth of codimension $r$, then $
f_0^d(x)=R=((1,0),(1,0),\dots,(1,0),\infty, \dots,\infty) $ (r
copies of $(1,0)$) for any $x\in \Sing(J,1)$; and if $X\subset W$
is reduced, pure dimensional and of codimension $r$, then
$((1,0),(1,0),\dots,(1,0),\infty, \dots ,\infty)$ (r copies of
$(1,0)$) is the value $R$ in property P4) of \ref{resfunct}.

\end{parrafo}

\section{On Hironaka's trick and proof of Lemma \ref{sobretruco}}

 \label{Hironaka}
 The purpose of this Section is to prove
Lemma \ref{sobretruco} which states that the function $\ord$,
introduced in \ref{eqord} for basic objects, can be naturally
defined in the setting of general basic objects.

Let \( (\mathcal{F},(W,E)) \) be an \( d \)-dimensional
 general basic object, and set an open covering  \(
 \{U_{\alpha}\}_{\alpha\in
\Lambda} \)
    of $W$ as in Definition
\ref{defbo}.

Recall that $( \mathcal{F},(W,E))$ defines a closed set $F
(\subset W)$, and that for each index $\alpha$ there is a closed
smooth d-dimensional subscheme $\widetilde{W_{\alpha}} \subset
U_{\alpha}$, and a basic object $(\widetilde{W_{\alpha}},
(B_{\alpha}, d_{\alpha}), \widetilde{E_{\alpha}})$ such that
$$F\cap U_{\alpha}=\Sing(B_{\alpha}, d_{\alpha}).$$

Assume that a point $x\in F$ appears in two such charts, namely
$x\in F\cap U_{\alpha}\cap U_{\beta}$. In order to simplify
notation set $$(\widetilde{W_{\alpha}}, (B_{\alpha}, d_{\alpha}),
\widetilde{E_{\alpha}})= (W',(B',d'),E')
$$ and

$$(\widetilde{W_{\beta}}, (B_{\beta}, d_{\beta}),
\widetilde{E_{\beta}})= (W'',(B'',d''),E'').$$

So \( x\in\Sing(B',d')=\Sing(B'',d'') \) and the claim in Lemma
\ref{sobretruco} is that:
 \begin{equation*}
\frac{\Fnu_{B'}(x)}{d'}=\frac{\Fnu_{B''}(x)}{d''}.
 \end{equation*}
(notation as is \ref{eqord}).

 \begin{proof}
 Set \( \omega'=\Fnu_{B'}(x) \) and \(
\omega''=\Fnu_{B''}(x) \).
 We shall prove the Lemma by constructing
infinitely many sequences of transformations of general basic
 objects. A sequence
 \begin{equation}
(\mathcal{F},(W,E))\stackrel{\Pi_{0}}{\longleftarrow}
(\mathcal{F}_{0},(W_{0},E_{0}))\stackrel{\Pi_{1}}{\longleftarrow}
(\mathcal{F}_{1},(W_{1},E_{1}))\stackrel{\Pi_{2}}{\longleftarrow}
\cdots\stackrel{\Pi_{k}}{\longleftarrow}
(\mathcal{F}_{k},(W_{k},E_{k}))
        \label{SeqGBOHiro}
 \end{equation}
of transformations of general basic objects defines sequences of
 transformations of basic
objects, say:
 \begin{multline}
(W',(B',d'),E')\stackrel{\Pi'_{0}}{\longleftarrow}
(W'_{0},(B'_{0},d'),E'_{0})
         \stackrel{\Pi'_{1}}{\longleftarrow}
(W'_{1},(B'_{1},d'),E'_{1})
         \stackrel{\Pi'_{2}}{\longleftarrow}
\cdots \\
         \cdots
        \stackrel{\Pi'_{k}}{\longleftarrow}
(W'_{k},(B'_{k},d'),E'_{k}),
         \label{SeqBasHiro1}
 \end{multline}

and
 \begin{multline}
(W'',(B'',d''),E'')\stackrel{\Pi''_{0}}{\longleftarrow}
(W''_{0},(B''_{0},d''),E''_{0})
\stackrel{\Pi''_{1}}{\longleftarrow}
(W''_{1},(B''_{1},d''),E''_{1})
\stackrel{\Pi''_{2}}{\longleftarrow} \cdots \\
        \cdots
\stackrel{\Pi''_{k}}{\longleftarrow}
(W''_{k},(B''_{k},d''),E''_{k}). \label{SeqBasHiro2}
 \end{multline}
 We take the first transformation \( \Pi_{0}
\) of (\ref{SeqGBOHiro}) to be a
 projection (as in \ref{defbo} iv)), so the first
transformations of
 (\ref{SeqBasHiro1}) and (\ref{SeqBasHiro2}) are
projections too.
 All the other transformation will be permissible
transformations (as in
 (as in \ref{defbo} iii))). For each index \( k >0\), sequence
 (\ref{SeqGBOHiro}) will be defined as follows:
 \begin{enumerate}

\item Identify \( L_{0}=\Pi_{0}^{-1}(x) \) with \(
\mathbb{A}^{1}_k \) and set \( x_{0} = 0 \in L_0 \).
         Note that \(
L_{0} \subset F_{0} \), the singular locus of \(
(\mathcal{F}_{0},(W_{0},E_{0})) \).
        \item Given an index \( s \geq
0 \), a line \( L_s \subset F_s \)
         and a point \( x_s \in L_s \),
define the transformation \(
         \Pi_{s+1} \) with center \( x_s \).
Now set:
         \begin{description}
      \item[i] \(L_{s+1}\subset
F_{s+1}\) as the strict transform
      of \( L_s \);
       \item[ii] \(
H_{s+1}(\in E_{s+1}) \) as the exceptional
       locus of \( \Pi_{s+1}
\);
      \item[iii] \( x_{s+1}=H_{s+1}\cap L_{s+1} \).

\end{description}
 \end{enumerate}
 In this way (1) together with (2)
provide a rule to construct
  a sequence (\ref{SeqGBOHiro}) of
length \( s \), for any \( s
 \geq 1 \). In this sequence \( L_{s}\subset F_{s} \) for any \( s \), so in
 particular \(
x_{s}\in F_{s} \), and by assumption:
 \begin{equation*}
x_{s}\in\Sing(B'_{s},d')=\Sing(B''_{s},d'')
         \qquad \forall s\geq
0.
 \end{equation*}
 Locally at $x_s$ there are expressions, as in (\ref{expressbo1}),
 say:

\begin{equation} \label{ExprHiro}
       (B'_{s})_{ x_{s}}=
\mathcal{I}(H'_{s})_{ x_{s}}^{a'_{s}}(\overline{B'}_{s})_{ x_{s}}
\qquad
        (B''_{s})_{ x_{s}}=
        \mathcal{I}(H''_{s})_{
x_{s}}^{a''_{s}}(\overline{B''}_{s})_{ x_{s}}.
 \end{equation}
 Note that
here  \( H'_{s}=H_{s}\cap W'_{s} \) and \(
 H''_{s}=H_{s}\cap W''_{s} \).
On may check, by induction on \( s
 \), that
 \begin{equation*}
a'_s = s(\omega' - d') \qquad a''_{s}=s(\omega''-d'').
 \end{equation*}

Since only the first term of this sequence is a projection, for \(
 s \geq
1 \), \( \dim(W'_{s})=\dim(W''_{s})=d+1 \). It follows that
\begin{equation*}
         \begin{split}
            \dim(F_s\cap
H_s)=d
      & \Leftrightarrow a'_{s}=s(\omega'-d')\geq d'  \\
& \Leftrightarrow a''_{s}=s(\omega''-d'')\geq d''.
         \end{split}
\end{equation*}
 Note that \( \dim{H'_{s}}=\dim{H''_{s}}=d \), so if \(
\dim(F_s\cap H_s)=d \) then \( F_{s}\cap H_{s}=H'_{s}=H''_{s} \).

Furthermore if \( \dim(F_s\cap H_s)=d \), then \( F_s\cap H_s \)
 is a
permissible center for the general basic object.
 In such case, set \( F_s
\cap H_s \) as a center of a
 transformation \( \Pi_{s+1} \). It turns out
that in (\ref{ExprHiro}),
 \begin{equation*}
a'_{s}=s(\omega'-d')-d', \qquad a''_{s}=s(\omega''-d'')-d''.
\end{equation*}
 Fix the index \( s \) and set, if possible, the center of
transformations \( \Pi_{s+j} \) to be \( F_{s+j}\cap H_{s+j} \),
for \( j\geq 0 \). Note that
\begin{equation*}
         \begin{split}
\dim(F_{s+j}\cap H_{s+j})=d
      & \Leftrightarrow
a'_{s+j}=s(\omega'-d')-jd'\geq d'  \\
            & \Leftrightarrow
a''_{s+j}=s(\omega''-d'')-jd''\geq d''.
\end{split}
\end{equation*}
And we conclude that
 \begin{equation*}
\begin{split}
             \dim(F_{s+j}\cap H_{s+j})=d\ \text{(in which
case is a permissible
             center)}
       & \Leftrightarrow j\leq
\ell'_{s}  \\
             & \Leftrightarrow j\leq \ell''_{s},
\end{split}
 \end{equation*}
where
\begin{equation*}
\ell'_{s}=\left\lfloor\frac{s(\omega'-d')}{d'}\right\rfloor
       \qquad
\ell''_{s}=\left\lfloor\frac{s(\omega''-d'')}{d''}\right\rfloor
\end{equation*}
 and \( \lfloor\cdot\rfloor \) denotes the integer part.

Finally note that
 \begin{gather*}
\frac{\Fnu_{B'}(x)}{d'}=\frac{w'}{d'}=
\lim_{s\to\infty}\frac{1}{s}\ell'_{s}+1,  \\
\frac{\Fnu_{B''}(x)}{d''}=\frac{w''}{d''}=
\lim_{s\to\infty}\frac{1}{s}\ell''_{s}+1.
 \end{gather*}
This expresses the rational numbers \(
 \dfrac{\Fnu_{B'}(x)}{d'} \) and \(
\dfrac{\Fnu_{B''}(x)}{d''} \)
  in terms of permissible sequences of the
general basic object  \(
 (\mathcal{F},(W,E)) \). Hence \( \dfrac{\Fnu_{B'}(x)}{d'}=\dfrac{\Fnu_{B''}(x)}{d''} \).
 \end{proof}

\end{document}